\documentclass[10pt]{article}

\usepackage{tikz}
\usepackage{adjustbox}

\usepackage[round]{natbib}

\usepackage{amsmath, amsfonts, amssymb, amsthm}

\usepackage{mathrsfs}

\usepackage{soul}

\usepackage{bm}
\usepackage{dsfont}

\usepackage{booktabs} 

\usepackage[colorlinks=true, linkcolor=red, urlcolor=blue, citecolor=blue]{hyperref}

\usepackage{geometry}
\usepackage{graphicx}
\usepackage{color}
\usepackage{nicefrac}

\usepackage{scrextend}

\usepackage{setspace}                                                   
\usepackage{titlesec}
\titleformat{\section}[hang]{\large\scshape}{\thesection.}{1em}{}

\newtheoremstyle{mytheoremstyle} 
    {0.3cm}                      
    {0cm}                        
    {\itshape}                   
    {}                           
    {\scshape}                   
    {: }                          
    {0em}                       
    {}  

\theoremstyle{mytheoremstyle}
\newtheorem{Theorem}{Theorem}
\newtheorem{Lemma}{Lemma}
\newtheorem{Corollary}{Corollary}
\newtheorem{Proposition}{Proposition}

\newtheoremstyle{myExampleRemarkstyle} 
    {0.3cm}                    
    {0cm}                           
    {\itshape}                   
    {}                           
    {\scshape}                   
    {: }                          
    {0em}                       
    {}  

\theoremstyle{myExampleRemarkstyle}
 
\newtheorem{Remark}{Remark}
\newtheorem{Assumption}{Assumption}

\renewcommand{\theAssumption}{\Alph{Assumption}}

\newtheoremstyle{simuStyle}
{0.3cm} 
{0cm} 
{} 
{} 
{\bfseries} 
{.} 
{0em} 
{} 

\theoremstyle{simuStyle}

\newtheoremstyle{stratStyle}
{0.3cm} 
{0cm} 
{} 
{} 
{\scshape} 
{: } 
{0em} 
{} 

\theoremstyle{stratStyle}

\DeclareSymbolFont{lettersA}{U}{txmia}{m}{it}
\DeclareMathSymbol{\real}{\mathord}{lettersA}{"92}
\DeclareMathSymbol{\field}{\mathord}{lettersA}{"83}

\def\real{{\rm I\!R}}

\DeclareMathOperator*{\argmin}{argmin}

\def\0{{\bf 0}}

\def\btheta{{\bm{\theta}}}

\def\bbeta{{\bm{\beta}}}

\def\bpi{{\bm{\pi}}}
\def\hbpi{\hat{{\bm{\pi}}}}

\def\bOmega{{\bm{\Omega}}}

\DeclareMathOperator*{\var}{var}

\DeclareMathOperator*{\argzero}{argzero}

\definecolor{mypurple}{RGB}{0,0,0}
\definecolor{myblue}{RGB}{0,87,120}
\definecolor{aqua}{RGB}{0,0,0}
\definecolor{myorange}{RGB}{0,0,0}
\definecolor{mygrey}{RGB}{255,255,255}

\definecolor{col12}{RGB}{175,35,36}
\definecolor{col14}{RGB}{0,104,71}
\definecolor{col15}{RGB}{43,99,173}

\usetikzlibrary{shapes.geometric}

\def\ph{\hat{\bm{\pi}}}
\def\bt{\bm{\theta}}
\def\bto{\bt_0}

\def\bD{\bm{\Delta}}
\def\hbt{\hat{\bt}_{(j,n,H)}}
\def\cn{\mathbf{c}(n)}
\def\bw{\bm{\omega}}
\def\bwo{\bw_0}

\def\Ln{\mathbf{L}(n)}

\def\0{\mathbf{0}}
\def\bI{\mathbf{I}}

\def\bb{\bm{\beta}}
\def\hbb{\hat{\bb}}
\def\hbbr{\hbb^{\text{ridge}}_{\lambda}}
\def\hbbo{\hbb^{\text{OLS}}}

\def\bX{\mathbf{X}}
\newcommand{\hp}[2]{\ph\left({#1},n,{#2}\right)}
\newcommand{\hpH}[1]{\frac{1}{H}\sum_{h=1}^H\ph\left({#1},n,\bw_{h+jH}\right)}
\renewcommand{\a}[1]{\mathbf{a}\left({#1}\right)}
\renewcommand{\r}[1]{\mathbf{r}\left({#1},n\right)}
\renewcommand{\v}[2]{\mathbf{v}\left({#1},n,{#2}\right)}
\newcommand{\vH}[1]{\frac{1}{H}\sum_{h=1}^H\mathbf{v}\left({#1},n,\bw_{h+jH}\right)}

\def\boxit#1{\vbox{\hrule\hbox{\vrule\kern3pt
          \vbox{\kern3pt#1\kern3pt}\kern3pt\vrule}\hrule}}

\definecolor{pinegreen}{rgb}{0.0, 0.47, 0.44}

\begin{document} 
	
	\begin{center}
		{\huge \textsc{On the Properties of Simulation-based\\ Estimators in High Dimensions}}\\
		\vspace{0.75cm}
		{\large \textsc{St\'ephane Guerrier}$^\S$, \textsc{Mucyo Karemera}$^\S$,\\ \textsc{Samuel Orso}$^\P$ \& \textsc{Maria-Pia~Victoria-Feser}$^\P$ }\\
		\vspace{0.5cm}
		{$^\S\,$Pennsylvania State University; $^\P\,$Research Center for Statistics, GSEM, University of Geneva}
	\end{center}
\vspace{0.25cm}

	\begin{addmargin}[3em]{3em}
		\footnotesize
		\textsc{Abstract:} 
Considering the increasing size of available data, the need for statistical methods that control the finite sample bias is growing. This is mainly due to the frequent settings where the number of variables is large and allowed to increase with the sample size bringing standard inferential procedures to incur significant loss in terms of performance. Moreover, the complexity of statistical models is also increasing thereby entailing important computational challenges in constructing new estimators or in implementing classical ones. A trade-off between numerical complexity (e.g. approximations of the likelihood function) and statistical properties (e.g. reduced finite sample bias) is often accepted. However, numerically efficient estimators that are altogether unbiased, consistent and asymptotically normal in high dimensional problems would generally be ideal. In this paper, we set a general framework from which such estimators can easily be derived for wide classes of models. This framework is based on the concepts that underlie simulation-based estimation methods such as indirect inference. The approach allows various extensions compared to previous results as it is adapted to possibly inconsistent estimators and is applicable to discrete models (e.g. logistic regression) and/or models with a large number of parameters (compared to the sample size). We consider an algorithm, namely the Iterative Bootstrap (IB), to efficiently compute simulation-based estimators by showing its convergence properties. Within this framework we also prove the properties of simulation-based estimators, more specifically the (finite sample) unbiasedness, consistency and asymptotic normality when the number of parameters is allowed to increase with the sample size. Therefore, an important implication of the proposed approach is that it allows to obtain unbiased estimators in finite samples. Finally, we study this approach when applied to three common models, namely logistic regression, negative binomial regression and regularized regression (lasso). The simulations not only confirm the theoretical results but, in many instances, show a gain in finite sample efficiency (mean squared error) compared to standard estimators and bias-reduced ones, leading to the conclusion that, in finite sample sizes, removing the bias does not necessarily increase the variance.
\vspace{0.25cm}


\vspace{0.25cm}

\noindent \textsc{Acknowledgement:} The authors are grateful to Guillaume Blanc, Roberto Molinari and Yuming Zhang for their comments that greatly helped to improve the content and the presentation of the paper.

	\end{addmargin}


\section{Introduction}
\label{seq:intro2}

In modern statistical analysis, an important challenge lies in the control of finite-sample bias of classical estimators. For example, the Maximum Likelihood Estimator (MLE), which is asymptotically unbiased under some regularity conditions, has finite-sample bias that can result in a significant inferential loss (see e.g. the review of \citealp{kosmidis2014bias}). This problem is typically magnified in situations where the number of variables $p$ is large and possibly allowed to increase with the sample size $n$.  For example, \cite{CaSu:18} show that the MLE for the logistic regression model can be severely biased in the case where $n$ and $p$ become increasingly large (given a fixed ratio). The idea of allowing the dimension $p$ of the data to increase with the sample size $n$ is motivated by the fact that many new technologies are now producing extremely large datasets with huge numbers of features. Moreover, along with this rapid growth in data size, model complexity is also increasing, thus creating important computational challenges in constructing new estimators or in implementing standard ones such as the MLE. To compute these estimators, approximate methods such as pseudo-likelihood functions or approximated estimating equations are increasingly used in practice as the resulting estimators are typically numerically easier to implement, although they might lead to biased and/or inconsistent estimators. Hence, estimators for parametric models can be biased for several reasons and, while various finite sample or asymptotic bias correction methods have been proposed, a widely applicable approach that corrects finite sample and asymptotic bias, in possibly high dimensional settings, is not available. This paper aims at making an additional step in this direction.

Finite-sample bias reduction methods have recently received substantive attention. One strategy consists in deriving an estimate of the bias which is then added to the estimator. Bias estimation is often based on simulation methods such as the jackknife \citep{Efro:82} or the bootstrap \citep{Efro:79}, or can be approximated using asymptotic expansions (see e.g. \citealp{CoVa:97}, \citealp{CoTU:08} and the references therein\footnote{Such an approach was originated by \cite{CoSn:68} and further developed by \cite{efron1975,Scha:83,GaPeTh:85,CoTsWe:86,CoMcCu:91,LiBr:96}.}). An alternative approach is to correct the bias of an estimator by modifying its associated estimating equations. For example, \cite{Firt:93} provides an adjustment of the MLE score function and this approach has been successively adapted and extended to Generalized Linear Models (GLM), among others, by \cite{MeMa:95,BuMaGr:02,KoFi:09,KoFi:11,Kosm:14}\footnote{Firth's correction has also been used to correct the finite sample bias of the MLE of other models; see e.g. \cite{PeKeGa:98,BuLeLe:07,Kosm:17}.}. In a similar fashion, \cite{KPSaSa:17} propose a median bias reduction adjustment for the MLE \citep[see also][]{HiTsMe:89,KyKoSa:18}.  Moreover, when the finite sample bias results from the presence of nuisance parameters (i.e. when the model includes latent variables), a possible approach is the use of the modified profile likelihood \citep{BaNi:83}, which consists in a simple modification to the profile likelihood, and may be considered as a generalization of marginal and conditional likelihoods (see e.g. \citealp{WANG2008,BeHa:09,Cadigan2012,BaBeSaSa:16} and the references therein\footnote{For example, \cite{FeReCo:91,DCMaStYo:96,Seve:98,Seve:00,Seve:02,Sato:03,CoxD:06,BeSa:06,PaSa:06,BrDaRe:07}.}). However, these finite sample bias reduction techniques have mainly been studied in low-dimensional settings and their performance in asymptotic regimes where both $n$ and $p$ can tend to infinity, is often unclear.

With complex models and high dimensional data, other estimators are often preferred to the MLE for several reasons which include, for example, its computational difficulty or the fact that the MLE doesn't exist. In addition, one may prefer estimators that ensure robustness against data contamination thereby moving away from the MLE. The choice of other types of estimators can often result in a potential asymptotic bias and, for this reason, different strategies are used to correct them. One of the most common strategies can be put into the framework of the Generalized Method of Moments \citep{Hans:82,Hall:05} which, in essence, defines estimators that are the minimizers of a measure of discrepancy between a (suitable) function of the data and its expectation under the assumed model \citep[see also][]{Gallant1997}. This function is typically the vector of empirical moments but, for example, can also be defined by a bounded score-type function for robustness reasons \citep[see e.g.][]{Hube:81,HaRoRoSt:86}\footnote{For more recent general accounts, see e.g. \citet{MaMaYo:06,HuRo:09,HeCaCoVF:09}.}. The resulting estimators are hence, under certain conditions, Fisher consistent by construction. When the expectation is not available in closed form,  numerical approximations are used, such as the Laplace approximation \citep[see e.g][]{TiKaKa:89}, but the potential bias due to the approximation is usually ignored when developing inference \citep[for examples of the use of the Laplace approximation in the frequentist framework, see e.g.][]{BrCl:93,huber2004estimation,RiVeLe:09}. Alternatively, the simulated method of moments \citep{mcfadden1989method, gallant1996moments, duffie1993simulated} uses simulations to approximate the expectation. A similar strategy for asymptotic bias correction is the method of indirect inference  of \cite{smith1993estimating} and \cite{gourieroux1993indirect} \citep[see also][]{Smit:08}. With this method, the  discrepancy is measured on the difference between the (inconsistent) estimator and an approximation of its expectation under the assumed model which is obtained by means of simulations (for more details, see Section \ref{sec:setting}). Unfortunately, the applicability of these strategies is often limited, for computational reasons (especially when $p$ is large), and their behaviour in high dimensional settings hasn't been clearly determined. 

In this paper, we propose a general framework (discussed in Section \ref{sec:setting}), based on the principle of indirect inference, for finite sample (and asymptotic) bias correction, which produces (computationally efficient) unbiased estimators for (possibly complex) parametric models in high dimensional settings. Moreover, as it will be illustrated by means of an analysis of several common models, the finite sample efficiency is also often improved. We show that (under some conditions) these estimators are consistent and asymptotically normal in Section \ref{Sec_mainRes}. In addition, we demonstrate that the proposed approach allows to obtain unbiased estimators in finite samples. At the same time we address the computational limitations of the indirect inference method by showing that such estimators can be computed using Iterative Bootstrap (IB) proposed by \cite{kuk1995asymptotically}. Some of the results presented in this paper are built upon the results developed by \cite{guerrier2018simulation} in the low-dimensional setting. However, this paper sets both a broader framework as well as new and stronger results compared to existing ones. For example, we extend the framework to inconsistent estimators as well as to discontinuous ones (with respect to the parameters) which provides an appreciable advantage for discrete data models. Therefore, this new framework can be applied to inconsistent estimators in the context of discrete models and/or models with a large number of parameters (compared to the sample size). We illustrate the applicability of these findings on the logistic and negative binomial regression as well as on regularized regression (such as the lasso) in a simulation study where $p$ is ``large'' (depending on the context) compared to $n$. The estimation of the parameters for the logistic regression is investigated in Section \ref{sec:logistic} while the estimation of those for the negative binomial is explored in Section \ref{sec:negbinomial}. Indeed, it is well known that the parameter estimates of the MLE are biased for logistic regression when the ratio $p/n$ is large \citep[see e.g.][]{Scha:83} as is the MLE of the dispersion parameter in the negative binomial regression model \citep[see e.g.][and the references therein]{SaPa:05,LlSm:07}. Furthermore, we also propose new consistent and unbiased estimators for the logistic regression model that are  both ``robust'' to separation as well as to data contamination. As a third example, we consider Regularized Regression Estimators (RRE), such as the ridge regression estimator and the lasso, and provide an RRE with reduced finite-sample bias and mean squared error whose behaviour is studied by means of simulations in low- and high dimensional settings (Section~\ref{Sec_lasso}). Finally, Section \ref{sec:conclude} concludes. 


\section{General setting}
\label{sec:setting}

In order to introduce the framework we intend to study, let us define $\mathbf{X}\left(\bm{\theta}, n, \bm{\omega}\right)\in\real^n$ as being a random sample generated under model $F_{\bm{\theta}}$ (possibly conditional on a set of fixed covariates), where $\bm{\theta}\in\bm{\Theta}\subset\real^p$ is the parameter vector of interest and $\bm{\omega}\in\real^m,\,m\geq n,$ represents a random variable explaining the source of randomness of the sample. More specifically, $\bm{\omega}$ can be considered as a random \emph{seed} that, once observed, produces the sample of size $n$ in a deterministic manner based on a given value of $\bm{\theta}$. Indeed, $\bm{\omega}$ can be conceived as a random variable issued from a model $G$, thereby justifying the definition of the random sample as $\mathbf{X}\left(\bm{\theta}, n, \bm{\omega}\right)$. With this in mind, denoting ``$\overset{d}{=}$'' as ``equality in distribution'', there is no requirement for $\bm{\omega}$ to be unique since it is possible to have $\mathbf{X}\left(\bm{\theta}, n, \bm{\omega}\right)\overset{d}{=}\mathbf{X}\left(\bm{\theta}, n, \bm{\omega}^\ast\right)$ even though $\bm{\omega}\overset{d}{\neq}\bm{\omega}^\ast$. For simplicity, throughout this paper we will consider $\bm{\omega}$ as only belonging to a fixed set $\bOmega \equiv \left\{\bm{\omega}_{l}\right\}_{l \in \mathbb{N}^\ast} \subset \real^m$, where $\mathbb{N}^\ast \equiv \mathbb{N}\setminus \{0\}$. In this setting, defining $\bm{\theta}_0 \in \bm{\Theta}$, we let $\bm{\omega}_0$ and $\bm{\theta}_0$ denote respectively the unknown \textit{random} vector and \textit{fixed} parameter vector used to generate the random sample $\mathbf{X}\left(\bm{\theta}_0, n, \bm{\omega}_0\right)$ that will be used to estimate $\bm{\theta}_0$. Knowing that the sample $\mathbf{X}\left(\bm{\theta}_0, n, \bm{\omega}_0\right)$ is generated from a random seed $\bm{\omega}_0$, we can also consider other samples that we denote as $\mathbf{X}\left(\bm{\theta}, n, \bm{\omega}\right)$, where $(\bm{\theta}, n, \bm{\omega}) \in \bm{\Theta} \times \mathbb{N}^\ast \times \bm{\Omega}$, which can therefore be simulated based on different values of $\bm{\omega}$. Notice that, based on this premise, we have that $\bm{\omega}_0 \not\in \bOmega$ which therefore implies that, in the context of this paper, it is not possible to generate samples using $\bm{\omega}_0$. Hence, the use of $\bm{\omega}$ allows to explicitly display the randomness of a sample. In addition, with this notation it is not only possible to clearly distinguish the observed samples and the simulated ones but also to define the difference (or equivalence) between two simulated samples, say  $\mathbf{X}\left(\bm{\theta}, n, \bm{\omega}_j\right)$ and  $\mathbf{X}\left(\bm{\theta}, n, \bm{\omega}_l\right)$.  

 Having defined the setting of reference for this paper, we can now focus on the estimation procedure aimed at retrieving the value of $\bm{\theta}_0$ from $\mathbf{X}\left(\bm{\theta}_0, n, \bm{\omega}_0\right)$. For this purpose, let us define $\hat{\bm{\pi}}\left(\mathbf{X}\left(\bm{\theta}_0, n, \bm{\omega}_0\right)\right)$ as a biased estimator of $\bm{\theta}_0$ which, despite being possibly inconsistent, is either readily available or can easily be computed. For simplicity of notation, from this point onward we will refer to the estimator $\hat{\bm{\pi}}\left(\mathbf{X}\left(\bm{\theta}_0, n, \bm{\omega}_0\right)\right)$ as $\hat{\bm{\pi}}\left(\bm{\theta}_0, n, \bm{\omega}_0\right)$. Having specified this, in order to correct the (finite sample and asymptotic) bias of $\hat{\bm{\pi}}(\bm{\theta}_0,n,\bm{\omega}_0)$, similarly to \cite{guerrier2018simulation}, we consider an approach based on the principle of indirect inference. 
 Therefore, for all $(j, n, H) \in \mathbb{N} \times \mathbb{N}^\ast \times \mathbb{N}^\ast$, we define the Just identified Indirect inference Estimator (JIE) as  
  \begin{equation}
    \hat{\bm{\theta}}_{(j,n,H)} \in \widehat{\bm{\Theta}}_{(j,n,H)} \equiv \argzero_{\bm{\theta} \in
      \bm{\Theta}} \;  \hat{\bm{\pi}}(\bm{\theta}_0, n, \bm{\omega}_0) - \frac{1}{H} \sum_{h = 1}^H  \hat{\bm{\pi}}(\bm{\theta}, n, \bm{\omega}_{h+jH}),
    \label{eq:indirectInf:hTimesN}
  \end{equation}
where $\hat{\bm{\pi}}(\bm{\theta}, n, \bm{\omega}_{j})$ is obtained on $\mathbf{X}\left(\bm{\theta}, n, \bm{\omega}_j \right)$. This definition assumes that the solution set $\widehat{\bm{\Theta}}_{(j,n,H)}$ is not empty for all $(j, n, H) \in \mathbb{N} \times \mathbb{N}^\ast \times \mathbb{N}^\ast$, which is reasonable when $\dim(\hat{\bm{\pi}}(\bm{\theta}, n, \bm{\omega}))= \dim(\bm{\theta})$. Henceforth, we will refer to $\hat{\bm{\pi}}(\bm{\theta}, n, \bm{\omega})$ as the auxiliary estimator.
\begin{Remark}
	\label{rem:ind:inf}
	The JIE given in \eqref{eq:indirectInf:hTimesN} is a special case of an indirect inference estimator which, using the previous notation, can be defined as 
	\begin{equation}
	    \hat{\bm{\theta}}^\ast_{(j,n,H)} \in \widehat{\bm{\Theta}}^\ast_{(j,n,H)} \equiv \argmin_{\bm{\theta} \in
	      \bm{\Theta}} \; \Big\| \hat{\bm{\pi}}(\bm{\theta}_0, n, \bm{\omega}_0) - \frac{1}{H} \sum_{h = 1}^H  \hat{\bm{\pi}}(\bm{\theta}, n, \bm{\omega}_{j+h-1}) \Big\|^2_{\bm{\Phi}} \; ,
	      \label{eq:indirectInf:rem}
	  \end{equation}
	where $\bm{\Phi}$ is a positive-definite matrix. However, since we suppose that $\widehat{\bm{\Theta}}_{(j,n,H)}\neq \emptyset$  and since $\dim(\hat{\bm{\pi}}(\bm{\theta}_0, n, \bm{\omega}_0)) = \dim(\bm{\theta})$ we have that $\widehat{\bm{\Theta}}_{(j,n,H)} = \widehat{\bm{\Theta}}^\ast_{(j,n,H)}$, which clearly shows the equivalence between (\ref{eq:indirectInf:hTimesN}) and (\ref{eq:indirectInf:rem}), for any positive-definite matrix $\bm{\Phi}$.
	\end{Remark}
\vspace{0.25cm}
Considering the definition of the JIE in \eqref{eq:indirectInf:hTimesN}, it is straightforward to notice that the optimization problem for the second term in the expression can be computationally demanding. However, \cite{guerrier2018simulation} showed that under suitable conditions the IB is a computationally efficient and numerically stable algorithm to compute indirect inference estimators such as the JIE. More specifically, the IB provides the sequence $\left\{\tilde{\bm{\theta}}^{(k)}_{(j,n,H)}\right\}_{k \in \mathbb{N}}$ defined as
	\begin{equation}
		\tilde{\bm{\theta}}^{(k)}_{(j,n,H)} \equiv \tilde{\bm{\theta}}^{(k-1)}_{(j,n,H)} + \left[ \hat{\bm{\pi}}(\bm{\theta}_0, n, \bm{\omega}_0) - \frac{1}{H} \sum_{h = 1}^H  \hat{\bm{\pi}}\left(\tilde{\bm{\theta}}^{(k-1)}_{(j,n,H)}, n, \bm{\omega}_{j+h-1}\right) \right],
		\label{eq:iterboot}
	\end{equation}
	where $\tilde{\bm{\theta}}^{(0)}_{(j,n,H)}\in\bm{\Theta}$. For example, if $\hat{\bm{\pi}}(\bm{\theta}_0, n, \bm{\omega}_0)\in\bm{\Theta}$, then this estimator can be used as the initial value of the above sequence. When this iterative procedure converges, we define $\tilde{\bm{\theta}}_{(j,n,H)}$ as the limit in $k$ of $\tilde{\bm{\theta}}^{(k)}_{(j,n,H)}$. Note that, it is possible that the sequence is stationary in that there exists $k^* \in \mathbb{N}$ such that for all $k \geq k^*$ we have $\tilde{\bm{\theta}}_{(j,n,H)} = \tilde{\bm{\theta}}^{(k)}_{(j,n,H)}$.

\begin{Remark}
  When $\hat{\bm{\pi}}(\bm{\theta}_0, n, \bm{\omega}_0)$ in (\ref{eq:iterboot}) is a consistent estimator of $\btheta_0$, the first step of the IB sequence (i.e. $k=1$) is clearly equivalent to the standard bootstrap bias-correction proposed by \cite{efron1994introduction} which can be used to significantly reduce the bias of a consistent estimator (under some appropriate conditions). Nevertheless, \cite{guerrier2018simulation} showed that the IB provides more substantial bias reduction and, as will be shown later, can also be applied to inconsistent estimators.
\end{Remark}

\vspace{0.25cm}

Having defined the JIE, the next section studies the properties of the JIE using the IB as an efficient algorithm to compute this estimator.

\section{Properties of Simulation-based Estimators}
\label{Sec_mainRes}


In this section we discuss the convergence of the IB sequence along with the properties of the JIE. However, since we build upon the results of \cite{guerrier2018simulation}, before discussing the mentioned properties we briefly highlight how the findings in this section deliver both a broader framework as well as new and stronger results compared to existing ones. A first example of such stronger results consists in the demonstration that the JIE is unbiased in finite samples as opposed to previous findings in which it was shown that, at best, the (elementwise) bias was of order $\mathcal{O}(n^{-3})$. Starting from this example, we can additionally state that this new framework allows for the following generalizations compared to existing approaches: \textit{(i)} the auxiliary estimator may have an asymptotic bias that depends on $\bm{\theta}$; \textit{(ii)} the auxiliary estimator may also be discontinuous on $\bm{\theta}$ (as is commonly the case when considering discrete data models); \textit{(iii)} the finite sample bias of the auxiliary estimator is allowed to have a more general expression with respect to previously defined bias functions; \textit{(iv)} some of the topological requirements on $\bm{\Theta}$ are relaxed \textit{(v)} $p$ is allowed to increase  with the sample size $n$.

Having specified how the findings of this section improve over existing results, we can now deliver the details of this study which addresses the convergence of the IB sequence (Section \ref{sec:boot}) together with the bias-correction properties (Section \ref{sec:bias}), consistency (Section \ref{sec:consist}) and asymptotic normality (Section \ref{sec:norm}) of the JIE.
This order of presentation is justified by the rationale that, before considering the properties of the JIE, one should make sure that it can be computed. Indeed, in many practical situations (for example when the model at hand is complex and/or when $p$ is relatively large) the JIE is difficult to compute when using, for example, numerical procedures to directly solve (\ref{eq:indirectInf:hTimesN}). It is therefore legitimate to first assess if the IB can be used as a suitable approach to compute the JIE before discussing its statistical properties. Our ultimate goal is to provide the sufficient conditions for unbiasedness, consistency and asymptotic normality of $\tilde{\bm{\theta}}_{(j,n,H)}$ (i.e. the limit of the IB sequence). This is done in Section \ref{sec:discussion}, where the assumptions used to derive the different properties of the IB and JIE are put in a global framework. In fact, Section \ref{sec:discussion} contains the main theoretical results of this paper. 

Various assumptions are considered and, in order to provide a general reference for the type of assumption, we employ the following conventions to name them:
\begin{itemize}
    \item Assumptions indexed by ``A'' refer to the topology of $\bm{\Theta}$;
    \item Assumptions indexed by ``B'' refer to the existence of the expected value of the auxiliary estimator;
    \item Assumptions indexed by ``C'' refer to the random component of the auxiliary estimator;
    \item Assumptions indexed by ``D'' refer to the bias of the auxiliary estimator.
\end{itemize}

 When needed, sub-indices will be used to distinguish assumptions of the same type (for example Assumptions \ref{assum:A:1} and \ref{assum:A:2} are different but both related to the topology of $\bm{\Theta}$) and assumptions with two sub-indices are used to denote an assumption that implies two separate assumptions (for example Assumption  \ref{assum:A:13} implies both \ref{assum:A:1} and \ref{assum:A:3}). 

\begin{Remark}
\label{asym:framework}
The asymptotic framework we use in this section is somewhat unusual as we always consider arbitrarily large but finite $n$ and $p$. Indeed, if $p$ is such that $p \to \infty$ as $n \to \infty$, our asymptotic results may not be valid when taking the limit in $n$ (and therefore in $p$) but are valid for all finite $n$ and $p$. The difference between the considered framework and others where limits are studied is rather subtle as, in practice, asymptotic results are typically used to derive approximations in finite samples for which the infinite (dimensional) case is not necessarily informative. Indeed, the topology of $\real^p$ for any finite $p$ and $\real^\infty$ are profoundly different. For example, the consistency of an estimator is generally dependent on the assumption that $\bm{\Theta}$ is compact. In the infinite dimensional setting closed and bounded sets are not necessarily compact. Therefore, this assumption becomes rather implausible for many statistical models and would imply, among other things, a detailed topological discussion of requirements imposed on $\bm{\Theta}$. Similarly, many of the mathematical arguments presented in this paper may not apply in the infinite dimensional case. Naturally, when $p \to c <\infty$ as $n \to \infty$, limits in $n$ are allowed. Although not necessary, we abstain from using statements or definitions such as ``$\,\xrightarrow{p}$'' (convergence in probability) or ``$\,\xrightarrow{d}\,$'' (convergence in distribution) which may lead one to believe that the limit in $n$ exists. This choice is simply made to avoid confusion.
\end{Remark}
\vspace{0.25cm}

\subsection{Convergence of the Iterative Bootstrap}
\label{sec:boot}

In order to study the properties of the sequence defined in (\ref{eq:iterboot}), we set the following assumption framework. Our first assumption concerns the topology of the parameter space $\bm{\Theta}$. 
	
 \setcounter{Assumption}{0}
    \renewcommand\theAssumption{\Alph{Assumption}$_1$}
    \begin{Assumption}
		\label{assum:A:1}
		Let $\bm{\Theta}$ be bounded and such that
		\begin{equation*}
			\argzero_{\bm{\theta} \in
		      \real^p \setminus \bm{\Theta}} \;  \hat{\bm{\pi}}(\bm{\theta}_0, n, \bm{\omega}_0) - \frac{1}{H} \sum_{h = 1}^H  \hat{\bm{\pi}}(\bm{\theta}, n, \bm{\omega}_{h+jH})
		      = \emptyset\,.
		\end{equation*}
\end{Assumption}
\vspace{0.25cm}

Assumption \ref{assum:A:1} is very mild and essentially ensures that no solution outside of $\bm{\Theta}$ exists for the JIE defined in (\ref{eq:indirectInf:hTimesN}). Therefore, the solution set $\widehat{\bm{\Theta}}_{(j,n,H)}$ may also be written as
\begin{equation*}
	    \widehat{\bm{\Theta}}_{(j,n,H)} \equiv \argzero_{\bm{\theta} \in
	      \real^p} \;  \hat{\bm{\pi}}(\bm{\theta}_0, n, \bm{\omega}_0) - \frac{1}{H} \sum_{h = 1}^H  \hat{\bm{\pi}}(\bm{\theta}, n, \bm{\omega}_{h+jH}).
\end{equation*}
In addition, Assumption \ref{assum:A:1} doesn't require $\bm{\Theta}$ to be a compact set since, for example, it can represent an open set. Next we impose some conditions on the auxiliary estimator $\hat{\bm{\pi}}(\bm{\theta},n, \bm{\omega})$ thereby defining the second assumption that has no associated sub-index as it will remain unchanged throughout the paper. Denoting $\mathbf{a}_k$ as being the $k$-th entry of a generic vector $\mathbf{a} \in \real^p$, the assumptions is as follows.
    \setcounter{Assumption}{1}
    \renewcommand{\theHAssumption}{otherAssumption\theAssumption}
    \renewcommand\theAssumption{\Alph{Assumption}}
\begin{Assumption}
	\label{assum:B:1}
	For all $\left(\bm{\theta}, \, n, \, \bm{\omega}\right) \in \bm{\Theta} \times \mathbb{N}^\ast \times \bOmega$, the following expectation exists and is finite, i.e.
		\begin{equation*}
			\bm{\pi} \left( \bm{\theta}, n\right) \equiv \mathbb{E}\left[\hat{\bm{\pi}}(\bm{\theta},n, \bm{\omega})\right] < \infty.
		\end{equation*}
		Moreover, $\bm{\pi}_k \left( \bm{\theta}\right) \equiv \displaystyle{\lim_{n \to \infty}} \bm{\pi}_k \left( \bm{\theta}, n\right)$ exists for all $k = 1, \ldots, p$.
\end{Assumption}
\vspace{0.25cm}
	The second assumption is likely  to be satisfied in the majority of practical situations (for instance if $\bm{\Theta}$ is bounded) and is particularly useful as it allows to decompose the estimator $\hat{\bm{\pi}}(\bm{\theta},n, \bm{\omega})$ into a non-stochastic component $\bm{\pi} \left( \bm{\theta}, n\right)$ and a random term $\mathbf{v} \left(\bm{\theta},n, \bm{\omega}\right)$. Indeed, using Assumption \ref{assum:B:1}, we can write:
	\begin{equation}
	   \hat{\bm{\pi}}(\bm{\theta},n,\bm{\omega})  = \bm{\pi} \left( \bm{\theta}, n\right)+ \mathbf{v} \left(
	    \bm{\theta},n, \bm{\omega}\right),
	  \label{bias:estim}
	\end{equation}
where $\mathbf{v} \left(\bm{\theta},n, \bm{\omega}\right) \equiv \hat{\bm{\pi}}(\bm{\theta},n,\bm{\omega}) - \bm{\pi} \left( \bm{\theta}, n\right)$ is a zero-mean random vector. It must be noticed that this assumption doesn't imply that the functions $\bm{\pi} \left( \bm{\theta}, n\right)$ and $\hat{\bm{\pi}}(\bm{\theta},n, \bm{\omega})$ are continuous in $\btheta$. Consequently, when using the IB to obtain the JIE (see Proposition \ref{thm:iter:boot} below), it is possible to make use of auxiliary estimators  that are not continuous in $\btheta$. In \cite{guerrier2018simulation}, the continuity of these functions was (implicitly) assumed while the weaker assumptions presented in this paper extend the applicability of this framework also to models for discrete data such as GLMs (see Sections \ref{sec:logistic} and \ref{sec:negbinomial} for the logistic and negative binomial regression models, respectively). We now move onto
	Assumption \ref{assum:C:1} (below) that imposes some restrictions on the random vector $\mathbf{v} \left(\bm{\theta},n, \bm{\omega}\right)$. 
    
    \setcounter{Assumption}{2}
    \renewcommand{\theHAssumption}{otherAssumption\theAssumption}
    \renewcommand\theAssumption{\Alph{Assumption}$_1$}
	\begin{Assumption}
		\label{assum:C:1}
		There exists a real $\alpha > 0$ such that for all $\left(\bm{\theta}, \, \bm{\omega}\right) \in \bm{\Theta} \times \bOmega$ and all $k = 1,\,\ldots,\,p$, we have
		\begin{equation*}
			\mathbf{v}_k(\bm{\theta},n, \bm{\omega}) = \mathcal{O}_{\rm p}(n^{-\alpha})\;\;\;\;
			 \text{and}\;\;\;\;
			 \lim_{n \to \infty} \; \frac{p^{\nicefrac{1}{2}}}{n^\alpha} = 0.
		\end{equation*}

	\end{Assumption}
	\vspace{0.25cm}
	
	Assumption \ref{assum:C:1} is frequently employed and typically very mild as it simply requires that the variance of $\mathbf{v} \left(\bm{\theta},n, \bm{\omega}\right)$ goes to zero as $n$ increases. For example, if $\hat{\bm{\pi}}(\bm{\theta},n,\bm{\omega})$ is $\sqrt{n}$-consistent (towards $\bm{\pi} \left( \bm{\theta}\right)$) then we would have $\alpha = \nicefrac{1}{2}$. Next, we consider the bias of $\hat{\bm{\pi}}(\bm{\theta},n,\bm{\omega})$ for $\bm{\theta}$ and we let $\mathbf{d}\left(\bm{\theta}, n\right) \equiv \mathbb{E} \left[\hat{\bm{\pi}}(\bm{\theta},n,\bm{\omega})\right] - \bm{\theta}$. Using this definition we can rewrite (\ref{bias:estim}) as follows: 
	\begin{equation}
	   \hat{\bm{\pi}}(\bm{\theta},n,\bm{\omega})  = \bm{\theta} + \mathbf{d} \left( \bm{\theta}, n\right)+ \mathbf{v} \left(
	    \bm{\theta},n, \bm{\omega}\right).
	  \label{bias:estim2}
	\end{equation}
	Moreover, the bias function $\mathbf{d}\left(\bm{\theta}, n\right)$ can always be expressed as follows
	\begin{equation}
	  \mathbf{d}\left(\bm{\theta}, n\right) = \mathbf{a}(\bm{\theta})+\mathbf{c}(n) + \mathbf{b} \left(\bm{\theta} , n\right),
	  \label{def:fct:bias:d}
	\end{equation}
	where $\mathbf{a}(\btheta)$ is defined as the asymptotic bias in the sense that $\mathbf{a}_k(\bm{\theta}) \equiv \displaystyle{\lim_{n \to \infty}} \; \mathbf{d}_k\left(\bm{\theta}, n\right)$ for $k = 1, \ldots, p$, while
	$\mathbf{c}(n)$ and $\mathbf{b} \left(\bm{\theta} , n\right)$ are used to represent the finite sample bias. 
	More precisely, $\mathbf{a}(\btheta)$ contains all the terms that are strictly functions of $\btheta$, $\mathbf{c}(n)$ the terms that are strictly functions of $n$ and $\mathbf{b} \left(\bm{\theta}, n\right)$ the rest. This definition implies that if $\mathbf{d}\left(\bm{\theta}, n\right)$ contains a constant term, it is included in $\mathbf{a}(\btheta)$.  Moreover, the function $\mathbf{b} \left(\bm{\theta} , n\right)$ can always be decomposed into a linear and a non-linear term in $\btheta$, i.e.
		\begin{equation}
			\mathbf{b}\left(\bm{\theta}, n\right) = \mathbf{L}(n) \bm{\theta} + \mathbf{r} \left(\bm{\theta} , n\right),
			\label{def:fct:bias:b}
		\end{equation}
		where $\mathbf{L}(n) \in \real^{p \times p}$ and $\mathbf{r}\left(\bm{\theta}, n\right)$ doesn't contain any linear term in $\bm{\theta}$. Denoting $\mathbf{A}_{k,l}$ as the entry in the $k$-th row and $l$-th column of a generic matrix $\mathbf{A} \in \real^{p \times p}$, Assumption \ref{assum:D:1} (below) imposes some restrictions on the bias function. 
    \setcounter{Assumption}{3}
    \renewcommand\theAssumption{\Alph{Assumption}$_1$}
	\begin{Assumption}
		\label{assum:D:1}
		The bias function $\mathbf{d}\left(\bm{\theta}, n\right)$ is such that:
		\begin{enumerate}
		    \item The function $\mathbf{a}(\bm{\theta})$ is a contraction map in that for all
		    $\bm{\theta}_1,\bm{\theta}_2\in\bm{\Theta}$ such that $\bm{\theta}_1 \neq \bm{\theta}_2$ we have 
		    \begin{equation*}
		        \Big\lVert \mathbf{a}(\bm{\theta}_2)-\mathbf{a}(\bm{\theta}_1)  \Big\rVert_2 <  \big\rVert \bm{\theta}_2-\bm{\theta}_1 \big\lVert_2\,.
		    \end{equation*}
		    
		\item There exist real $\beta, \gamma > 0$ such that for all $\bm{\theta} \in \bm{\Theta}$ and all $k,l = 1,\,\ldots, \, p$, we have
		\begin{equation*}
		\begin{aligned}
			  \mathbf{L}_{k,l}(n) = \mathcal{O}(n^{-\beta}), \;\;\;  \mathbf{r}_k\left(\bm{\theta}, n\right) = \mathcal{O}(n^{-\gamma}),\;\;\; 
			  \lim_{n \to \infty} \; \frac{p}{n^\beta} = 0,  \;\;\;\ \text{and} \;\;\;
			  \lim_{n \to \infty} \; \frac{p^{\nicefrac{1}{2}}}{n^\gamma} = 0.
		 \end{aligned}
		\end{equation*}
		\end{enumerate}
	\end{Assumption}
    \vspace{0.25cm}
    
	The first part of Assumption \ref{assum:D:1} is reasonable provided that the asymptotic bias is relatively ``small''  compared to $\btheta$ (up to a constant term). For example, if $\mathbf{a}(\btheta)$ is sublinear in $\btheta$, i.e. $\mathbf{a}(\btheta) = \mathbf{M} \btheta + \mathbf{s}$, then the first part of Assumption \ref{assum:D:1} would be satisfied if the Frobenius norm is such that $||\mathbf{M}||_F < 1$ since
	\begin{equation*}
	    \Big\lVert \mathbf{a}(\bm{\theta}_2)-\mathbf{a}(\bm{\theta}_1)  \Big\rVert_2 = \Big\lVert \mathbf{M}(\bm{\theta}_2-\bm{\theta}_1)  \Big\rVert_2 \leq  \big\lVert\mathbf{M}\big\rVert_F \big\lVert \bm{\theta}_2-\bm{\theta}_1 \big\rVert_2 < \big\lVert \bm{\theta}_2-\bm{\theta}_1 \big\rVert_2 .
	\end{equation*}
	Moreover, the function ${\bm{\pi}}(\bm{\theta})$ (as defined in Assumption \ref{assum:B:1}) is often called the asymptotic binding function in the indirect inference literature \citep{gourieroux1993indirect}. To ensure the consistency of such estimators, it is typically required for ${\bm{\pi}}(\bm{\theta})$ to be continuous and injective. In our setting, this function is given by ${\bm{\pi}}(\bm{\theta}) = \btheta + \mathbf{a}(\btheta)$ and its continuity and injectivity are directly implied by the first part of Assumption \ref{assum:D:1}. Indeed, since $\mathbf{a}(\bm{\theta})$ is a contraction map, it is continuous. Moreover, taking $\bm{\theta}_1,\bm{\theta}_2\in\bm{\Theta}$ with
	${\bm{\pi}}(\bm{\theta}_1) = {\bm{\pi}}(\bm{\theta}_2)$, then $
	    \big\lVert {\mathbf{a}}(\bm{\theta}_1) - {\mathbf{a}}(\bm{\theta}_2) \big\rVert = \big\lVert \bm{\theta}_1 - \bm{\theta}_2 \big\rVert$,
	which is only possible if $\btheta_1 = \btheta_2$. Thus, $\bm{\pi}(\bm{\theta})$ is injective.
	
	\begin{Remark}
	\label{Remark:kickIB}
	   In situations where one may suspect that $\mathbf{a}(\bm{\theta})$ is not a  contraction map, a possible solution is to modify the sequence considered in the IB as follows:
	   \begin{equation*}
	       	\tilde{\bm{\theta}}^{(k)}_{(j,n,H)} \equiv \tilde{\bm{\theta}}^{(k-1)}_{(j,n,H)} + \varepsilon_k\left[ \hat{\bm{\pi}}(\bm{\theta}_0, n, \bm{\omega}_0) - \frac{1}{H} \sum_{h = 1}^H  \hat{\bm{\pi}}\left(\tilde{\bm{\theta}}^{(k-1)}_{(j,n,H)}, n, \bm{\omega}_{h+jH}\right) \right],
	   \end{equation*}
	   with $\varepsilon_k \in (0, 1]$ for all $k \in \mathbb{N}$. If $\varepsilon_k=\varepsilon$ (i.e. a constant), $\mathbf{a}(\bm{\theta})$ does not need to be a contraction map. Indeed, if $\bm{\Theta}$ is bounded and $\mathbf{a}(\bm{\theta})$ is differentiable, it is always possible to find an $\varepsilon$ such that $\varepsilon \, \mathbf{a}(\bm{\theta})$ is a contraction map. A formal study on the influence of $\varepsilon_k$ on the IB algorithm is, however, left for further research.  
	\end{Remark}
\vspace{0.25cm}
While more general, the second part of Assumption \ref{assum:D:1} would be satisfied, for example, if $\mathbf{b}(\bm{\theta},n)$ is a sufficiently smooth function in $\bm{\theta}/n$, thereby allowing a Taylor expansion, as considered, for example, in \cite{guerrier2018simulation}. Moreover, this assumption is typically less restrictive than the approximations that are commonly used to describe the bias of estimators. Indeed, a common assumption is that the bias of a consistent estimator (including the MLE), can be expanded in a power series in $n^{-1}$ (see e.g. \citealp{kosmidis2014bias} and \citealp{hall1988bootstrap} in the context of the iterated bootstrap), i.e.,
\begin{equation}
    \mathbf{d}(\bm{\theta},n)  = \sum_{j = 1}^{k} \frac{\mathbf{h}^{(j)}(\bm{\theta})}{n^{j}} + \mathbf{g}(\btheta, n),
    \label{Eq_bais-standard}
\end{equation}
where $\mathbf{h}^{(j)}(\bm{\theta})$ is $\mathcal{O}(1)$ elementwise, for $j = 1, \ldots, k$, and $\mathbf{g}(\btheta, n)$ is $\mathcal{O}\left(n^{-(k+1)}\right)$ elementwise, for some $k \geq 1$. The bias function $\mathbf{d}(\bm{\theta},n)$ given in (\ref{Eq_bais-standard}) clearly satisfies the requirements of Assumption~\ref{assum:D:1}.

	Moreover, under the form of the bias postulated in (\ref{Eq_bais-standard}), we have that $\beta, \gamma \geq 1$. If the auxiliary estimator is $\sqrt{n}$-consistent, we have $\alpha = \nicefrac{1}{2}$, and therefore the requirements of Assumptions \ref{assum:C:1} and \ref{assum:D:1}, i.e.
	\begin{equation*}
	    \lim_{n \to \infty} \; \max \left( \frac{p^{\nicefrac{1}{2}}}{n^{\min(\alpha,\gamma)}},  \frac{p}{n^\beta} \right) = 0,
	\end{equation*}
	are satisfied if
	\begin{equation*}
	    \lim_{n \to \infty} \; \frac{p}{n} = 0.
	\end{equation*}

	We now study the convergence properties of the IB (defined in (\ref{eq:iterboot})) when used to obtained the JIE presented in (\ref{eq:indirectInf:hTimesN}). In Lemma \ref{lemma:unique:iter:boot}, we show that when $n$ is sufficiently large, the solution space $\widehat{\bm{\Theta}}_{(j,n,H)}$ of the JIE contains only one element. In other words, this result ensures that the function $\frac{1}{H} \sum_{h = 1}^H  \hat{\bm{\pi}}(\bm{\theta}, n, \bm{\omega}_{h+jH})$  in (\ref{eq:indirectInf:hTimesN}) is injective for fixed (but possibly large) sample size $n$ and fixed $H$. Lemma \ref{lemma:unique:iter:boot} formally states this result and its proof is given in Appendix \ref{proof:lemma:1}.
	%
    %
    
    %
    
    %
    
	\begin{Lemma}
		\label{lemma:unique:iter:boot}
		Under Assumptions \ref{assum:A:1}, \ref{assum:B:1}, \ref{assum:C:1} and \ref{assum:D:1}, for all $(j, H) \in \mathbb{N} \times \mathbb{N}^\ast$ and $n \in \mathbb{N}^\ast$ sufficiently large, the solution space $\widehat{\bm{\Theta}}_{(j,n,H)}$ is a singleton.
	\end{Lemma}
	\vspace{0.25cm}
    Building on the identification result given in the above lemma, the following proposition states the equivalence of the JIE in (\ref{eq:indirectInf:hTimesN}) with the limit of the IB sequence in (\ref{eq:iterboot}).
	\begin{Proposition}
		\label{thm:iter:boot}
		Under Assumptions \ref{assum:A:1}, \ref{assum:B:1}, \ref{assum:C:1} and \ref{assum:D:1}, for all $(j, H) \in \mathbb{N} \times \mathbb{N}^\ast$ and $n \in \mathbb{N}^\ast$ sufficiently large, the sequence $\left\{\tilde{\bm{\theta}}^{(k)}_{(j,n,H)}\right\}_{k \in \mathbb{N}}$ is such that
		\begin{equation*}
    	\tilde{\bm{\theta}}_{(j,n,H)} = \lim_{k \to \infty} \; \tilde{\bm{\theta}}^{(k)}_{(j,n,H)} = \hat{\bm{\theta}}_{(j,n,H)}.
		\end{equation*}
		Moreover, there exists a real $ \epsilon \in (0, \, 1)$ such that for any $k \in \mathbb{N}^\ast$ 
		\begin{equation*}
			\left\lVert\tilde{\bm{\theta}}^{(k)}_{(j,n,H)} - \hat{\bm{\theta}}_{(j,n,H)}\right\rVert_2  =\mathcal{O}_{\rm p}({p}^{\nicefrac{1}{2}}\,\epsilon^k).
		\end{equation*}
	\end{Proposition}
	\vspace{0.25cm}
	The proof of this result is given in Appendix \ref{proof:thm:1}. An important consequence of Proposition \ref{thm:iter:boot} is that the IB provides a computationally efficient algorithm to solve the optimization problem that defines the JIE in (\ref{eq:indirectInf:hTimesN}). In practical settings, the IB can often be applied to the estimation of complex models where standard optimization procedures used to solve (\ref{eq:indirectInf:hTimesN}) may fail to converge numerically (see e.g. \citealp{guerrier2018simulation}). Moreover, the IB procedure is generally much faster, since a consequence of Proposition \ref{thm:iter:boot} is that $\tilde{\bm{\theta}}^{(k)}_{(j,n,H)}$ converges to $\hat{\bm{\theta}}_{(j,n,H)}$ (in norm) at an exponential rate. However, the convergence of the algorithm may be slower when $p$ is large.

\subsection{Bias Correction Properties of the JIE}

\label{sec:bias}

 In this section we study the finite sample bias of the JIE $\hat{\bm{\theta}}_{(j,n,H)}$. For this purpose, we first consider a slightly modified assumption framework where Assumptions \ref{assum:A:1} and \ref{assum:D:1} are modified, while Assumption \ref{assum:B:1} remains unchanged and Assumption \ref{assum:C:1} becomes irrelevant. 
 Firstly, the modified version of Assumption \ref{assum:A:1} is as follows:
    \setcounter{Assumption}{0}
    \renewcommand{\theHAssumption}{otherAssumption\theAssumption}
    \renewcommand\theAssumption{\Alph{Assumption}$_2$}
    \begin{Assumption}
    \label{assum:A:2}
     $\btheta_0 \in \operatorname{Int}(\bm{\Theta})$ and $\bm{\Theta}$ is a convex and compact set such that
    \begin{equation*}
			\argzero_{\bm{\theta} \in
		      \real^p \setminus \bm{\Theta}} \;  \hat{\bm{\pi}}(\bm{\theta}_0, n, \bm{\omega}_0) - \frac{1}{H} \sum_{h = 1}^H  \hat{\bm{\pi}}(\bm{\theta}, n, \bm{\omega}_{h+jH})
		      = \emptyset\,.
	\end{equation*}
    \end{Assumption}
    \vspace{0.25cm}
    This topological assumption is quite mild although stronger than necessary. Indeed, the convexity condition and the fact that $\bto$ is required to be in the interior of $\bm{\Theta}$ are convenient to ensure that expansions can be made between $\bto$ and an arbitrary point in $\bm{\Theta}$. The rest of this assumption ensures that certain random variables that we will consider will be bounded. Clearly, Assumption \ref{assum:A:2} implies Assumption \ref{assum:A:1}. Moving our focus to Assumption \ref{assum:D:1}, the latter is modified by imposing a different set of restrictions on the bias function.

    \setcounter{Assumption}{3}
    \renewcommand{\theHAssumption}{otherAssumption\theAssumption}
    \renewcommand\theAssumption{\Alph{Assumption}$_2$}
    \begin{Assumption}
    \label{assum:D:2}
    Consider the functions $\mathbf{d}(\bm{\theta}, n)$, $\mathbf{a}(\btheta)$, $\mathbf{b}(\bm{\theta}, n)$ and $\mathbf{r}(\bm{\theta}, n)$ defined in \eqref{def:fct:bias:d} and \eqref{def:fct:bias:b} and  let $\mathbf{R} (\bm{\theta}, n) \equiv \frac{\partial}{\partial \, \btheta^T}  \mathbf{r}(\bm{\theta}, n) \in \real^{p \times p}$ be the Jacobian matrix of $\mathbf{r} (\bm{\theta}, n)$ in $\btheta$.  These functions are, for all $\btheta \in \bm{\Theta}$, such that:
    \begin{enumerate}
        \item There exists a matrix $\mathbf{M} \in \real^{p \times p}$ and a vector $\mathbf{s} \in \real^{p}$ such that $\mathbf{a}(\btheta) = \mathbf{M}\btheta + \mathbf{s}$.
        \item There exists a sample size $n^* \in \mathbb{N}^*$ and a real $\gamma > 0$ such that for all $n \in \mathbb{N}^*$ satisfying $n \geq n^*$, all $\btheta \in \bm{\Theta}$ and all $k,l = 1,\ldots, p$, we have 
    \begin{equation*}
	  \begin{aligned}
			 \mathbf{r}_k\left(\bm{\theta}, n\right) = \mathcal{O}(n^{-\gamma}) \;\;\; \text{and} \;\;\; \mathbf{R}_{k,l} (\bm{\theta}, n)  \text{ is continuous in }\btheta\in\bm{\Theta}.
		\end{aligned}
	\end{equation*}


    \item There exists a sample size $n^* \in \mathbb{N}^*$ such that for all $n \in \mathbb{N}^*$ satisfying $n \geq n^*$:
    \begin{itemize}
        \item the matrix $(\mathbf{M} + \mathbf{L}(n)+\mathbf{I})^{-1}$ exists,
        \item $p^2 \,n^{-\gamma} < 1$.
    \end{itemize}
    
    \end{enumerate}
    \end{Assumption}
    \vspace{0.25cm}
    
     Clearly, neither Assumption \ref{assum:D:1} nor Assumption \ref{assum:D:2} imply each other. The first part of Assumption \ref{assum:D:2} imposes a restrictive form for the bias function $\mathbf{a}(\btheta)$. If $\hat{\bm{\pi}}(\bm{\theta}_0, n, \bm{\omega})$ is a consistent estimator of $\bm{\theta}_0$, then this restriction is automatically satisfied. However, the matrix $\mathbf{M}$ may be such that $\lVert\mathbf{M}\rVert_F \geq 1$, thereby relaxing the contraction mapping hypothesis on $\mathbf{a}(\btheta)$ given in Assumption \ref{assum:D:1}. {\color{black}The second part of Assumption \ref{assum:D:2} implies that $\mathbf{r} (\bm{\theta}, n)$ is continuous in $\btheta\in\bm{\Theta}$. Another consequence, taking into account Assumption \ref{assum:A:2}, is that $\mathbf{r}_k(\hbt, n)$ and $\mathbf{R}_{k,l} (\hbt, n)$ are bounded random variables for all $k,l=1,\dots,p$.}
    The last of part of Assumption \ref{assum:D:2} requires that the matrix $(\mathbf{M} + \mathbf{L}(n)+\mathbf{I})^{-1}$ exists when $n$ is sufficiently large. This requirement is quite general and is, for example, satisfied if $\mathbf{a}(\btheta)$ is relatively ``small'' compared to $\bm{\theta}$ or if $\hat{\bm{\pi}}(\bm{\theta}_0, n, \bm{\omega})$ is a consistent estimator of $\bm{\theta}_0$. Interestingly, this part of the assumption can be interpreted as requiring that the matrix $(\mathbf{M} + \mathbf{I})^{-1}$ exists, which directly implies that the binding function $\bm{\pi}(\btheta)$ is injective. Finally, in many practical settings it is reasonable\footnote{For example using the bias function proposed in (\ref{Eq_bais-standard}) and assuming that the first term of the expansion $\mathbf{h}^{(1)}(\btheta)$ is linear.} to assume that $\gamma =2$ and the condition that $p^2 \,n^{-\gamma} < 1$ is satisfied if $p/n < 1$ for sufficiently large $n$, which is less strict than the conditions in Assumption \ref{assum:D:1}.

     Under these new conditions, the unbiasedness of $\hat{\bm{\theta}}_{(j,n,H)}$ is established in Proposition \ref{THM:bias} below in the case where $p$ is bounded. The proof of this result is given in Appendix \ref{proof:thm:bias}.
    
    \begin{Proposition}  
		\label{THM:bias}
		Under Assumptions \ref{assum:A:2}, \ref{assum:B:1} and \ref{assum:D:2}, there exists a sample size $n^* \in \mathbb{N}^*$ such that for all $n \in \mathbb{N}^*$ satisfying $n \geq n^*$ and for all $(j, H) \in \mathbb{N} \times \mathbb{N}^\ast$, we have 
		\begin{equation*}
		\Big\lVert	\mathbb{E} \left[\hbt \right] - \bto  \Big\rVert_2 = 0.  
		\end{equation*}
	\end{Proposition}
	\vspace{0.25cm}
	
	It must be underlined that this proposition doesn't deliver any asymptotic argument as it is valid for any $n\in\mathbb{N}^*$ such that $n\geq n^*$. Therefore, this result shows that the JIE is unbiased in finite (but possibly large) sample sizes. {\color{black}Indeed, the requirement of $n\geq n^*$ is linked to Assumption \ref{assum:D:2} and if the second and third part of this assumption are satisfied for a given $n^*$, then the result holds for this $n^*$.} In practical settings, the value $n^*$ appears to be very small as in the majority of the simulated settings considered for this paper (some of which are shown in Sections \ref{sec:logistic}, \ref{sec:negbinomial}, \ref{Sec_lasso}) the resulting estimator does indeed appear to be unbiased.
    Moreover, in the proof of Proposition \ref{THM:bias} we actually show that for all $k\in\mathbb{N}^*$ we have $\Big\lVert	\mathbb{E} \left[\hbt \right] - \bto  \Big\rVert_2 ~=~\mathcal{O}\left((p^2n^{-\gamma})^k\right)$, which implies that $\Big\lVert	\mathbb{E} \left[\hbt \right] - \bto  \Big\rVert_2 = 0$. This result generalizes the results of \cite{guerrier2018simulation} who, under far more restrictive assumptions, showed that $\mathbb{E} [\hat{\bm{\theta}}_{(j,n,H)} ] - \bm{\theta}_0$ is $\mathcal{O}\left(n^{-\delta}\right)$, elementwise, where $\delta \in (2, 3]$.
	
	At this point, we must highlight that Proposition \ref{THM:bias} relies on the assumption that the asymptotic bias $\mathbf{a}(\btheta)$ is sublinear, which can be quite restrictive. Nevertheless, in Section \ref{sec:consist}, we show that, under some additional conditions, $\hat{\bm{\theta}}_{(j,n,H)}$ is a consistent estimator of $\bm{\theta}_0$. Therefore, a possible approach that would guarantee an unbiased estimator is to obtain the JIE $\hat{\bm{\theta}}_{(j,n,H)}$ from an inconsistent auxiliary estimator $\hat{\bm{\pi}}(\bm{\theta},n,\bm{\omega})$ and then to compute a new JIE with auxiliary estimator $\hat{\bm{\theta}}_{(j,n,H)}$. In practice however, this computationally intensive approach is probably unnecessary since the (one step) JIE appears to eliminate the bias almost completely when considering inconsistent auxiliary estimators $\hat{\bm{\pi}}(\bm{\theta},n,\bm{\omega})$ (see e.g. Section \ref{sec:logistic}). In fact, we conjecture that for a larger class of asymptotic bias functions, Proposition \ref{THM:bias} remains true, although the verification of this conjecture is left for further research.

\subsection{Consistency of the JIE}
\label{sec:consist}

We now study the consistency properties of the JIE and, as for the previous section, we again slightly modify the assumption framework. More in detail, we impose stronger requirements on $\alpha$ and $\gamma$ in  Assumptions \ref{assum:C:1} and \ref{assum:D:2}. 

	\setcounter{Assumption}{2}
	\renewcommand{\theHAssumption}{otherAssumption\theAssumption}
	\renewcommand\theAssumption{\Alph{Assumption}$_2$}
	\begin{Assumption}
	\label{assum:C:2}
	Preserving the same definitions and requirements of Assumption \ref{assum:C:1}, we additionally require that the constant $\alpha$ is such that
	\begin{equation*}
	   \lim_{n \to \infty} \frac{p^{\nicefrac{5}{2}}}{n^\alpha}=0.
	\end{equation*}
	\end{Assumption}
	\vspace{0.25cm}
	
	\setcounter{Assumption}{3}
	\renewcommand{\theHAssumption}{otherAssumption\theAssumption}
	\renewcommand\theAssumption{\Alph{Assumption}$_3$}
	\begin{Assumption}
	\label{assum:D:3}
	Preserving the same definitions and requirements of Assumption \ref{assum:D:2}, we additionally require that the constant $\gamma$ is such that
	\begin{equation*}
	   \lim_{n \to \infty} \frac{p^{\nicefrac{5}{2}}}{n^\gamma}=0.
	\end{equation*}
	\end{Assumption}
	\vspace{0.25cm}
	
	Using these assumptions, we present the following corollary to Proposition \ref{THM:bias}. 
	
	\begin{Corollary}
	\label{coro:consist}
	Under Assumptions \ref{assum:A:2}, \ref{assum:B:1}, \ref{assum:C:2}, and \ref{assum:D:3}, $\hat{\bm{\theta}}_{(j,n,H)}$ is a consistent estimator of $\bm{\theta}_0$ for all $(j, H) \in \mathbb{N} \times \mathbb{N}^\ast$ in that for all $\varepsilon > 0$ and all $\delta > 0$, there exists a sample size $n^\ast \in\mathbb{N}^\ast$ such for all $n \in \mathbb{N}^\ast$ satisfying $n \geq n^\ast$ we have:
	\begin{equation*}
	    \Pr \left(|| \hat{\bm{\theta}}_{(j,n,H)} - \bm{\theta}_0 ||_2 \geq \varepsilon \right) \leq   \delta.
	\end{equation*}
	\end{Corollary}
	\vspace{0.25cm}
	The proof of this corollary is given in Appendix \ref{proof:coro:consist}. As discussed in Remark \ref{asym:framework}, this corollary (similarly to Proposition \ref{THM:consistency} given further on) implies that when $p \to c <\infty$ as $n \to \infty$ we can simply write:
	\begin{equation*}
	   \lim_{n \to \infty} \Pr \left(|| \hat{\bm{\theta}}_{(j,n,H)} - \bm{\theta}_0 ||_2 > \varepsilon \right) = 0.
	\end{equation*}
	However, in the case where $p \to \infty$ as $n \to \infty$, the statement is weaker in the sense that we cannot conclude that the limit exist but only that $ \Pr (|| \hat{\bm{\theta}}_{(j,n,H)} - \bm{\theta}_0 ||_2 > \varepsilon )$ is arbitrarily small.
	
	Moreover, Corollary \ref{coro:consist} shows that in the situation where $\alpha = \nicefrac{1}{2}$ and $\gamma = 2$, the estimator $\hat{\bm{\theta}}_{(j,n,H)}$ is consistent for $\bm{\theta}_0$ under the conditions of Proposition \ref{THM:bias} and if
	\begin{equation*}
	    \lim_{n \to \infty} \; \frac{p^5}{n} = 0.
	\end{equation*}
	Nevertheless, this requirement can (in some sense) be relaxed as done in Proposition \ref{THM:consistency} which is based on Assumptions \ref{assum:A:3} and \ref{assum:D:4} (given below), and whose proof is provided in Appendix~\ref{proof:thm:2}.
	
	\setcounter{Assumption}{0}
	\renewcommand{\theHAssumption}{otherAssumption\theAssumption}
	\renewcommand\theAssumption{\Alph{Assumption}$_3$}
	\begin{Assumption}
	\label{assum:A:3}
	Preserving the same definitions and requirements given in Assumption \ref{assum:A:1}, we add the condition that $\bm{\Theta}$ is closed.
	\end{Assumption}
	\vspace{0.25cm}
    
	\setcounter{Assumption}{3}
	\renewcommand{\theHAssumption}{otherAssumption\theAssumption}
	\renewcommand\theAssumption{\Alph{Assumption}$_4$}
	\begin{Assumption}
	\label{assum:D:4}
	Preserving the same definitions given in Assumption \ref{assum:D:1} and defining $c_n \equiv \displaystyle{\max_{j=1,\dots,p}} \mathbf{c}_j(n)$, we require that:
	\begin{enumerate}
	    \item $\mathbf{a}(\btheta)$ is continuous and such that the function $\btheta + \mathbf{a}(\btheta)$ is injective,
	    
	    \item the constant $\beta$ and the sequence $\left\{c_n\right\}_{n\in\mathbb{N}^*}$ are such that
	\begin{equation*}
	   \lim_{n \to \infty} \; \frac{p^{\nicefrac{3}{2}}}{n^\beta} = 0, \;\;\;\;\text{and} \;\;\;\; \lim_{n \to \infty} \; p^{\nicefrac{1}{2}}c_n = 0.
	\end{equation*}
	\end{enumerate}
    The other requirements of Assumption \ref{assum:D:1} remain unchanged.
	\end{Assumption}
	\vspace{0.25cm}
	
	As stated in the above Assumption \ref{assum:D:4}, some requirements are the same as those in Assumption \ref{assum:D:1}, namely that there exist real $\beta, \gamma > 0$ such that for all $\bm{\theta} \in \bm{\Theta}$, we have, for all $k,l = 1,\,\ldots, \, p$,
		\begin{equation*}
		\begin{aligned}
			  \mathbf{L}_{k,l}(n) = \mathcal{O}(n^{-\beta}), \;\;\;  \mathbf{r}_k\left(\bm{\theta}, n\right) = \mathcal{O}(n^{-\gamma}) \;\;\;\ \text{and} \;\;\;
			  \lim_{n \to \infty} \; \frac{p^{\nicefrac{1}{2}}}{n^\gamma} = 0.
		 \end{aligned}
		\end{equation*}

	It can be noticed that Assumption \ref{assum:A:3} ensures that $\bm{\Theta}$ is compact while, compared to Assumption \ref{assum:D:1}, Assumption \ref{assum:D:4} firstly relaxes the condition on $\mathbf{a}(\btheta)$, then imposes stronger requirements on $\beta$ and finally imposes a condition on the vector $\mathbf{c}(n)$. Clearly, Assumption \ref{assum:A:1} implies Assumption \ref{assum:A:3} while, on the contrary, Assumptions \ref{assum:D:1} and \ref{assum:D:4} don't mutually imply each other. In the situation where $\alpha = \nicefrac{1}{2}$, $\beta = 1$, $\gamma = 2$ and $\mathbf{c}(n) = \mathbf{0}$, Assumption \ref{assum:D:4} (on which Proposition \ref{THM:consistency} is based) is satisfied if
	\begin{equation*}
	    \lim_{n \to \infty} \; \frac{p^{\nicefrac{3}{2}}}{n} = 0,
	\end{equation*}
	which relaxes (in some sense) the condition required in Corollary \ref{coro:consist}. 

	\begin{Proposition}
	\label{THM:consistency}
	Under Assumptions \ref{assum:A:3}, \ref{assum:B:1}, \ref{assum:C:1} and \ref{assum:D:4}, $\hat{\btheta}_{(j,n,H)}$ is a consistent estimator of $\btheta_0$ for all $(j,H) \in \mathbb{N} \times \mathbb{N}^\ast$ in the same sense as in Corollary \ref{coro:consist}. 
	\end{Proposition}
	\vspace{0.25cm}
	Having discussed the consistency of the JIE, the next section tackles the asymptotic normality of this estimator.

\subsection{Asymptotic Normality of the JIE}
\label{sec:norm}

In order to study the asymptotic normality of the JIE we need to again modify the assumption framework one final time. Indeed, to derive this result, some of the assumptions are modified in order to incorporate additional requirements on the auxiliary estimator. 

Before introducing the next assumption we define various quantities. First, we let
 \begin{equation*}
    \bm{\Sigma}(\btheta,n) \equiv \var \left( \sqrt{n} \mathbf{v}(\bm{\theta}, n, \bm{\omega}) \right),
\end{equation*}
where $\bm{\Sigma}(\btheta,n)$ is nonsingular.
Then, using this definition, we introduce the following quantity
\begin{equation*}
    Y(\btheta, n, \bm{\omega}, \mathbf{u}) \equiv \sqrt{n} \mathbf{u}^T \bm{\Sigma}(\btheta,n)^{-\nicefrac{1}{2}} \mathbf{v}(\bm{\theta}, n, \bm{\omega}),
\end{equation*}
where $\mathbf{u} \in \real^p$ is such that $||\mathbf{u} ||_2 = 1$. Clearly, from the definition of $\mathbf{v}(\bm{\theta}, n, \bm{\omega})$ we have that $\mathbb{E}[Y(\btheta, n, \bm{\omega}, \mathbf{u})] = 0$. Therefore, without loss of generality we can always decompose $Y(\btheta, n, \bm{\omega}, \mathbf{u})$ as follows:
\begin{equation}
    Y(\btheta, n, \bm{\omega}, \mathbf{u}) =  Z(\btheta, n, \bm{\omega}, \mathbf{u}) + \delta_n W(\btheta, n, \bm{\omega}, \mathbf{u}),
    \label{eq:gauss:approx}
\end{equation}
where $W(\btheta, n, \bm{\omega}, \mathbf{u})$ and $Z(\btheta, n, \bm{\omega}, \mathbf{u})$ are zero mean random variables and $\big(\delta_n\big) \in \mathfrak{D}$ where
\begin{equation*}
	    \mathfrak{D} \equiv \left\{ \left( \delta_n \right)_{n \in \mathbb{N}^\ast} \, \big| \; \delta_n \in \real \setminus (-\infty, 0) \; \forall n, \; \delta_l \geq \delta_k \;\; \text{if} \;\; l < k \;\; \text{and} \;\; \lim_{n \to \infty} \delta_n = 0\right\}.
\end{equation*}
In the assumption below we restrict the behaviour of $\mathbf{v}(\bm{\theta}, n, \bm{\omega})$ and in particular we require that it satisfies a specific Gaussian approximation based on the decomposition considered in (\ref{eq:gauss:approx}).

\setcounter{Assumption}{2}
\renewcommand{\theHAssumption}{otherAssumption\theAssumption}
\renewcommand\theAssumption{\Alph{Assumption}$_3$}

\begin{Assumption}
\label{assum:C:3}
Preserving the requirements of Assumption \ref{assum:C:1} we add the following conditions on $\mathbf{v}(\bm{\theta}, n, \bm{\omega})$. With $n \in \mathbb{N}^*$, there exists a sample size $n^\ast \in \mathbb{N}^*$ such that:
\begin{enumerate}
\item For all $\bm{\theta}\in \bm{\Theta}$ and all $n \geq n^\ast$, the matrix $\bm{\Sigma}(\btheta,n)$ exists, is nonsingular and is such that $\bm{\Sigma}_{k,l}(\btheta,n) = \mathcal{O}(1)$ for all $k,l=1,\dots,p$.
\item For all $\left(\btheta,\bm{\omega}\right)\in\bm{\Theta}\times\bm{\Omega}$ and all $n \geq n^\ast$, the Jacobian matrix $\mathbf{V}(\bm{\theta}, n , \bm{\omega}) \equiv \frac{\partial}{\partial \, \btheta^T} \mathbf{v}\left(\bm{\theta}, n, \bm{\omega}\right)$ exists and is continuous in $\btheta \in \bm{\Theta}$. 
\item  Considering the decomposition in (\ref{eq:gauss:approx}), then for all $\mathbf{u} \in \real^p$ such that $||\mathbf{u} ||_2 = 1$, there exist sequences $\big(\delta_n\big) \in \mathfrak{D}$, as well as a random variable $W(\btheta, n, \bm{\omega}, \mathbf{u}) = \mathcal{O}_{\rm p}(1)$ with existing second moment, such that
$Z(\btheta, n, \bm{\omega}, \mathbf{u})$ is a standard normal random variable.
\end{enumerate}
\end{Assumption}

The first requirement of Assumption \ref{assum:C:3} is quite mild and commonly assumed as it simply requires that $\bm{\Sigma}(\btheta, n)$ is a suitable covariance matrix. {\color{black}In addition, it implies that $\mathbf{v}(\bm{\theta}, n, \bm{\omega}) = \mathcal{O}_{\rm p} \left(n^{-\nicefrac{1}{2}}\right)$, elementwise}. Similarly to the discussion following Assumption \ref{assum:D:2}, our second requirement {\color{black} ensures that $\mathbf{V}_{k,l}(\hbt, n , \bm{\omega})$ is a bounded random variable for all $k,l=1,\dots,p$ and all $(n, \bm{\omega}) \in \mathbb{N}^* \times \bm{\Omega}$.} 
The third condition of Assumption \ref{assum:C:3} describes how ``close'' $\mathbf{v}(\bm{\theta}, n, \bm{\omega})$ is to a multivariate normal distribution. Such an assumption is quite strong and may not always be satisfied. In the case where $p \to c < \infty$ as $n \to \infty$, our requirement on the distribution of $\mathbf{v}(\bm{\theta}, n, \bm{\omega})$ can simply be expressed as
\begin{equation*}
    \sqrt{n} \mathbf{v}(\bm{\theta}, n, \bm{\omega}) \xrightarrow{d} \mathcal{N}(\mathbf{0}, \bm{\Sigma}(\btheta)),
\end{equation*}
where $\bm{\Sigma}(\btheta) \equiv \displaystyle \lim_{n \to \infty} \bm{\Sigma}(\btheta, n)$. Finally, Assumption \ref{assum:C:3} clearly implies Assumption \ref{assum:C:1} but not necessarily Assumption~\ref{assum:C:2}.

In the next assumption, we impose an additional requirement on the bias of the auxiliary estimator and for this purpose we introduce the following notation. Let $\mathbf{f}(\btheta) \,  : \bm{\Theta} \to \real^p$ and $\mathbf{F}(\btheta) \equiv \frac{\partial}{\partial \, \btheta^T} \mathbf{f}\left(\bm{\theta}\right) \in \real^{p \times p}$ be its Jacobian matrix. Then, we define $\mathbf{F}\big(\btheta^{(\mathbf{f})}\big)$ such that when using the multivariate mean value theorem between $\hat{\btheta}_{(j,n,H)}$ and ${\btheta}_{0}$ we obtain
\begin{equation}
    \mathbf{f} \left(\hat{\btheta}_{(j,n,H)}\right) =  \mathbf{f} \left({\btheta}_{0}\right) + \mathbf{F}\left(\btheta^{(\mathbf{f})}\right) \left( \hat{\btheta}_{(j,n,H)} - {\btheta}_{0} \right).
    \label{def:MVT:multi_vari}
\end{equation}
Therefore, $\btheta^{(\mathbf{f})}$ corresponds to a set of $p$ vectors lying in the segment $(1 - \lambda) \hat{\btheta}_{(j,n,H)} + \lambda {\btheta}_{0}$ for $\lambda \in [0,1]$ (with respect to the function $\mathbf{f}$). Keeping the latter notation in mind (for a generic function $\mathbf{f}$), we provide the following assumption.

\setcounter{Assumption}{3}
\renewcommand{\theHAssumption}{otherAssumption\theAssumption}
\renewcommand\theAssumption{\Alph{Assumption}$_5$}
\begin{Assumption}
\label{assum:D:5}
There exists a sample size $n^\ast \in \mathbb{N}^\ast$ such that for all $n \in \mathbb{N}^\ast$ satisfying $n \geq n^\ast$ and all $\btheta\in\bm{\Theta}$, the Jacobian matrices $\mathbf{A}(\bm{\theta}) \equiv \frac{\partial}{\partial \, \btheta^T} \mathbf{a}\left(\bm{\theta}\right)$ and $\mathbf{R}(\bm{\theta}, n) \equiv \frac{\partial}{\partial \, \btheta^T} \mathbf{r}\left(\bm{\theta}, n\right)$ exist and are continuous. Moreover, we require that for all $\mathbf{u}\in\real^p$ with $||\mathbf{u}||_2 = 1$, there exists a sequence $(\delta_n) \in\mathfrak{D}$ such that for all $(H, \bm{\omega}) \in \mathbb{N}^\ast \times \bm{\Omega}$
\begin{equation}
    \sqrt{n}\left(1+\frac{1}{H}\right)^{-\nicefrac{1}{2}}\mathbf{u}^T\bm{\Sigma}(\bto,n)^{-\nicefrac{1}{2}}\Big[\mathbf{B}(\bto,n)-\mathbf{B}\left(\btheta^{(\mathbf{a}, \mathbf{r})},n\right)\Big]\left(\hbt-\bto\right) = \mathcal{O}_{\rm p}(\delta_n),
    \label{eq:assum:d:5}
\end{equation}
where $\mathbf{B}(\bt,n) \equiv \mathbf{I} + \mathbf{A}(\bt) + \mathbf{L}(n) + \mathbf{R}(\bt,n)$ and $\mathbf{B}\big(\btheta^{(\mathbf{a}, \mathbf{r})},n\big) \equiv \mathbf{I} + \mathbf{A}(\btheta^{(\mathbf{a})}) + \mathbf{L}(n) + \mathbf{R}\left(\btheta^{(\mathbf{r})},n\right)$.
\end{Assumption}
\vspace{0.25cm}

Assumption \ref{assum:D:5} allows us to quantify how ``far'' the matrices $\mathbf{B}(\bto,n)$ and $\mathbf{B}\big(\hbt, n\big)$ are from each other. In the case where $p \to c < \infty$ as $n \to \infty$, the consistency of $\hbt$ and the continuity of $\mathbf{B}(\bt) \equiv \displaystyle\lim_{n \to \infty} \mathbf{B}(\bt,n)$ in $\bm{\theta}_0$ would be sufficient so that Assumption \ref{assum:D:5} would not be needed to establish the asymptotic distribution of $\hbt$. However, when $p \to \infty$ as $n \to \infty$ the requirement in (\ref{eq:assum:d:5}) can be strong and difficult to verify for a specific model and auxiliary estimator. Having stated this, the following proposition defines the distribution of the JIE estimator.

\begin{Proposition}
\label{Thm:Gauss:approx}
Under Assumptions \ref{assum:A:2}, \ref{assum:B:1}, \ref{assum:C:3} and \ref{assum:D:5}, for all $\mathbf{u} \in \real^p$ such that $||\mathbf{u}||_2 = 1$, there exist a sample size $n^\ast \in \mathbb{N}^\ast$ and a sequence $\left( \delta_n\right)\in\mathfrak{D}$ such that for all $n \in \mathbb{N}^\ast$ satisfying $n \geq n^\ast$ we have
\begin{equation*}
    \sqrt{n}\left(1+\frac{1}{H}\right)^{-\nicefrac{1}{2}}\mathbf{u}^T\bm{\Sigma}(\bto,n)^{-\nicefrac{1}{2}}\mathbf{B}(\bto, n)\left(\hbt-\bto\right) \overset{d}{=} Z + \delta_n\mathcal{O}_{\rm p}\left( \max\left(1,\frac{p^2}{\sqrt{H}}\right)\right), 
\end{equation*}
where $Z \sim \mathcal{N}(0,1)$.
\end{Proposition}
\vspace{0.25cm}

The proof of this result is given in Appendix \ref{proof:Gauss:approx}. Proposition \ref{Thm:Gauss:approx} provides an approximation of the distribution of $\hbt$ that can be used in the following way. As discussed in Remark \ref{asym:framework}, when $p \to c < \infty$ as $n \to \infty$, Proposition \ref{Thm:Gauss:approx} simply states that
\begin{equation*}
    \sqrt{n} \left(\hbt-\bto\right) \xrightarrow{d} \mathcal{N}\left(\mathbf{0}, \bm{\Xi}_H\right),
\end{equation*}
where
\begin{equation}
    \bm{\Xi}_H \equiv \left(1+\frac{1}{H}\right) \mathbf{B}(\bto)^{-1} \bm{\Sigma}(\bto) \mathbf{B}(\bto)^{-1},
    \label{Eq_ass-var-JIE}
\end{equation}
provided that $\mathbf{B}(\bm{\theta})$ is continuous in $\bto$ and $\mathbf{B}(\bm{\theta}_0)$ is non-singular. The estimation of $\bm{\Xi}_H$  is discussed, for example, in \cite{gourieroux1993indirect} (see also \citealp{genton2000robust}). When $p \to \infty$ as $n \to \infty$, one needs to additionally assume that $H=\mathcal{O}(p^4)$ in order for the Gaussian approximation in Proposition \ref{Thm:Gauss:approx} to hold. This therefore suggests that the ``quality'' (and validity) of the approximation depends on $H$ when $p$ is large. However, the conditions of Proposition \ref{Thm:Gauss:approx} are sufficient but may not be necessary thereby implying that  $H=\mathcal{O}(p^4)$ may not always be needed as a condition for this approximation.

\begin{Remark}
In the proof of Proposition \ref{Thm:Gauss:approx} we show that 
\begin{equation}
    \sqrt{n}\left(1+\frac{1}{H}\right)^{-\nicefrac{1}{2}}\mathbf{u}^T\bm{\Sigma}(\bto,n)^{-\nicefrac{1}{2}}\mathbf{B}(\bt^{(\mathbf{a}, \mathbf{r})}, n)\left(\hbt-\bto\right) \overset{d}{=} Z + \delta_n\mathcal{O}_{\rm p}\left( \max\left(1,\frac{p^2}{\sqrt{H}}\right)\right). 
    \label{eq:gauss:approx:sm}
\end{equation}
By combining this result with Assumption \ref{assum:D:5} the statement of the proposition directly follows. As a result, Assumption \ref{assum:D:5} is not strictly necessary for the Gaussian approximation of the distribution of $\hbt$. However, while this assumption is strong, it is quite convenient to deliver a reasonable approximation of $\mathbf{B}(\bt^{(\mathbf{a}, \mathbf{r})}, n)$ when $p$ is allowed to diverge, since the applicability of (\ref{eq:gauss:approx:sm}) is quite limited in practice.
\end{Remark}
\vspace{0.25cm}


\subsection{Main Results}
\label{sec:discussion}

\begin{figure}[!ht]
\centering
\begin{adjustbox}{width=.85\textwidth}

\begin{tikzpicture}

\node (rectD14) at (4.625,-2.75) [fill = blue!10, draw = blue!45, very thick,  minimum width=0.8cm, minimum height=0.8cm, rounded corners=0.2cm] {};

\node (rectD15) at (6.875,-2.75) [fill = orange!10, draw = orange!70, very thick, minimum width=0.8cm, minimum height=0.8cm, rounded corners=0.2cm] {};

\node (rectA13) at (4.625,2.75) [fill = blue!10, draw = blue!45, very thick,  minimum width=0.8cm, minimum height=0.8cm, rounded corners=0.2cm] {};

\node (rectD12) at (2.375,-2.75) [fill = green!20, draw = green!80!black!80, very thick, minimum width=0.8cm, minimum height=0.8cm, rounded corners=0.2cm] {};

\node (rect1) at (-1,0) [myblue, draw, thick, fill = mygrey, minimum width=1cm, minimum height=4cm] {};
\node (rect2) at (1.25,0) [myblue, draw, thick, fill = mygrey, minimum width=1cm, minimum height=4cm] {};
\node (rect3) at (3.5,0) [myblue, draw, fill = mygrey, thick, minimum width=1cm, minimum height=4cm] {};
\node (rect4) at (5.75,0) [myblue, draw, fill = mygrey, thick, minimum width=1cm, minimum height=4cm] {};
\node (rect5) at (8,0) [myblue, draw, fill = mygrey, thick, minimum width=1cm, minimum height=4cm] {};

\node (rectA2) at (1.25,1.5) [fill = green!20, draw = green!80!black!80, very thick, minimum width=0.8cm, minimum height=0.8cm, rounded corners=0.2cm] {};

\node [myblue, text width=1.8cm,align=center] (thm2) at (1.25,4.1) {\scriptsize Proposition \ref{THM:bias}\\Unbiasedness};
\node [myblue, text width=1.8cm,align=center] (coro1) at (-1,4.1) {\scriptsize Corollary \ref{coro:consist}\\Consistency};
\node [myblue, text width=2.2cm,align=center] (thm1) at (3.5,4.1) {\scriptsize Proposition \ref{thm:iter:boot}\\IB convergence};
\node [myblue, text width=1.8cm,align=center] (thm3) at (5.75,4.1) {\scriptsize Proposition \ref{THM:consistency}\\Consistency};
\node [myblue, text width=2.2cm,align=center] (thm4) at (8,4.1) {\scriptsize Proposition \ref{Thm:Gauss:approx}\\Asym. normality};

\draw [myblue, -] (coro1) -- (rect1);
\draw [myblue, -] (thm2) -- (rect2);
\draw [myblue, -] (thm1) -- (rect3);
\draw [myblue, -] (thm3) -- (rect4);
\draw [myblue, -] (thm4) -- (rect5);


\node [aqua] (D) at (2.375,-2.75) {\ref{assum:D:12}};

\node [mypurple] (Dast) at (4.625,-2.75) {\ref{assum:D:14}};

\node [mypurple] (Dstar) at (6.875,-2.75) {\ref{assum:D:15}};

\node [mypurple] (Astar12) at (4.625,2.75) {\ref{assum:A:13}};



\node (rectC3) at (8,-0.5) [fill = orange!10, draw = orange!70, very thick,  minimum width=0.7cm, minimum height=0.7cm, rounded corners=0.2cm] {};

\node (rectB1) at (8,0.5) [fill = orange!10, draw = orange!70, very thick,  minimum width=0.7cm, minimum height=0.7cm, rounded corners=0.2cm] {};

\node (rectC12) at (5.75,-0.5) [fill = blue!10, draw = blue!45, very thick,  minimum width=0.7cm, minimum height=0.7cm, rounded corners=0.2cm] {};

\node (rectB2) at (5.75,0.5) [fill = blue!10, draw = blue!45, very thick, minimum width=0.7cm, minimum height=0.7cm, rounded corners=0.2cm] {};

\node (rectA4) at (8,1.5) [fill = orange!10, draw = orange!70, very thick,  minimum width=0.7cm, minimum height=0.7cm, rounded corners=0.2cm] {};

\node (rectC1) at (3.5,-0.5) [fill = green!20, draw = green!80!black!80, very thick, minimum width=0.7cm, minimum height=0.7cm, rounded corners=0.2cm] {};

\node (rectB2) at (3.5,0.5) [fill = green!20, draw = green!80!black!80, very thick, minimum width=0.7cm, minimum height=0.7cm, rounded corners=0.2cm] {};

\node [myorange] (A2) at (1.25,1.5) {\ref{assum:A:2}};
\node [myorange] (B1) at (1.25,0.5) {\ref{assum:B:1}};
\node [myorange] (D2) at (1.25,-1.5) {\ref{assum:D:2}};

\node [myorange] (A22) at (-1,1.5) {\ref{assum:A:2}};
\node [myorange] (B2) at (-1,0.5) {\ref{assum:B:1}};
\node [myorange] (C2) at (-1,-0.5) {\ref{assum:C:2}};
\node [myorange] (D3) at (-1,-1.5) {\ref{assum:D:3}};

\node [myorange] (A1) at (3.5,1.5) {\ref{assum:A:1}};
\node [myorange] (B3) at (3.5,0.5) {\ref{assum:B:1}};
\node [myorange] (C1) at (3.5,-0.5) {\ref{assum:C:1}};
\node [myorange] (D1) at (3.5,-1.5) {\ref{assum:D:1}};

\node [myorange] (A3) at (5.75,1.5) {\ref{assum:A:3}};
\node [myorange] (B4) at (5.75,0.5) {\ref{assum:B:1}};
\node [myorange] (C12) at (5.75,-0.5) {\ref{assum:C:1}};
\node [myorange] (D4) at (5.75,-1.5) {\ref{assum:D:4}};

\node [myorange] (A4) at (8,1.5) {\ref{assum:A:2}};
\node [myorange] (B5) at (8,0.5) {\ref{assum:B:1}};
\node [myorange] (C3) at (8,-0.5) {\ref{assum:C:3}};
\node [myorange] (D5) at (8,-1.5) {\ref{assum:D:5}};

\draw [myorange, ->] (D3) -- (D2);
\draw [myorange, ->] (rectA2) -- (A1);
\draw [myorange, ->] (C2) -- (rectC1);
\draw [myorange, ->] (A3) -- (A1);
\draw [myorange, ->] (rectC3) -- (rectC12);
\draw [myorange, ->] (D5) -- (D4);
\draw [myorange, ->] (rectA4) edge[bend left = -15] (rectA13);

\draw [mypurple, ->] (rectA13.south west) -- (A1);
\draw [mypurple, ->] (rectA13.south east) -- (A3);

\draw [mypurple, ->] (rectD12.north west) -- (D2);
\draw [mypurple, ->] (rectD12.north east) -- (D1);
\draw [mypurple, ->] (rectD14.north east) -- (D4);
\draw [mypurple, ->] (rectD14.north west) -- (D1);
\draw [mypurple, ->] (rectD15.north west) -- (D1);
\draw [mypurple, ->] (rectD15.north east) -- (D5);

\draw [myblue, double, <-] (6.3,0.15) -- (7.45,0.15);
\draw [myblue, double, <-] (0.7,0.15) -- (-0.45,0.15);

\node (rectred) at (-1,-4) [fill = green!20, draw = green!80!black!80, thick, minimum width=0.4cm, minimum height=0.4cm, rounded corners=0.1cm] {};
\node (rectblue) at (-1,-4.5) [fill = blue!10, draw = blue!45, thick,  minimum width=0.4cm, minimum height=0.4cm, rounded corners=0.1cm] {};
\node (rectgreen) at (-1,-5) [fill = orange!10, draw = orange!70, thick, minimum width=0.4cm, minimum height=0.4cm, rounded corners=0.1cm] {};

\node [myblue, text width=8.8cm] (thm2) at (3.8,-4) {\scriptsize Assumptions used in Theorem \ref{thm:main:1} (IB - Unbiasedness)};

\node [myblue, text width=8.8cm] (thm2) at (3.8,-4.5) {\scriptsize Assumptions used in Theorem \ref{thm:main:2} (IB - Consistency)};

\node [myblue, text width=8.8cm] (thm2) at (3.8,-5) {\scriptsize Assumptions used in Theorem \ref{thm:main:3} (IB - Asymptotic normality)};

\end{tikzpicture}
\end{adjustbox}
\caption{Illustration of the implication links (arrows) of the different assumptions used in Sections \ref{sec:boot} to \ref{sec:norm} and construction of the sufficient conditions 
\ref{assum:A:13}, \ref{assum:D:12}, \ref{assum:D:14} and \ref{assum:D:15},  needed for Theorems \ref{thm:main:1} to \ref{thm:main:3}. The double arrows are used to denote implication of boxes (i.e. a set of assumptions) while simple arrows are used to represent implications between assumptions. The color boxes are used to highlight the assumptions used in Theorems \ref{thm:main:1} to \ref{thm:main:3}.}
\label{Fig_assumptions-mess}
\end{figure}
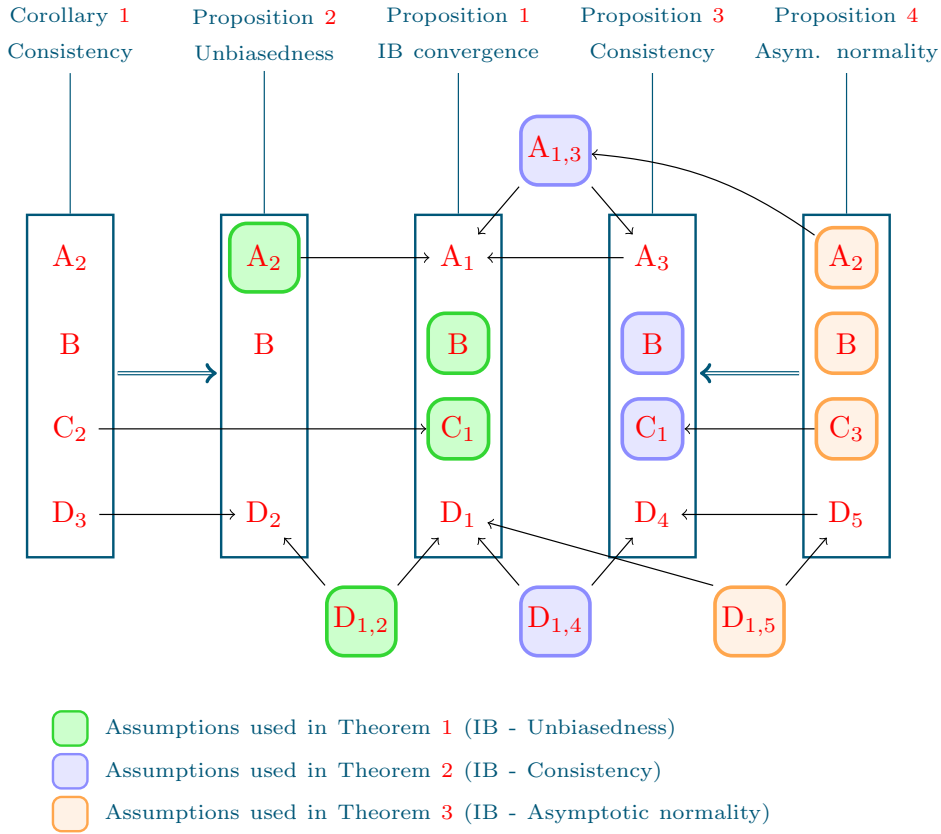

The results presented so far can now be combined to deliver the three main theorems of this paper. These refer to the statistical properties of $\tilde{\bm{\theta}}_{(j,n,H)}$ defined in (\ref{eq:iterboot}) (i.e. the limit in $k$ of $\tilde{\bm{\theta}}^{(k)}_{(j,n,H)}$), namely its unbiasedness (in finite samples), its consistency and its asymptotic normality. In order to state these properties we formulate a new set of assumptions which simply consists in a combination of those given in the previous sections. More precisely, we construct one assumption on the topology of $\bm{\Theta}$ and three assumptions on the form of the bias of the auxiliary estimator. Given the amount of assumptions and the different assumption frameworks, as a support to the reader, Figure \ref{Fig_assumptions-mess} provides an overview of the different frameworks that include the new assumptions (delivered further on) and their different implications and combination structure. For each new framework structure we state its relative theorem whose proofs are omitted given that they are a direct consequence of the results in the previous sections. 

Based on the above premise, we consider the new assumption needed to derive the unbiasedness of $\tilde{\btheta}_{(j,n,H)}$. This assumption concerns the form of the bias of the auxiliary estimator. Before introducing Assumption \ref{assum:D:12}, we briefly recall that the bias function $\mathbf{d}(\btheta,n)$ can be decomposed as $\mathbf{d}(\btheta,n) = \mathbf{a}(\btheta) + \mathbf{c}(n) +  \mathbf{L}(n)\btheta + \mathbf{r}(\btheta,n),$
where $\mathbf{a}(\btheta)$, $\mathbf{c}(n)$ and $\mathbf{r}(\btheta,n) \in \real^p$ and $ \mathbf{L}(n) \in \real^{p \times p}$.

%

\setcounter{Assumption}{3}
\renewcommand\theAssumption{\Alph{Assumption}$_{1,2}$}
\begin{Assumption}
\label{assum:D:12}
	The bias function $\mathbf{d}\left(\bm{\theta}, n\right)$ is such that:
		\begin{enumerate}
		    \item There exists a matrix $\mathbf{M} \in \real^{p \times p}$ with $|| \mathbf{M} ||_F < 1$ and a vector $\mathbf{s} \in \real^{p}$ such that \begin{equation*}
		        \mathbf{a}(\btheta) = \mathbf{M}\btheta + \mathbf{s}.
		    \end{equation*}
		    
		    \item The Jabcobian matrix $\mathbf{R}(\bm{\theta}, n) \equiv \frac{\partial}{\partial \, \btheta^T} \mathbf{r}\left(\bm{\theta}, n\right)$ exists for all $(\bm{\theta}, n) \in \bm{\Theta} \times \mathbb{N}^*$ and $\mathbf{R}_{k,l}\left(\bm{\theta}, n\right)$ is continuous in $\btheta \in \bm{\Theta}$ for any $k,l=1,\dots,p$. Moreover, there exist real $\beta, \gamma > 0$ such that for all $\bm{\theta} \in \bm{\Theta}$ and any $k,l=1,\dots,p$, we have
		    \begin{equation*}
		    \begin{aligned}
			  \mathbf{L}_{k,l}(n) = \mathcal{O}(n^{-\beta}), \;\;\;  \mathbf{r}_k\left(\bm{\theta}, n\right) = \mathcal{O}(n^{-\gamma}),\;\;\; 
		      \text{and} \;\;\;
			  \lim_{n \to \infty} \; \frac{p}{n^\beta} = 0.
		    \end{aligned}
		    \end{equation*}
		  
		    \item There exists a sample size $n^*\in \mathbb{N}^*$ such that for all $n \in \mathbb{N}^*$ satisfying $n \geq n^*$
		    \begin{itemize}
                \item the matrix $(\mathbf{M} + \mathbf{L}(n)+\mathbf{I})^{-1}$ exists,
                 \item $p^2 \,n^{-\gamma} < 1$.
		    \end{itemize}
		\end{enumerate}
\end{Assumption}

With the assumptions we can now deliver Theorem \ref{thm:main:1} that states the unbiasedness of  $\tilde{\btheta}_{(j,n,H)}$.

\begin{Theorem}
\label{thm:main:1}
 Under Assumptions \ref{assum:A:2}, \ref{assum:B:1}, \ref{assum:C:1}, \ref{assum:D:12}, $\tilde{\btheta}_{(j,n,H)}$ is an unbiased estimator of $\bm{\theta}_0$ in that for all $n\in\mathbb{N}^*$ such that $n\geq n^*$ and for all $(j, H) \in \mathbb{N} \times \mathbb{N}^\ast$, we have 
		\begin{equation*}
		\Big\lVert	\mathbb{E} \left[\tilde{\bm{\theta}}_{(j,n,H)} \right] - \bto  \Big\rVert_2 = 0. 
		\end{equation*}
\end{Theorem}

We now focus on the consistency of $\tilde{\btheta}_{(j,n,H)}$ and, for this purpose, we state the following assumptions.

\setcounter{Assumption}{0}
\renewcommand{\theHAssumption}{otherAssumption\theAssumption}
\renewcommand\theAssumption{\Alph{Assumption}$_{1,3}$}
\begin{Assumption}
\label{assum:A:13}
Let $\bm{\Theta}$ be compact and such that
		\begin{equation*}
			\argzero_{\bm{\theta} \in
		      \real^p \setminus \bm{\Theta}} \;  \hat{\bm{\pi}}(\bm{\theta}_0, n, \bm{\omega}_0) - \frac{1}{H} \sum_{h = 1}^H  \hat{\bm{\pi}}(\bm{\theta}, n, \bm{\omega}_{h+jH})
		      = \emptyset\,.
		\end{equation*}
\end{Assumption}

\setcounter{Assumption}{3}
\renewcommand\theAssumption{\Alph{Assumption}$_{1,4}$}
\begin{Assumption}
\label{assum:D:14}
    The bias function $\mathbf{d}\left(\bm{\theta}, n\right)$ is such that:
		\begin{enumerate}
		    \item The function $\mathbf{a}(\btheta)$ is continuous and such that the function $\btheta + \mathbf{a}(\btheta)$ is injective.
		    
		    \item There exist real $\beta, \gamma > 0$ such that for all $\bm{\theta} \in \bm{\Theta}$ and any $k,l=1,\dots,p$, we have
		    \begin{equation*}
		    \begin{aligned}
			  \mathbf{L}_{k,l}(n) = \mathcal{O}(n^{-\beta}), \;\;\;  \mathbf{r}_k\left(\bm{\theta}, n\right) = \mathcal{O}(n^{-\gamma})\;\;\; 
		      \text{and} \;\;\;
			  \lim_{n \to \infty} \; \frac{p}{n^\beta} = 0.
		    \end{aligned}
		    \end{equation*}
		  
		    \item Defining $c_n \equiv \displaystyle{\max_{j=1,\dots,p}} \mathbf{c}_j(n)$ for all $n \in \mathbb{N}^*$, we require that  the constant $\beta$ and the sequence $\left\{c_n\right\}_{n\in\mathbb{N}^*}$ are such that
            \begin{equation*}
	        \lim_{n \to \infty} \; \frac{p^{\nicefrac{3}{2}}}{n^\beta} = 0, \;\;\;\;\text{and} \;\;\;\; \lim_{n \to \infty} \; p^{\nicefrac{1}{2}}c_n = 0.
	        \end{equation*}
		\end{enumerate}
\end{Assumption}

Again, using these new assumptions we can deliver the next theorem that states the consistency properties of $\tilde{\btheta}_{(j,n,H)}$.

\begin{Theorem}
\label{thm:main:2}
	Under Assumptions \ref{assum:A:13}, \ref{assum:B:1}, \ref{assum:C:1}, and \ref{assum:D:14}, $\tilde{\btheta}_{(j,n,H)}$ is a consistent estimator of $\bm{\theta}_0$ in that for all $(j, H) \in \mathbb{N} \times \mathbb{N}^\ast$, for all $\varepsilon > 0$ and all $\delta > 0$ there exists a sample size $n^\ast \in\mathbb{N}^\ast$ such that for all $n \in \mathbb{N}^\ast$ satisfying $n \geq n^\ast$:
	\begin{equation*}
	    \Pr \left(\big\lVert \tilde{\btheta}_{(j,n,H)} - \bm{\theta}_0 \big\rVert_2 > \varepsilon \right) \leq   \delta.
	\end{equation*}
\end{Theorem}

Finally we For the asymptotic normality of $\tilde{\btheta}_{(j,n,H)}$, we state the following assumptions.

\setcounter{Assumption}{3}
\renewcommand\theAssumption{\Alph{Assumption}$_{1,5}$}
\begin{Assumption}
\label{assum:D:15}
    The bias function $\mathbf{d}\left(\bm{\theta}, n\right)$ is such that:
		\begin{enumerate}
		    \item The function $\mathbf{a}(\bm{\theta})$ is a contraction map in that for any $\bm{\theta}_1,\bm{\theta}_2\in\bm{\Theta}$ such that $\bm{\theta}_1 \neq \bm{\theta}_2$ we have 
		    \begin{equation*}
		        \big\lVert \mathbf{a}(\bm{\theta}_2)-\mathbf{a}(\bm{\theta}_1)  \big\rVert_2 <  \big\rVert \bm{\theta}_2-\bm{\theta}_1 \big\lVert_2\,.
		    \end{equation*}
		    
		\item There exist real $\beta, \gamma > 0$ such that for all $\bm{\theta} \in \bm{\Theta}$ and any $k,l = 1,\,\ldots, \, p$, we have
		\begin{equation*}
		\begin{aligned}
			  \mathbf{L}_{k,l}(n) = \mathcal{O}(n^{-\beta}), \;\;\;  \mathbf{r}_k\left(\bm{\theta}, n\right) = \mathcal{O}(n^{-\gamma}),\;\;\; 
			  \lim_{n \to \infty} \; \frac{p}{n^\beta} = 0,  \;\;\;\ \text{and} \;\;\;
			  \lim_{n \to \infty} \; \frac{p^{\nicefrac{1}{2}}}{n^\gamma} = 0.
		 \end{aligned}
		\end{equation*}
		%
		
       \item There exists a sample size $n^\ast \in \mathbb{N}^\ast$ such that the Jacobian matrices $\mathbf{A}(\bm{\theta}) \equiv \frac{\partial}{\partial \, \btheta^T} \mathbf{a}\left(\bm{\theta}\right)$ and $\mathbf{R}(\bm{\theta}, n) \equiv \frac{\partial}{\partial \, \btheta^T} \mathbf{r}\left(\bm{\theta}, n\right)$ exist and are continuous for all $n \in \mathbb{N}^*$ satisfying $n \geq n^*$ and all $\btheta\in\bm{\Theta}$. Moreover, we require that for all $\mathbf{u}\in\real^p$ with $||\mathbf{u}||_2 = 1$, there exists a sequence $(\delta_n) \in\mathfrak{D}$ such that for all $(H, \bm{\omega}) \in \mathbb{N}^\ast \times \bm{\Omega}$
        \begin{equation*}
         \sqrt{n}\left(1+\frac{1}{H}\right)^{-\nicefrac{1}{2}}\mathbf{u}^T\bm{\Sigma}(\bto,n)^{-\nicefrac{1}{2}}\Big[\mathbf{B}(\bto,n)-\mathbf{B}\left(\btheta^{(\mathbf{a}, \mathbf{r})},n\right)\Big]\left(\tilde{\btheta}_{(j,n,H)}-\bto\right) = \mathcal{O}_{\rm p}(\delta_n),
        \end{equation*}
        where $\mathbf{B}(\bt,n) \equiv \mathbf{A}(\bt) + \mathbf{R}(\bt,n)$ and $\mathbf{B}\big(\btheta^{(\mathbf{a}, \mathbf{r})},n\big) \equiv \mathbf{A}(\btheta^{(\mathbf{a})}) + \mathbf{R}\left(\btheta^{(\mathbf{r})},n\right)$.
    \end{enumerate}

\end{Assumption}

Theorem \ref{thm:main:3} below states the asymptotic normality of  $\tilde{\btheta}_{(j,n,H)}$.

\begin{Theorem}
\label{thm:main:3}
Under Assumptions \ref{assum:A:2}, \ref{assum:B:1}, \ref{assum:C:3} and \ref{assum:D:15}, for all $\mathbf{u} \in \real^p$ such that $||\mathbf{u}||_2 = 1$, there exist a sample size $n^\ast \in \mathbb{N}^\ast$ and a sequence $\left( \delta_n\right)\in\mathfrak{D}$ such that for all for all $n \in \mathbb{N}^*$ satisfying $n \geq n^*$ we have
\begin{equation*}
    \sqrt{n}\left(1+\frac{1}{H}\right)^{-\nicefrac{1}{2}}\mathbf{u}^T\bm{\Sigma}(\bto,n)^{-\nicefrac{1}{2}}\mathbf{B}(\bto, n)\left(\tilde{\btheta}_{(j,n,H)}-\bto\right) \overset{d}{=} Z + \delta_n\mathcal{O}_{\rm p}\left( \max\left(1,\frac{p^2}{\sqrt{H}}\right)\right), 
\end{equation*}
where $Z \sim\mathcal{N}(0,1)$.
\end{Theorem}

\section{Classical and Robust Bias-Corrected Estimators for the Logistic Regression Model}
\label{sec:logistic}

The aim of this section is to use the general framework developed in the previous sections to the case of parameter's estimation for the logistic regression model \citep{NeWe:72,McCuNe:89}, which is one of the most frequently used model for binary response variables conditionally on a set of predictors. However, it is well known that in some quite frequent practical situations, the MLE is biased or its computation can become very unstable, especially when performing some type of resampling scheme for inference. The underlying reasons are diverse, but the main ones are the possibly large $p/n$ ratio, separability (leading to regression slopes estimates of infinite value) and data contamination (robustness). 

The first two sources are often confounded and practical solutions are continuously sought to overcome the difficulty in performing ``reasonable'' inference. For example, in medical studies, the biased MLE together with the problem of separability, has led to a rule of thumb called the number of Events Per Variable (EPV), that is the number of occurrences of the least frequent event over the number of predictors, which is used in practice to choose the maximal number of predictors one is ``allowed'' to use in a logistic regression model (see e.g. \citealp{AuSt:17} and the references therein\footnote{In particular, see also \cite{PeCoKeHoFe:96,BuGrHa:97,StEiHa:99,Gree:00,ViMC:07,CoCoAgGAPe:11,LYLES2012,Pavlou2015,GrMaAl:16}.}).  

The problem of separation or near separation in logistic regression is linked to the existence of the MLE which is not always guaranteed \citep[see][]{Silv:81,AlAn:84,SaDu:86,Zeng:17}. Alternative conditions for the existence of the MLE have been recently developed in \cite{CaSu:18} \citep[see also][]{SuCa:18}. In order to detect separation, several approaches have been proposed (see for instance, \citealp{LeAl:89,Kolo:97,CHRISTMANN2001}). The bias correction proposed by \cite{Firt:93} has the additional natural property that it is not subject to the problem of separability\footnote{It has been proposed in several places as an alternative to the MLE; see e.g. \cite{HeSc:02,BULL2002,Heinze2006ACI,BuLeLe:07,HePu:10,Wang:14,GrMa:15}.}. However, as noticed in for example \cite{Pavlou2016,Rahman2017,PuHeNoLuGe:17}, in the case of rare events, i.e. when the response classes are unbalanced (for example a positive outcome $\bm{Y}_i=1$ occurs say only $5\%$ of the times), the corrected MLE can suffer from a bias towards one-half in the predicted probabilities. This is because it can been seen as a penalized MLE with Jeffreys invariant prior \citep{Jeff:46}, which is also the case with other penalized MLE proposed for the logistic regression model \citep[see e.g.][]{LCHo:92,gelman2008,Cadigan2012}.  

Moreover, the MLE is known to be sensitive to slight model deviations that take the form of outliers in the data, so that several robust estimators for the logistic regression model and more generally for GLM, have been proposed (see e.g. \cite{CaRo:01b,Cize:08,HeCaCoVF:09}  and the references therein\footnote{Actually, many robust estimators have been proposed for the logistic regression in particular, and for the GLM in general; see e.g. \cite{Preg:81,Preg:82,StCaRu:86,Copa:88,KuStCa:89,Morg:92,CaPe:93,Chri:94,BiYo:96,MaBaLi:97,KORDZAKHIA2001,ADIMARI2001,MPVF-RobBin-02,Rousseeuw2003,CROUX2003,GERVINI2005,HOBZA2008,BeYo:11,HOSSEINIAN2011,Tabatabai2014robust} \citep[and for categorical covariates, one can use the MGP estimator of ][]{MPVF-Grouped:97}. Another type of robust estimators is based on weights that are proportional to the model density through density power divergence measures  \citep[see][]{Wind:95,basu1998robust,Choi2000robust,Jones2001robust}, and for the (polytomous) logistic regression model, see e.g. \cite{Blondell2005robust,Ghosh2016,Castilla2018robust}. The quantity of proposals might indicate that considering robustness properties only, for these models, is not sufficient to provide suitable estimators in practical settings.}). 

Despite all the proposals for finite sample bias correction, separation and data contamination problems, no estimator is so far able to handle the three potential sources of biases altogether. In this section, we propose a simple adaptation of available estimators that we put in the framework of the JIE. Although our choices are certainly not the best ones, they at least guarantee, at a reasonable computational cost, estimators that have a reduced finite sample bias comparable, for example, to the implementation of \cite{Firt:93} by \cite{KoFi:09} \citep[see also][]{kosmidis2010} in favourable settings. Moreover, their performance in terms of mean squared error appears to improve in contaminated and unbalanced settings. 

\subsection{Bias corrected estimators}

Consider the logistic regression model with (observed) response $\bm{Y}(\bm{\beta}_0, n, \bm{\omega}_0)$ (with $\bm{Y}_i(\bm{\beta}_0, n, \bm{\omega}_0)$, for $i=~1,\ldots,n$, as its elements) and linear predictor $\mathbf{X}\boldsymbol{\beta}$, $\mathbf{X}$ being an $n\times p$ matrix of fixed covariates with row $\mathbf{x}_i,i=1,\ldots,n$, and with logit link $\mathbb{E}[\bm{Y}_i(\bm{\beta}, n, \bm{\omega})] \equiv \mu_i(\bm{\beta})=\exp(\mathbf{x}_i\boldsymbol{\beta})/(1+\exp(\mathbf{x}_i\boldsymbol{\beta}))$, for all $\bm{\omega} \in \bOmega$. The MLE for $\bbeta$ is given by
\begin{eqnarray}
  \hat{\bm{\pi}}(\bm{\beta}_0, n, \bm{\omega}_0)& \equiv &\argzero_{\bbeta}  \frac{1}{n}\sum_{i=1}^n\mathbf{x}_i\left[\bm{Y}_i(\bm{\beta}_0, n, \bm{\omega}_0)-\mu_i(\bm{\beta})\right]\nonumber \\ & = &
  \argzero_{\bbeta}  \frac{1}{n}\sum_{i=1}^n\mathbf{x}_i\left[\bm{Y}_i(\bm{\beta}_0, n, \bm{\omega}_0)-\frac{\exp(\mathbf{x}_i\boldsymbol{\beta})}{1+\exp(\mathbf{x}_i\boldsymbol{\beta})}\right],
  \label{Eq_MLE-logistic}
\end{eqnarray}
and can be used as auxiliary estimator in (\ref{eq:indirectInf:hTimesN}) together with the IB in (\ref{eq:iterboot}). The MLE in (\ref{Eq_MLE-logistic}) can also be written as an $M$-estimator \citep{Hube:81} with corresponding $\psi$-function given by
\begin{equation}
 \psi_{\hbpi}\left(\bm{Y}_i(\bm{\beta}_0, n, \bm{\omega}_0);\mathbf{x}_i,\bbeta\right)= \mathbf{x}_i\left[\bm{Y}_i(\bm{\beta}_0, n, \bm{\omega}_0)-\mu_i(\bbeta)\right] .
 \label{Eq_psy-MLE-logistic}
\end{equation}
To avoid the (potential) problem of separation, we follow the suggestion of \cite{Rousseeuw2003} to transform the observed responses $\bm{Y}_i(\bm{\beta}_0, n, \bm{\omega}_0)$ to get pseudo-values (ps)
\begin{equation}
    \tilde{\bm{Y}}_i(\bm{\beta}_0, n, \bm{\omega}_0) = (1-\delta)\bm{Y}_i(\bm{\beta}_0, n, \bm{\omega}_0) +\delta,
    \label{Eq_pseudo-val}
\end{equation}
for all $i=1,\dots,n$, where $\delta \in \left[0,0.5\right)$ is a (fixed) ``small'' (i.e. close to 0). For a discussion on the choice of $\delta$ and also possible asymmetric transformations see for example \cite{Rousseeuw2003}. Note that this transformation is deterministic, hence not subject to sampling error. We call the resulting estimator the JIE-MLE-ps.

As robust estimator, we consider the robust $M$-estimator proposed by \cite{CaRo:01b}, with general $\psi$-function, for the GLM, given by
\begin{equation}
    \psi_{\hbpi}\left(\bm{Y}_i(\bm{\beta}_0, n, \bm{\omega}_0) | \mathbf{x}_i,\bbeta\right)=\psi_c\left(r_i\left(\bm{Y}_i(\bm{\beta}_0, n,\bm{\omega}_0) | \mathbf{x}_i,\bbeta\right)\right)w\left(\mathbf{x}_i\right)V^{-1/2}\left(\mu_i(\bbeta)\right)(\partial/\partial\bbeta)\mu_i(\bbeta)-a\left(\bbeta\right),
    \label{Eq_rob-glm}
\end{equation}
with $r_i\left(\bm{Y}_i(\bm{\beta}_0, n,\bm{\omega}_0);\mathbf{x}_i,\bbeta\right)=\left(\bm{Y}_i(\bm{\beta}_0, n,\bm{\omega}_0)-\mu_i(\bbeta)\right)/V^{-1/2}\left(\mu_i(\bbeta)\right)$, the Pearson residuals and with consistency correction factor
\begin{equation}
    a\left(\bbeta\right)=\frac{1}{n}\sum_{i=1}^n\mathbb{E}\left[    \psi_c\left(r_i\left(Y_i | \mathbf{x}_i,\bbeta\right)\right)w\left(\mathbf{x}_i\right)V^{-1/2}\left(\mu_i(\bbeta)\right)(\partial/\partial\bbeta)\mu_i(\bbeta)\right],
    \label{Eq_rob-consist}
\end{equation}
where the expectation is taken over the (conditional) distribution of the responses $\bm{Y}_i$ (given $\mathbf{x}_i$). For the logistic regression model, we have  $V\left(\mu_i(\bbeta)\right)=\mu_i(\bbeta)(1-\mu_i(\bbeta))$. We 
compute the robust auxiliary estimator $\hbpi(\bm{\beta}_0, n, \bm{\omega}_0)$ on the pseudo-values (\ref{Eq_pseudo-val}), using the implementation in the R \emph{glmrob} function of the \emph{robustbase} package \citep{robustbase2018}, with for $\psi_c$ in (\ref{Eq_rob-glm}) being the Huber's loss function (with default parameter $c$) (see \citealp{huber1964robust}) and $w\left(\mathbf{x}_i\right)=\sqrt{1-h_{ii}}$, $h_{ii}$ being the diagonal element of $\mathbf{X}\left(\mathbf{X}^T\mathbf{X}\right)^{-1}\mathbf{X}^T$. The resulting robust estimator is taken as the auxiliary estimator in (\ref{eq:indirectInf:hTimesN}) and, using the  IB in (\ref{eq:iterboot}), one obtains a robust JIE that we call JIE-ROB-ps. Both auxiliary estimators are not consistent estimators, but both JIE enjoy the properties of being consistent and have a reduced finite sample bias, with JIE-ROB-ps being, additionally,  also robust to data contamination. 

\subsection{Simulation study}
\label{Sec_sim-logistic}

We perform a simulation study to validate the properties of the JIE-MLE-ps and JIE-ROB-ps and compare their finite sample performance to other well established estimators. In particular, as a benchmark, we also compute the MLE, the bias reduced MLE (MLE-BR) using the R \emph{brglm} function (with default parameters) of the \emph{brglm} package \citep{brglm2017}, as well as the robust estimator (\ref{Eq_rob-glm}) using the R \emph{glmrob} function without data transformation (ROB). We consider four situations that can occur with real data, that is, balanced outcomes classes (Setting I) and unbalanced outcome classes (Setting II) with and without data contamination. Setting II is created by adapting the value of the intercept $\beta_0$. We also consider a moderately large model with $p=20$ and chose $n$ as to provide EPV of respectively $5$ and $1.5$. The parameters setting for the simulations are provided in Table \ref{tab:sim-logistic} below. 

\begin{figure}[!ht]
    \centering
    \includegraphics[width=15cm]{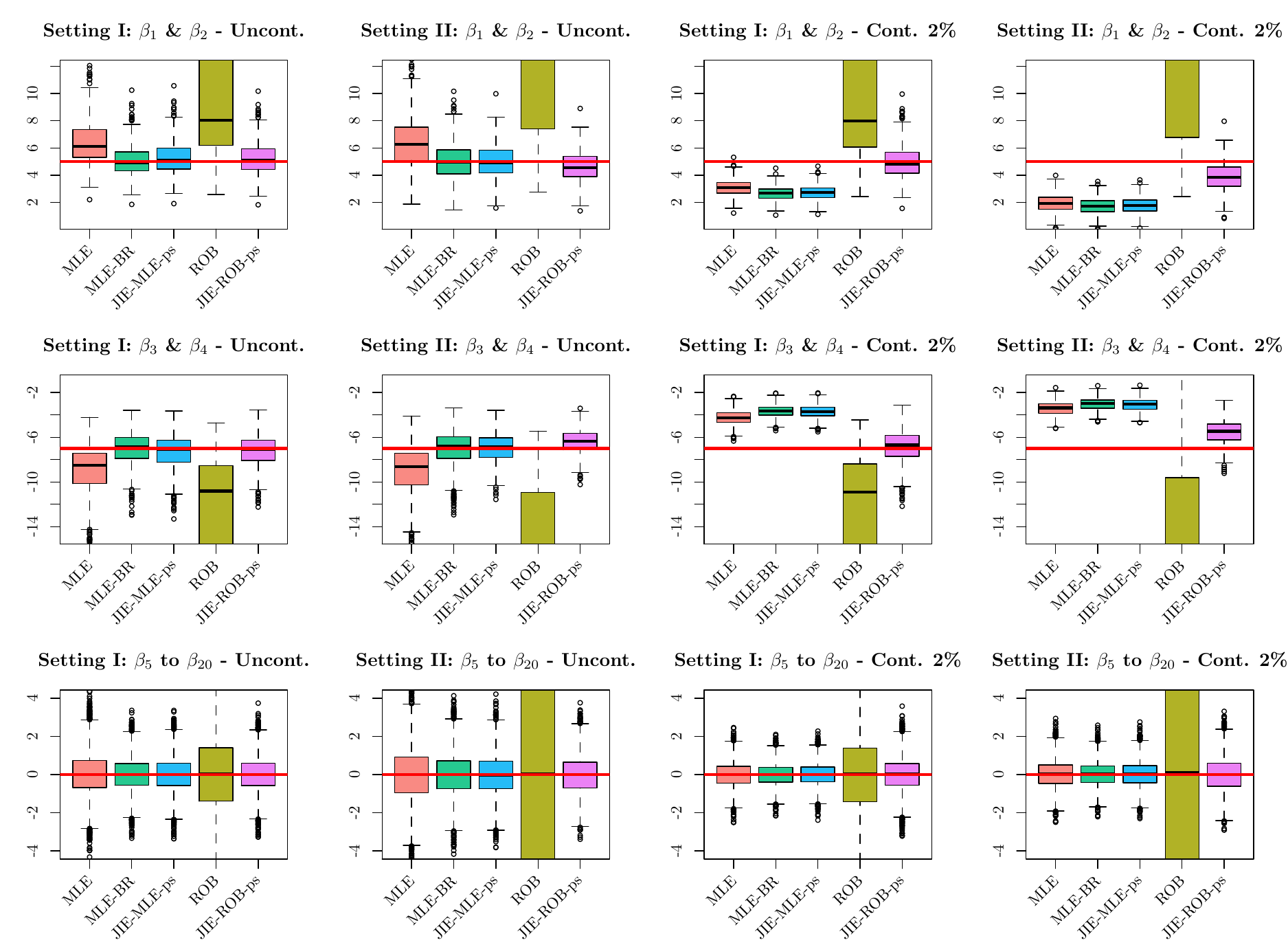}
    \caption{Finite sample distribution of estimators for the logistic regression model using the simulation settings presented in Table \ref{tab:sim-logistic}. The estimators are the MLE (MLE), the Firth's bias reduced MLE (MLE-BR), the JIE (JIE-MLE-ps) with the MLE computed on the pseudo values (\ref{Eq_pseudo-val}) as auxiliary, the robust estimator in (\ref{Eq_rob-glm}) (ROB) and the JIE (JIE-ROB-ps) with the robust estimator  computed on the pseudo values (\ref{Eq_pseudo-val}) as auxiliary.}
    \label{fig:sim-logistic-bxp}
\end{figure}

\begin{figure}[!ht]
    \centering
    \includegraphics[width=15cm]{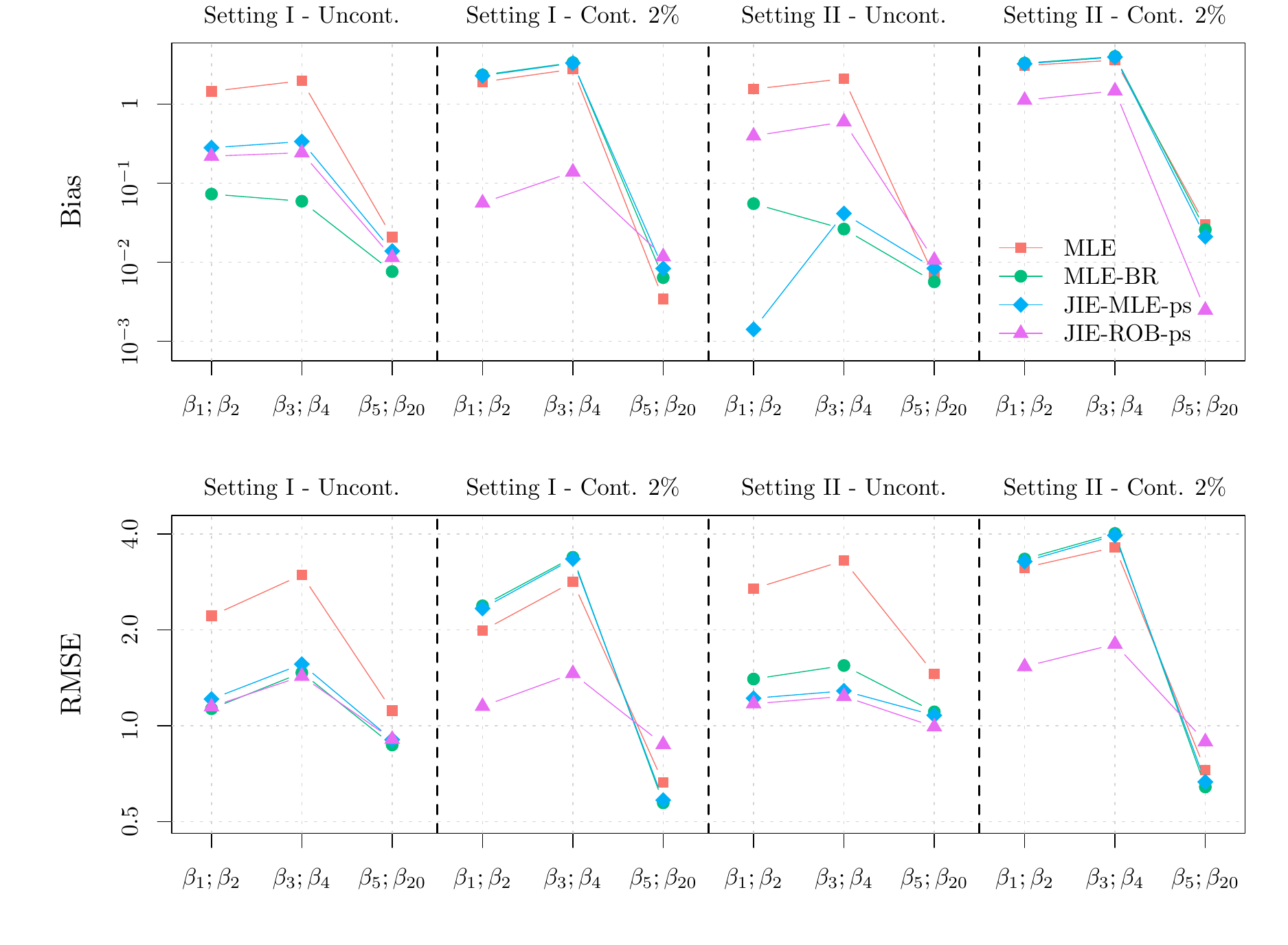}
    \caption{Finite sample bias and RMSE of estimators for the logistic regression model using the simulation settings presented in Table \ref{tab:sim-logistic}. The estimators are the MLE (MLE), the Firth's bias reduced MLE (MLE-BR), the JIE (JIE-MLE-ps) with the MLE computed on the pseudo values (\ref{Eq_pseudo-val}) as auxiliary and the JIE (JIE-ROB-ps) with the robust estimator in (\ref{Eq_rob-glm}) computed on the pseudo values (\ref{Eq_pseudo-val}) as auxiliary. The bias and RMSE of the robust estimator in (\ref{Eq_rob-glm}) (ROB) was omitted in order to avoid an unsuitable scaling of the graphs.}
    \label{fig:sim-logistic-summary}
\end{figure}

\begin{table}[!hb]
    \centering
    \begin{tabular}{lrr}
\toprule
Parameters &  Setting I & Setting II \\
\midrule
$p=$ & $20$ & $20$ \\
$n=$ & $200$ & $300$ \\
$\sum_{i=1}^ny_i\approx$ & $100$ & $270$ \\
EPV $\approx$ & $5$ & $1.5$ \\
$H=$ & $500$ & $400$ \\
$\beta_0=$ & $0$ & $5$ \\
$\beta_1=\beta_2=$ & $5$ & $5$ \\
$\beta_3=\beta_4=$ & $-7$ & $-7$ \\
$\beta_5=\ldots=\beta_{20}=$ & $0$ & $0$ \\
$\delta=$ & $0.05$ & $0.05$ \\
\bottomrule
    \end{tabular}
    \caption{Simulation settings for the logistic regression model. The number of simulations is $1000$ in both settings.}
    \label{tab:sim-logistic}
\end{table}
The covariates were simulated independently from a $\mathcal{N}(0,4/\sqrt{n})$, in order to ensure that the size of the log-odds ratio $\mathbf{x}_i\bbeta$ does not increase with $n$, so that $\mu_i(\bbeta)$ is not trivially equal to either $0$ or $1$ (see e.g. \citealp{SuCa:18}). To contaminate the data, we chose a rather extreme miss classification error to show a noticeable effect on the different estimators, which consists in permuting 2\% of the responses with corresponding larger (smaller) probabilities (expectations). The simulation results are presented in Figure \ref{fig:sim-logistic-bxp} as boxplots of the finite sample distribution, and in Figure \ref{fig:sim-logistic-summary} as the bias and Root Mean Squared Error (RMSE) of the different estimators.

The finite sample distributions presented in Figure \ref{fig:sim-logistic-bxp}, as well as the summary statistics given by the bias and RMSE presented in Figure \ref{fig:sim-logistic-summary}, allow us to draw the following conclusions that support the theoretical results. In Setting I, where the outcome classes are balanced, as expected the MLE is biased (except when the slope parameters are zero) and the MLE-BR, JIE-MLE-ps and JIE-ROB-ps are all unbiased which is not the case for the robust estimator ROB. Moreover, the variability of all estimators is comparable, except for ROB which makes it rather inefficient in this setting. With 2\% of contaminated data (missclassification error), the only unbiased estimator is JIE-ROB-ps and its behaviour remains stable compared to the uncontaminated data setting. This is in line with a desirable property of robust estimators, that is stability with or without (slight) data contamination. In Setting II where the outcome classes are unbalanced, with an EPV of $1.5$, without data contamination, all estimators have a performance loss in terms of RMSE except for the JIE-MLE-ps which performs (slightly) better than the MLE-BR. With contaminated data, although the best performance in terms of bias and RMSE is achieved  the JIE-ROB-ps, all other estimators are biased. However, as argued above, a better proposal for a robust, bias reduced and consistent estimator, as an alternative to JIE-ROB-ps, could in principle be proposed, but this is left for further research.

\section{Bias-Corrected Over-Dispersion Estimator for the Negative Binomial Regression Model}
\label{sec:negbinomial}

When analyzing count responses in a causal relationship, the Poisson distribution within the GLM framework is a natural option. However, data often exhibit overdispersion defined as the extra variation which is not explained by the Poisson distribution alone \citep[see e.g.][]{McCuNe:89,Cons:89,Jorg:97}. One approach to modelling count data with overdispersion is by means of a mixing distribution to the Poisson mean parameter, an approach that dates back to \cite{GrYu:20} who used the gamma as Poisson mixing distribution to obtain the Negative Binomial (NB) distribution\footnote{See also \cite{MaWoSt:81,Bres:84,Enge:84,Lawl:87,BaPa:88,Pieg:90,PaBa:98,YoSm:05,Hilb:11,CaTr:13,Hilb:14}.}. Basically, accounting for overdispersion is done by including an extra random effect variable, and conditionally on this random effect, a Poisson distribution is assumed for the response. A marginal model is then obtained by integrating the conditional distribution over the random effect distribution\footnote{Mixing distributions based on the Poisson are very popular models for count data;  for reviews, simulations and applications of different models, see e.g. \cite{Joe2005GeneralizedPD,KaXe:05,RiStAk:08,NiKa:08,ZhJo:09,HeCoHe:18}.}. In general, with models that account for overdispersion, it is well known that the MLE can be seriously biased in finite sample dimensions encountered in practice, and some authors have proposed finite sample adjustments of order $\mathcal{O}\left(n^{-1}\right)$; see for example \cite{SiRoBS:11} and the references therein\footnote{One can in particular mention \cite{BoCo:98,CoBo:01,JoKn:04,SaPa:05,LlSm:07,CoTU:08}.}. This finite sample bias is especially large for the MLE of the overdispersion parameter, which is widely used as an important measure in medical, biological and environmental studies\footnote{See e.g. \cite{ByAlGiPe:03,LSScKoGe:05,PaSa:07,RoSm:08,Kris:11}.}. In this section, we use the general framework developed in this paper in order to provide a JIE with reduced finite sample bias (and RMSE) for the parameters of the NB regression model.

The NB regression model with response $\bm{Y}$ (with elements $\bm{Y}_i\in \{0,1,\ldots,m\}$, $i=1,\ldots,n$) and linear predictor $\mathbf{X}\boldsymbol{\beta}$, $\mathbf{X}$ being an $n\times p$ matrix of fixed covariates with row $\mathbf{x}_i$, $i=1,\ldots,n$, and with exponential link $\mathbb{E}[\bm{Y}_i]=\mu_i\equiv\exp(\mathbf{x}_i\boldsymbol{\bbeta})$, has conditional (on $\mathbf{x}_i$, $i=1,\ldots,n$) density given by
\begin{equation}
    f(\bm{Y}_i = y\vert \mathbf{x}_i;\bbeta,\alpha)=\frac{\Gamma\left(y+\alpha\right)}{y!\Gamma\left(\alpha\right)}\left(\frac{\alpha}{\alpha+\mu_i}\right)^{\alpha}\left(\frac{\mu_i}{\alpha+\mu_i}\right)^{y}.
    \label{Eq_NB-reg-density}
\end{equation}
The variance function is $V(\mu_i)=\mu_i+\mu_i^2/\alpha$ and hence
decreasing values of $\alpha$ correspond to increasing levels of dispersion and as $\alpha\rightarrow \infty$, one gets the Poisson distribution. 

We propose to construct two estimators. One is the JIE in (\ref{eq:indirectInf:hTimesN}) with, as auxiliary (consistent) estimator, the MLE $\hat{\bm{\pi}}\left(\btheta_0, n, \bm{\omega}_0\right)$ ($\btheta_0=\left(\bbeta_0^T,\alpha_0\right)^T$) associated to the density (\ref{Eq_NB-reg-density}). To compute the MLE, we use the R \emph{glm.nb} function of the \emph{MASS} package \citep{VeRi:02}. The JIE is computed using the  IB in (\ref{eq:iterboot}), and we call it the JIE-MLE. 
%
%
The other estimator we propose is the JIE based on an alternative (inconsistent) auxiliary estimator computed in two steps. The slope coefficients $\bbeta$ are first estimated using the MLE $\hat{\bm{\pi}}_{\bbeta}\left(\btheta_0, n, \bm{\omega}_0\right)$ for the Poisson model, namely, using the $\psi$-function in (\ref{Eq_psy-MLE-logistic}) with $\mu_i(\bbeta)=\exp(\mathbf{x}_i\boldsymbol{\bbeta})$
%
%
and the overdispersion parameter $\alpha$ is estimated in a second step using the variance expression for the $\mu_i$'s, i.e.
\begin{equation*}
\hat{\bpi}_{\alpha}\left(\btheta_0, n, \bm{\omega}_0\right)=\frac{\sum_{i=1}^n\mu_i^2\left(\hat{\bm{\pi}}_{\bbeta}\left(\btheta_0, n, \bm{\omega}_0\right)\right)}{\sum_{i=1}^n\left(\bm{Y}_i\left(\btheta_0, n, \bm{\omega}_0\right)-\mu_i\left(\hat{\bm{\pi}}_{\bbeta}\left(\btheta_0, n, \bm{\omega}_0\right)\right)\right)^2-\sum_{i=1}^n\mu_i\left(\hat{\bm{\pi}}_{\bbeta}\left(\btheta_0, n, \bm{\omega}_0\right)\right)}.
\end{equation*}
The resulting JIE is computed using the IB in (\ref{eq:iterboot}) and we call it the JIE-MLE-P. Its main advantage is its computational efficiency (only two steps) which, in high dimensional settings, could result in the only numerically viable estimator  for the NB regression model.

\subsection{Simulation study}

\begin{figure}[ht]
    \centering
    \includegraphics[width=15cm]{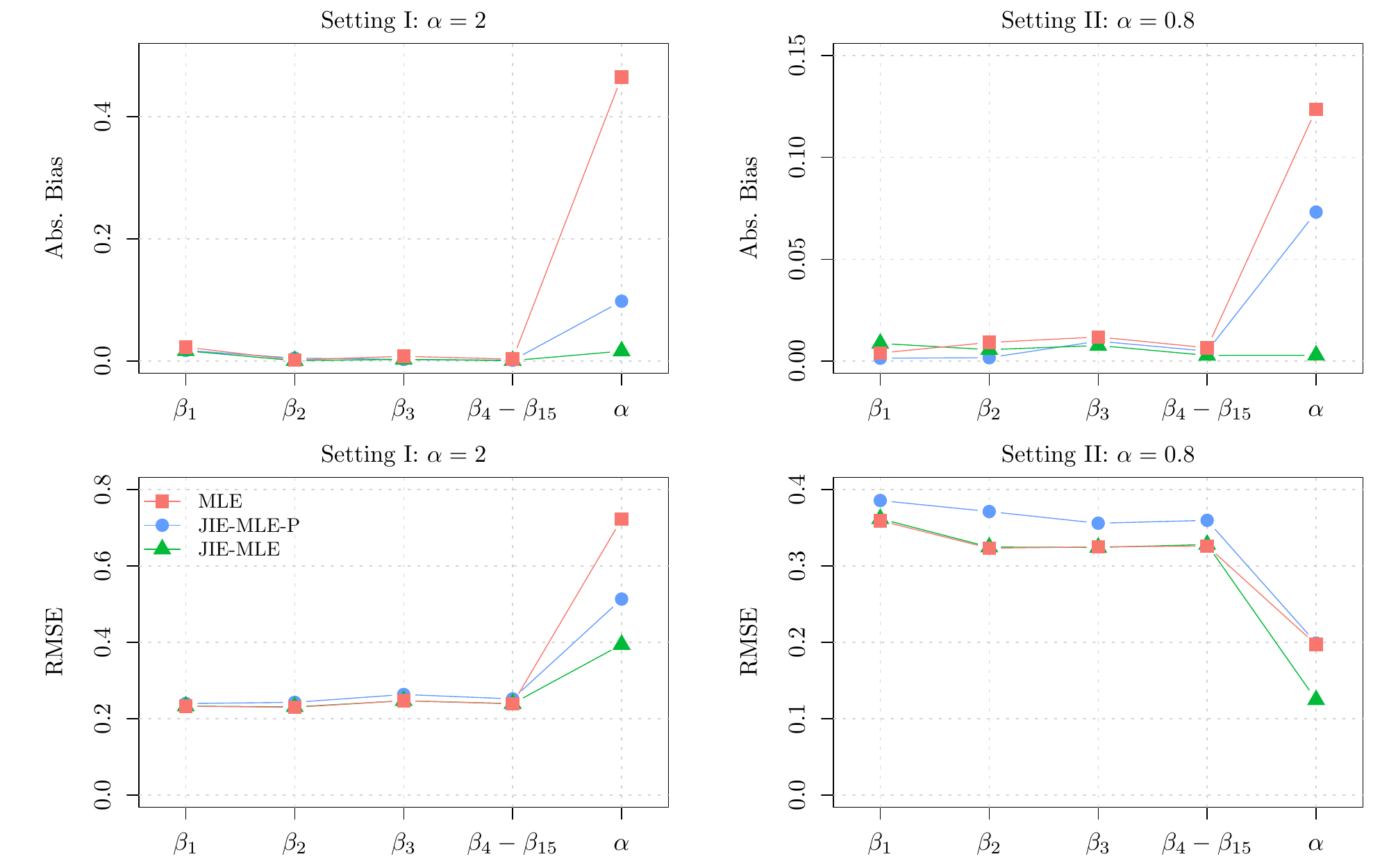}
    \caption{Finite sample bias and RMSE of estimators for the negative binomial regression model using the simulation settings presented in Table \ref{tab:sim-NB}. The estimators are the MLE (MLE) and the JIE with auxiliary estimator the MLE (JIE-MLE) and the naive estimator (JIE-MLE-P).}
    \label{fig:sim-NB-summary}
\end{figure}

We perform a simulation study to validate the properties of the JIE-MLE and the JIE-MLE-P. We compare their finite sample performance to the MLE. We consider two values for the overdispersion parameter, namely highly overdispersed data with $\alpha=0.8$ and less dispersed data with $\alpha=2$. The parameter's setting for the simulations are provided in Table \ref{tab:sim-NB}. The covariates were simulated independently from a uniform distribution between 0 and 1.  The simulation results are presented in Figure \ref{fig:sim-NB-summary} through the bias and RMSE of the three estimators.

\begin{table}[!hb]
    \centering
    \begin{tabular}{lrr}
    \toprule
Parameters &  Setting I & Setting II \\
\midrule
$p =$ & $15$ & $15$    \\
$n =$ & $150$ & $150$ \\
$H =$ & $500$ & $500$  \\
$\beta_0 =$ & $1$ & $1$  \\
$\beta_1 =$ & $2$ & $2$ \\
$\beta_2 =$ & $-1$ & $-1$  \\
$\beta_3=\ldots=\beta_{15} =$ & $0$ & $0$ \\
$\alpha =$ & $2$ & $0.8$ \\
\bottomrule
    \end{tabular}
    \caption{Simulation settings for the negative binomial regression model. The number of simulations is $1000$ in both settings.}
    \label{tab:sim-NB}
\end{table}

The summary statistics given by the bias and RMSE presented in Figure \ref{fig:sim-NB-summary} allow us to draw the following conclusions that support our theoretical results. While the MLE for the slope parameters $\bm{\beta}$ is not noticeably biased (or just slightly in the more overdispersed case), the MLE for the overdispersion parameter $\alpha$ suffers from an important finite sample bias. Both JIE, on the other hand, have quite smaller finite sample biases without increasing the variability, as reported by the simulations RMSE. The bias is nil with the JIE-MLE and  strongly reduced with the \mbox{JIE-MLE-P}, and for both, the RMSE reduction is noticeably larger compred to the MLE. This conclusion supports the theoretical results, since the MLE is a consistent estimator, and the bias therefore disappears when it is used as auxiliary estimator for the JIE. The relatively large bias of the MLE has as a consequence that, when performing inference using the associated asymptotic results, for example when building confidence intervals, their coverage is below the nominal level, given that the overdispersion parameter $\alpha$ is generally overestimated, and this happens in both more and less overdispersed data. Using the JIE-MLE, seems to greatly reduce the risk of falling in these situations by largely reducing the bias of the overdispersion parameter.

\section{Bias-corrected Regularized Linear Regression Estimators}
\label{Sec_lasso}

In his seminal article, \cite{Tibs:96} provides two main justifications for the use of regularized regression, i.e. least square estimation under constraints on the regression parameters. The first is prediction accuracy, since when $p$ is large, the Ordinary Least Squares (OLS) estimator, although unbiased (in the linear case), can have a large variance, and prediction accuracy can be improved by means of parameter's shrinkage. Shrinkage induces a bias, but might decrease the Mean Squared Error (MSE). The second reason is interpretation, since with a large number of predictors, one often prefers to identify a smaller subset that exhibits the strongest signal. Here again, constraining the OLS, introduces an estimation bias.

Among the traditional (convex) Regularized linear Regression Estimators (RRE), the ridge regression estimator \citep{HoKe:70}, based on Tikhonov regularization method \citep{Tikh:43}, and the lasso \citep{Tibs:96} are the most well known. For the linear regression model with response $\bm{Y}$ (with elements $\bm{Y}_i, i=1,\ldots,n$) and an $n\times p$ design matrix $\mathbf{X}$ with row $\mathbf{x}_i$, $i=1,\ldots,n$, a RRE $\hat{\bpi}_{\lambda}\left(\bbeta_0, n, \bm{\omega}_0\right)$ for the slope coefficients $\bbeta$\footnote{For simplicity of exposition, we consider the case when the residual variance $\sigma^2$ is either known or its estimation is not an important issue, and leave a deeper consideration to further research.}, for two common ones,  can be defined as
\begin{equation}
  \hat{\bpi}_{\lambda}\left(\bbeta_0, n, \bm{\omega}_0\right)=\argmin_{\bbeta} \;\frac{1}{n} \big\lVert\bm{Y}\left(\bbeta_0, n, \bm{\omega}_0\right)- \mathbf{X}\bbeta\big\rVert^2_2+\lambda\lVert\bbeta\rVert_{\phi}^{\phi}\,,
  \label{Eq_RRE}
\end{equation}
with  $\lambda$ chosen adequately. If $\phi=2$ then $\hat{\bpi}_{\lambda}\left(\bbeta_0, n, \bm{\omega}_0\right)$ is the ridge estimator and if $\phi=1$, it is the lasso (up to a reparametrization of $\lambda$). One of the very advantageous features of the lasso is that it automatically sets some of the estimated slope parameters to 0, so that it is also used as a model selection method \citep[the properties of the lasso estimator have been recently studied in particular by][among many others]{greenshtein2004,LeLiWa:06,zhang2008,bickel2009}. Actually, RRE are suitable in high dimension when $p>n$, since by shrinking the solution space with the penalty $\lambda\lVert\bbeta\rVert_{\phi}^{\phi}$, they are not ill posed as is the case, for example, with the OLS. For a detailed overview see for example \cite{BuhlmannBook2011}.

Debiasing RRE as recently attracted more attention. It takes place in the framework of inference after selection (see e.g. \citealp{Dezeure2017} and the references therein\footnote{In particular, see also \cite{wasserman2009,MeMeBu:09,Meinshausen2010,ShSa:13,buhlmann2013,ZaZa:14,lockhart2014,vandegeer2014,javanmard14,Buhlmann2014,Mein:15,BeChKa:15,dezeure2015,Zhang2017}.}). 
\emph{De-biased} (or \emph{de-sparsified}) RRE are asymptotically normal under less stringent conditions than the lasso estimator, so that reliable post model selection inference can be performed.  
In this section, we study, by means of simulations, the finite sample performance of the JIE, computed using the IB, with the lasso as auxiliary estimator (with fixed values for $\lambda$). Given the simplicity of the implementation, our aim is to study, practically, the reduction in bias and MSE of the JIE, and leave the formal aspects for further research. The study takes place in high ($p>n$) as well as in low ($p<n$) dimensional settings where the OLS is also available.

Interestingly, when $p<n$, if the JIE uses the ridge as the auxiliary estimator, then it is the OLS. Indeed, if we let $\hbbo$ and $\hbbr$ denote respectively the OLS and the ridge estimator and we consider the linear regression model 
\begin{equation*}
    \bm{Y} = \bX\bb_0 + \bm{\varepsilon},\quad \bm{\varepsilon}\sim\mathcal{N}(\0,\sigma_0^2 \bI_n).
\end{equation*}
Then, if the matrix $\left(\bX^T\bX\right)^{-1}$ exists the relationship  between the ridge and the OLS is simply
\begin{equation}
    \hbbr = \Big[\bI_p+\lambda\left(\bX^T\bX\right)^{-1}\Big]^{-1}\hbbo.
    \label{eq:ols:ridge}
\end{equation}
Let $\hat{\bm{\theta}}_{j,n,H}$ denote the JIE based on $\hbbr$, then we would expect $\hat{\bm{\theta}}_{j,n,H}$ to simply lead to the OLS (up to the simulation error) by the linearity of (\ref{eq:ols:ridge}). Naturally, in this case the use of simulated samples is not required as the expected of $\hbbr$ is known. However, if we consider relatively large value of $H$ we could simply write:
\begin{equation*}
    \hat{\bm{\theta}}_{(j,n,H)} \approx \argzero_{\bm{\beta}} \; \Big[\bI_p+\lambda\left(\bX^T\bX\right)^{-1}\Big]^{-1} \left( \hbbo - \bm{\beta}\right) = \hbbo.
\end{equation*}
This simple form of equivalence between the OLS and the ridge doesn't apply to the lasso, i.e. when $\phi=1$ in (\ref{Eq_RRE}). Hence, using an RRE as auxiliary estimator for the JIE provides a new estimator that has not been, to our knowledge, studied so far. The simulation studies presented below show that the JIE is an estimator with most of the slope coefficients remaining at zero, but having a finite sample distribution with reduced MSE compared to the lasso and to a \emph{debiased-lasso} estimator.

\subsection{Simulation study}

In order to study the performance of the JIE with the lasso as auxiliary estimators (JIE-lasso), we compare it to available benchmarks. Indeed, we consider the lasso itself but also the  \emph{debiased-lasso} (lasso-D) of~\citet{javanmard14} as an $R$ implementation is readily available\footnote{We use the recommended value of $\mu = 2\sqrt{\log(p)/n}$.}. When $p<n$, we also compute the OLS.

\begin{figure}[ht]
    \centering
    \includegraphics[width=10cm]{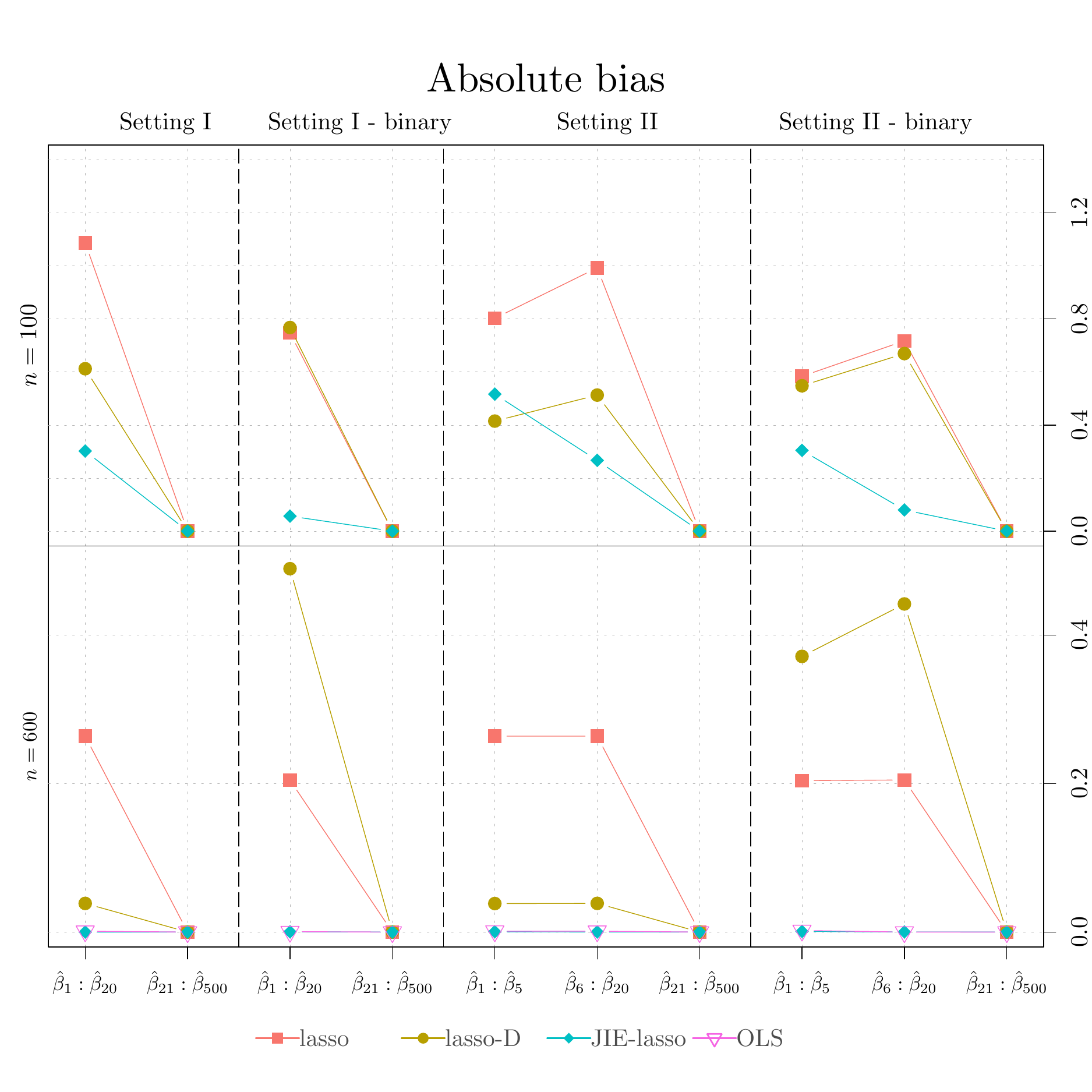}
    \caption{Absolute bias of the lasso, the \emph{debiased lasso} ( (lasso-D) of~\cite{javanmard14}, the OLS and the JIE with the lasso as auxiliary estimator (JIE-lasso). The penalty $\lambda$ in (\ref{Eq_RRE}) is set to $\lambda = 2.5\sqrt{\log(p)/n}$ for the lasso. The IB is used to compute the JIE-lasso with $H=250$. For each simulation setting, $1,000$ samples are generated.}
    \label{fig:sim-lasso-abias}
\end{figure}

\begin{figure}[ht]
    \centering
    \includegraphics[width=10cm]{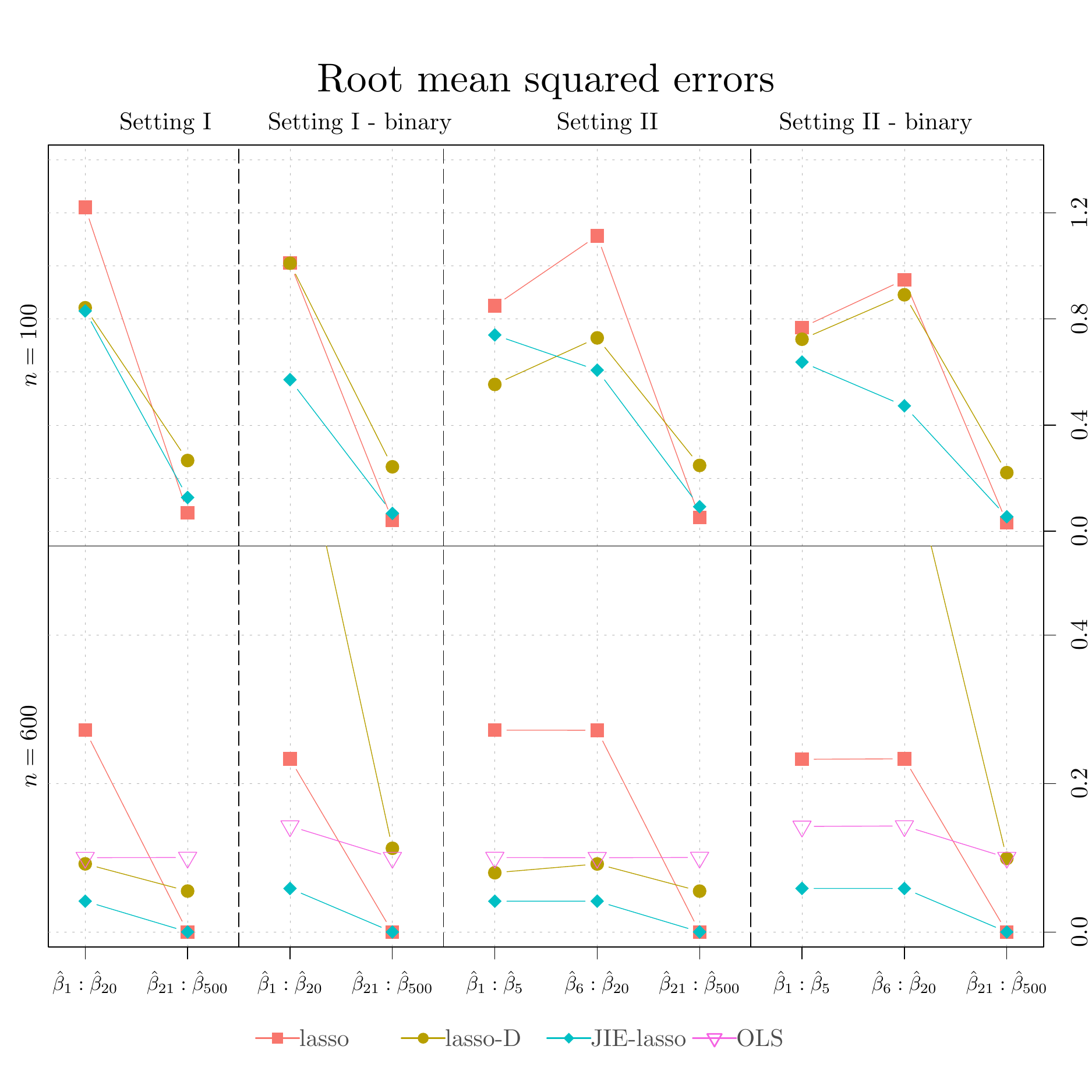}
    \caption{Root mean squared errors of the lasso, the \emph{debiased lasso} ( (lasso-D) of~\cite{javanmard14}, the OLS and the JIE with the lasso as auxiliary estimator (JIE-lasso). The penalty $\lambda$ in (\ref{Eq_RRE}) is set to $\lambda = 2.5\sqrt{\log(p)/n}$ for the lasso. The IB is used to compute the JIE-lasso with $H=250$. For each simulation setting, $1,000$ samples are generated.}
    \label{fig:sim-lasso-rmse}
\end{figure}

\begin{figure}[ht]
    \centering
    \includegraphics[width=10cm]{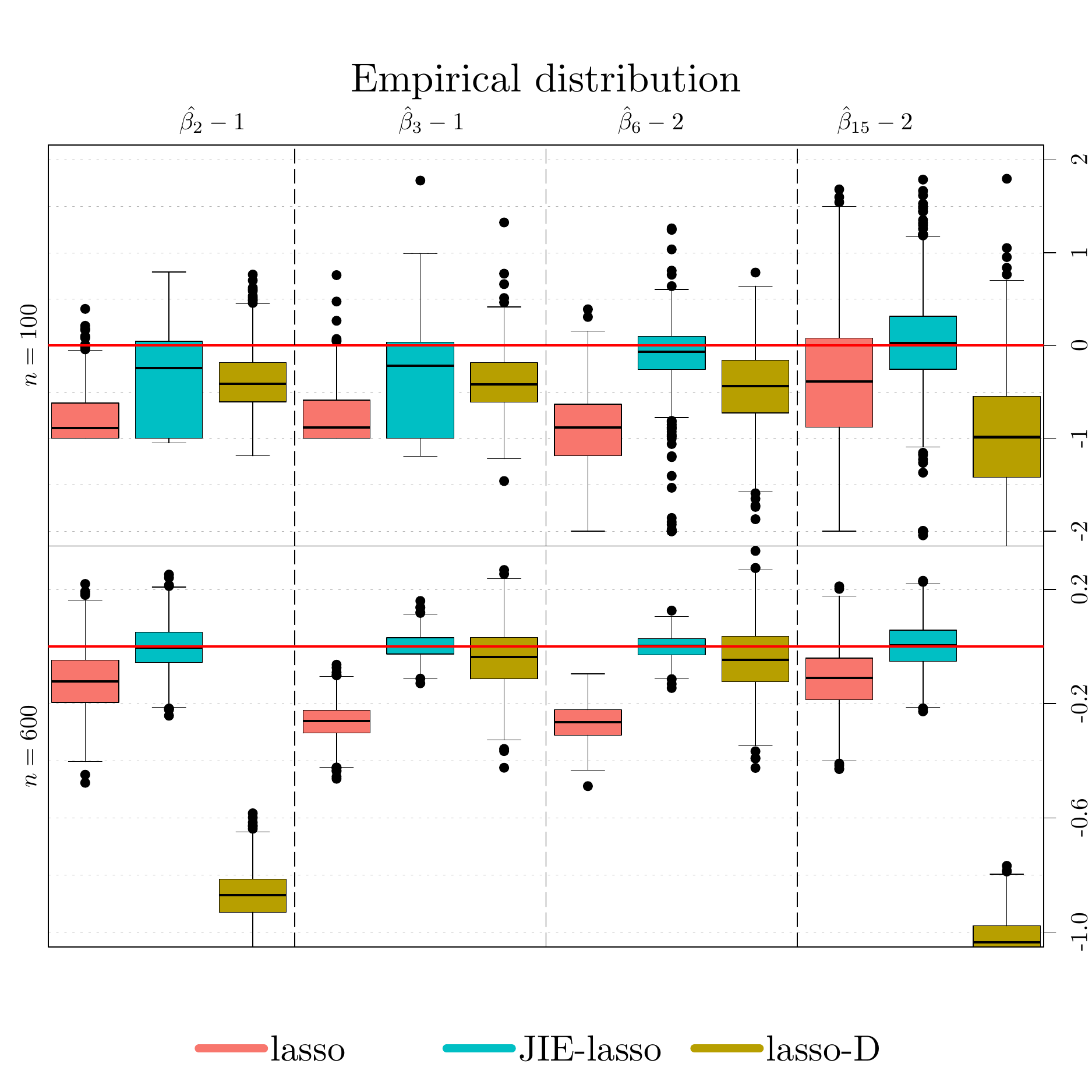}
    \caption{Empirical distribution of regression estimators, in Setting II, using the lasso, the \emph{debiased lasso} (lasso-D) of~\cite{javanmard14} and the JIE-lasso. The penalty $\lambda$ in (\ref{Eq_RRE}) is set to $\lambda = 2.5\sqrt{\log(p)/n}$ for the lasso. The IB with $H=250$ is used to compute the JIE and $1,000$ samples are generated. Both $\hat{\beta}_6$ and $\hat{\beta}_{15}$ are estimating a coefficient with ``strong signal'' ($\beta_6=\beta_{15}=2$), with $X_{15}$ being binary. Both $\hat{\beta}_2$ and $\hat{\beta}_3$ are estimating a coefficient with ``low signal'' ($\beta_2=\beta_3=1$), with $X_{2}$ being binary.}
    \label{fig:sim-lasso-dist}
\end{figure}

For the simulation setting, we use and extend the simulation framework of \cite{DebiasDebiasLasso2017}, who studies the conditions under which the bootstrap can lower the bias of the \emph{debiased-lasso} estimator in high dimensional settings for the linear regression. Namely, we consider two sample sizes, $n = 100$ and $n=600$ and a model's dimension of $p = 500$. For the design matrix, we simulate, as in \cite{DebiasDebiasLasso2017}, uncorrelated standard normal variables, but also the case of discrete (binary) covariates, by setting $X_j\sim \mbox{Bernoulli}(0.5), j=1,2,15,\ldots,20$. For the slope parameters, we consider a sparse setting, with high signal and a mixture of high and low signals, i.e.
\begin{enumerate}
    \item High signal: $\beta_j= 2, j=1,\ldots,20$ and $\beta_j=0,j=21,\ldots,p$. 
    \item High and low signal: $\beta_j = 1,j=1,\ldots,5$, $\beta_j=2, j=6,\ldots,20$ and $\beta_j=0,j=21,\ldots,p$. 
\end{enumerate}
For the penalty $\lambda$ of the lasso, a value $\lambda\asymp\sqrt{\log(p)/n}$ is generally recommended \citep[see e.g.][]{BuhlmannBook2011}, here we set it to $\lambda = 2.5\sqrt{\log(p)/n}$. To compute the JIE-lasso, we use the IB with $H=250$. The performances are measured in terms of absolute bias, RMSE and finite sample distribution, obtained from $1,000$ simulations, and are presented in Figures~\ref{fig:sim-lasso-abias}, \ref{fig:sim-lasso-rmse} and \ref{fig:sim-lasso-dist}.

In terms of finite sample bias (Figure \ref{fig:sim-lasso-abias}), across all simulation settings when $p>n$, the JIE-lasso has the smallest bias. It is comparable with the one of the  \emph{debiased lasso} in Setting 2, with uncorrelated standard normal covariates (no binary ones) and with low signal. The \emph{debiased lasso} has a reduced bias compared to the lasso, except when some covariates are binary. When $p<n$, the OLS and the JIE-lasso have a nill bias in all settings, which was expected for the OLS, but not necessarily for the JIE-lasso.

What is also striking is that the drastic bias reduction provided by the JIE-lasso, does not come with a too strong proportion of slope estimators that are not nill. In Table \ref{tab:lasso_zeros} are reported the proportion of estimators that are set to zero for both the lasso and the JIE-lasso. Although when $n=100$, the JIE-lasso sets more estimates to a non nil value, the proportion that is left to zero is still high. When $n=600$, this proportion is similar to the one of the lasso and near the correct one, but the JIE has a nil bias and reduced RMSE. The \emph{debiased lasso} does not provide estimates that are nil, but it has a smaller finite sample bias when compared to the lasso only. This is the case when the covariates are uncorrelated standard normal, but the bias becomes excessively large (compared to the lasso) when some covariates are binary.

\begin{table}[ht]
    \centering
    \begin{tabular}{lcccccccc}
    \toprule
    & \multicolumn{4}{c}{$n=100$} & \multicolumn{4}{c}{$n=600$} \\
    \cmidrule(r){2-5} \cmidrule(r){6-9}
Setting & I & I-binary & II & II-binary & I & I-binary & II & II-binary \\
    \midrule 
    lasso & 91.43\% & 93.18\% & 92.61\% & 93.92\% & 95.99\% & 95.99\% & 95.99\% & 95.99\% \\
    JIE-lasso & 76.78\% & 79.75\% & 79.25\% & 81.94\% & 95.50\% & 95.81\% & 95.50\% & 95.81\%  \\
    \bottomrule
    \end{tabular}
    \caption{Percentage of times estimators are exactly zero on average across the $1,000$ simulations for the two settings and the two sample sizes. The true proportion of zero in both settings is $96\%$.}
    \label{tab:lasso_zeros}
\end{table}

In terms of RMSE (Figure \ref{fig:sim-lasso-rmse}), the conclusions are similar. When $p>n$, the RMSE of the JIE-Lasso and the \emph{debiased lasso} are comparable when the covariates are standard normal, but the RMSE of the \emph{debiased lasso} is of the same order as the one of the lasso when some covariates are binary. In the $p<n$ case, the RMSE is clearly the smallest for the JIE-lasso compared to the OLS. With the finite sample distributions of estimators of some of the slope coefficients provided in Figure \ref{fig:sim-lasso-dist}, one can see again the JIE-lasso, an estimator very simple to implement, provides an advantageous alternative to the lasso and the \emph{debiased lasso}, especially outside the standard settings, like with binary covariates. When the signal is low and the sample size small, even if the RMSE of the JIE-lasso is smaller, the behaviour of the estimator is less satisfactory.

\section{Conclusions} 
\label{sec:conclude}

The JIE, an estimator based on an other auxiliary estimator,  that we propose in this paper enjoys, under suitable and arguably reasonable conditions, the property of being unbiased in finite (sufficiently large) samples, consistent and having an asymptotic normal distribution for inference. These properties are derived in high dimensional settings, and although all conditions can be put under the same umbrella, they have been separately set for each of the desirable properties and settings. Moreover, we show that the IB is a numerically efficient and convergent algorithm that can be used to compute the JIE. In practice, the crucial issue is then  how to define a \emph{suitable} auxiliary estimator and several strategies can be considered. If finite sample unbiased estimators are sought, then the safest option is to chose a consistent auxiliary estimator (e.g. the MLE). If the issue is numerical in the sense that it is not viable to compute a consistent estimator (especially in high dimensions), then the auxiliary estimator should be chosen to be simple to compute and possibly ``not too far'' from the consistent one. For constrained estimators that are naturally biased such as the lasso, they can be taken as auxiliary estimators for the JIE in order to obtain de-baised and constrained estimators.

In (\ref{eq:indirectInf:hTimesN}), we are implicitly assuming that to compute the JIE, synthetic samples can be generated under the model $F_{\bm{\theta}}$, which fully defines the data generating process. It is possible to extend this setting to some semiparametric cases, but we leave this problem for further research. 
Finally, since the efficiency of the JIE depends on the choice of the auxiliary estimator, the asymptotic variance of the JIE provided in Theorem \ref{thm:main:3}, could be used to chose among competing auxiliary estimators.

\newpage

\setcounter{equation}{0}
\renewcommand{\theequation}{\Alph{section}.\arabic{equation}}
\appendix	
\begin{center}
		\huge \textsc{Appendices}
\end{center}

\section{Proof of Lemma \ref{lemma:unique:iter:boot}}
\label{proof:lemma:1}
	
We define the function 
		\begin{equation}
			\begin{aligned}
				T_{(j,n,H)} \;\;:\;\; & \bm{\Theta}  & \longrightarrow  &\;\;\;\; \real^p\\
				& \bm{\theta}  & \longmapsto & \;\;\;\; T_{(j,n,H)}(\bm{\theta}),\\
			\end{aligned}
			\label{eq:define:T}
		\end{equation}
		where 
		\begin{equation*}
			T_{(j,n,H)}(\bm{\theta}) \equiv \bm{\theta} + \hat{\bm{\pi}}(\bm{\theta}_0, n, \bm{\omega}_0) - \frac{1}{H} \sum_{h = 1}^H  \hat{\bm{\pi}}(\bm{\theta}, n, \bm{\omega}_{h+jH}).
		\end{equation*}
		We define the set of fixed points of $T_{(j,n,H)}$ as follows
		\begin{equation*}
			\bm{\Theta}^\ast_{(j,n,H)} = \left\{\bm\theta \in \bm{\Theta} \;\; \Big| \;\;T_{(j,n,H)}(\bm{\theta}) = \bm{\theta} \right\}.
		\end{equation*}
		Clearly, we have that 
		\begin{equation}
			\bm{\Theta}^\ast_{(j,n,H)} = \widehat{\bm{\Theta}}_{(j,n,H)},
			\label{eq:setfixed_point}
		\end{equation}
		as defined in (\ref{eq:indirectInf:hTimesN}). In the rest of the proof we will start by showing that the function $T_{(j,n,H)}$ admit a unique fixed point in $\real^p$ by applying Banach fixed-point theorem. Then, using (\ref{eq:setfixed_point}) and Assumption \ref{assum:A:1} we will be able to conclude that the set $\widehat{\bm{\Theta}}_{(j,n,H)}$ only contains this fixed point. 
		Let us start by using (\ref{bias:estim}) and, Assumptions \ref{assum:B:1}, \ref{assum:C:1} and \ref{assum:D:1}, we have that
		\begin{equation*}
			\begin{aligned}
	\hat{\bm{\pi}}(\bm{\theta}_0, n, \bm{\omega}_0) &= \bm{\pi} \left( \bm{\theta}_0, n\right)+ \mathbf{v} \left(\bm{\theta}_0,n, \bm{\omega}_0\right)
				=  \bm{\theta}_0 + \mathbf{a}(\btheta_0) + \mathbf{c}(n) + \mathbf{L}(n) \bm{\theta}_0 + \bm{\delta}^{(1)},
			\end{aligned}
		\end{equation*}
		where $\bm{\delta}_i^{(1)} = \mathcal{O}\left(n^{-\gamma}\right) +  \mathcal{O}_{\rm p}\left(n^{-\alpha}\right)$ for $i = 1,\,\ldots,\, p$ and
		\begin{equation*}
			\begin{aligned}
				\frac{1}{H} \sum_{h = 1}^H  \hat{\bm{\pi}}(\bm{\theta}, n, \bm{\omega}_{h+jH}) &= \bm{\pi} \left( \bm{\theta}, n\right) + \frac{1}{H} \sum_{h = 1}^H  \mathbf{v}(\bm{\theta}, n, \bm{\omega}_{h+jH})\\
				&= \bm{\theta} + \mathbf{a}(\btheta) + \mathbf{c}(n) + \mathbf{L}(n) \bm{\theta} + \bm{\delta}^{(2)},
			\end{aligned}
		\end{equation*}
		where $\bm{\delta}^{(2)}_i =	\mathcal{O}\left(n^{-\gamma}\right) + \mathcal{O}_{\rm p}\left(H^{-\nicefrac{1}{2}} n^{-\alpha}\right)$ for $i = 1,\,\ldots,\, p$. Therefore, we have
		\begin{equation*}
			T_{(j,n,H)}(\bm{\theta}) = \bm{\theta}_0 + \mathbf{a}(\btheta_0) - \mathbf{a}(\btheta) + \mathbf{L}(n) \left(\bm{\theta}_0 - \bm{\theta}\right) +  \bm{\delta}^{(1)} - \bm{\delta}^{(2)}.
		\end{equation*}
		Next, we consider the following quantity for $\bm{\theta}_1, \bm{\theta}_2 \in \bm{\Theta}$, 
		\begin{equation*}
			\begin{aligned}
				\left\lVert T_{(j,n,H)}(\bm{\theta}_1) - T_{(j,n,H)}(\bm{\theta}_2)\right\rVert_2^2 &= \left\lVert\mathbf{a}(\btheta_2) - \mathbf{a}(\btheta_1) + \mathbf{L}(n) \left(\bm{\theta}_2 - \bm{\theta}_1\right) +  \bm{\delta}^{(3)} \right\rVert_2^2\\
				 &\leq \left\lVert\mathbf{a}(\btheta_2) - \mathbf{a}(\btheta_1) \right\rVert_2^2 + \left\lVert\mathbf{L}(n) \left(\bm{\theta}_2 - \bm{\theta}_1\right) \right\rVert_2^2 + \left\lVert\bm{\delta}^{(3)} \right\rVert_2^2, 
			\end{aligned}
		\end{equation*}
		where $\bm{\delta}_i^{(3)} = \mathcal{O}\left(n^{-\gamma}\right) +  \mathcal{O}_{\rm p}\left(n^{-\alpha}\right)$ for $i = 1,\,\ldots,\, p$. Then, using the fact that 
		\begin{equation*}
		    \left\lVert\mathbf{L}(n) \left(\bm{\theta}_2 - \bm{\theta}_1\right) \right\rVert_2^2 \leq \left\lVert\mathbf{L}(n)^T \mathbf{L}(n)\right\rVert_F \,  \left\lVert\bm{\theta}_2 - \bm{\theta}_1 \right\rVert_2^2
		\end{equation*}
		where $\lVert\cdot\rVert_F$ is the Frobenius norm, we obtain 
		\begin{equation*}
			\begin{aligned}
				\left\lVert T_{(j,n,H)}(\bm{\theta}_1) - T_{(j,n,H)}(\bm{\theta}_2) \right\rVert_2^2 \leq \left\lVert\mathbf{a}(\btheta_2) - \mathbf{a}(\btheta_1) \right\rVert_2^2 + \left\lVert\mathbf{L}(n)^T \mathbf{L}(n)\right\rVert_F \,  \left\lVert\bm{\theta}_2 - \bm{\theta}_1 \right\rVert_2^2 + \left\lVert\bm{\delta}^{(3)} \right\rVert_2^2 .
			\end{aligned}
		\end{equation*}

		We now consider each term of the above equation separately. First, since $\bm{a}(\btheta)$ is a contraction map, there exists a $\varepsilon\in (0,1)$ such that 
	    \begin{equation*}
	    \left\lVert\mathbf{a}(\btheta_2) - \mathbf{a}(\btheta_1)\right\rVert_2 \leq \varepsilon\left\lVert\btheta_2 - \btheta_1\right\rVert_2.   
	    \end{equation*}
		Secondly, we have
		\begin{equation*}
			\lVert\mathbf{L}(n)^T \mathbf{L}(n)\rVert_F \,  \left\lVert\bm{\theta}_2 - \bm{\theta}_1 \right\rVert_2^2 = \Delta \; \left\lVert\bm{\theta}_2 - \bm{\theta}_1 \right\rVert_2^2,
		\end{equation*}
		where $\Delta = \mathcal{O}\left(p^2 n^{-2\beta} \right)$. Indeed, by writing $\mathbf{B} = \mathbf{L}(n)^T \mathbf{L}(n)$ we obtain
		\begin{equation*}
			\begin{aligned}
				\Delta &= \lVert \mathbf{B} \rVert_F = \sqrt{\sum_{j = 1}^p \sum_{l = 1}^p B_{j,l}^2} \leq p \max_{j,l = 1,\, \ldots, \, p} |B_{j,l}| = p \max_{j,l = 1,\, \ldots, \, p} |\sum_{k = 1}^p \mathbf{L}_{k,j}(n) \mathbf{L}_{k,l}(n)|\\
	 & \leq p^2	\max_{j,l,k = 1,\, \ldots, \, p} |\mathbf{L}_{k,j}(n) \mathbf{L}_{k,l}(n)| = \mathcal{O}\left(p^2 n^{-2\beta}\right).
			\end{aligned}
	\end{equation*}
	Finally, we have
	\begin{equation*}
		\lVert\bm{\delta}^{(3)} \rVert_2^2 = \sum_{j = 1}^p \left(\bm{\delta}_j^{(3)}\right)^2 \leq p \max_{j = 1,\,\ldots,\,p} \left(\bm{\delta}_j^{(3)}\right)^2 = \mathcal{O}\left(p n^{-2 \gamma}\right) + \mathcal{O}_{\rm p}\left(p n^{-2 \alpha}\right).
	\end{equation*}
	%
	%
    By combining these results, we have
	\begin{equation*}
			\left\Vert T_{(j,n,H)}(\bm{\theta}_1) - T_{(j,n,H)}(\bm{\theta}_2) \right\Vert_2^2 \leq (\varepsilon^2+\Delta) \left\Vert\bm{\theta}_2 - \bm{\theta}_1 \right\Vert_2^2 + \mathcal{O}\left(p n^{-2 \gamma}\right) + \mathcal{O}_{\rm p}\left(p n^{-2 \alpha}\right).
	\end{equation*}
	Since $\alpha, \beta, \gamma> 0$ by Assumptions \ref{assum:C:1} and \ref{assum:D:1} and $\bm{\Theta}$ is bounded by Assumption \ref{assum:A:1}, then for sufficiently large $n$ we have that for all $\bm{\theta}_1, \, \bm{\theta}_2 \in \bm{\Theta}$ 
	\begin{equation}
	\label{T:contract}
		\left\Vert T_{(j,n,H)}(\bm{\theta}_1) - T_{(j,n,H)}(\bm{\theta}_2) \right\Vert_2 <  \left\Vert\bm{\theta}_2 - \bm{\theta}_1\right\Vert_2.
	\end{equation}
	Considering $T_{(j,n,H)}$ as a function from $\real^p$ to itself we can apply Banach fixed-point theorem implying that there exists a unique fixed point $\bm{\theta}^\ast_{(j,n,H)} \in \real^p$ such that
	\begin{equation*}
		T_{(j,n,H)} \left(\bm{\theta}^\ast_{(j,n,H)}\right) = \bm{\theta}^\ast_{(j,n,H)}\,.
	\end{equation*}
	By Assumption \ref{assum:A:1} we have that $\bm{\theta}^\ast_{(j,n,H)} \in \bm{\Theta}$, which implies by (\ref{eq:setfixed_point}) that
	\begin{equation*}
		\widehat{\bm{\Theta}}_{(j,n,H)} = \left\{ \bm{\theta}^\ast_{(j,n,H)} \right\} = \left\{ \hat{\bm{\theta}}_{(j,n,H)} \right\},
	\end{equation*}
	which concludes the proof. \hfill \qedsymbol
	
	\section{Proof of Proposition \ref{thm:iter:boot}}
	\label{proof:thm:1}
	
	We start by recalling that  $\left\{ \tilde{\bm{\theta}}^{(k)}_{(j,n,H)} \right\}_{k \in \mathbb{N}}$ is defined as
	\begin{equation*}
		\tilde{\bm{\theta}}^{(k)}_{(j,n,H)} = T_{(j,n,H)} \left(\tilde{\bm{\theta}}^{(k-1)}_{(j,n,H)}\right)
	\end{equation*}
	where $\tilde{\bm{\theta}}^{(0)}_{(j,n,H)} \in \bm{\Theta}$ and the function $T_{(j,n,H)}$ is defined in (\ref{eq:define:T}).
	Using the same arguments used in the proof of Lemma \ref{lemma:unique:iter:boot} we have that $T_{(j,n,H)}$ allows to apply Banach fixed-point theorem implying that 
	\begin{equation*}
		\lim_{k \to \infty}  \; \tilde{\bm{\theta}}^{(k)}_{(j,n,H)} = \hat{\bm{\theta}}_{(j,n,H)},
	\end{equation*}
	which concludes the first part of the proof.
	
	For the second part, consider a $k\in \mathbb{N}^\ast$. By the inequality (\ref{T:contract}) of the proof of Lemma \ref{lemma:unique:iter:boot}, there exists an $\epsilon \in (0, \, 1)$ such that 
	\begin{equation*}
\left\Vert T_{(j,n,H)}\left(\tilde{\bm{\theta}}^{(k)}_{(j,n,H)}\right) - T_{(j,n,H)}\left(\hat{\bm{\theta}}_{(j,n,H)}\right) \right\Vert_2 \leq \epsilon \left\Vert\tilde{\bm{\theta}}^{(k)}_{(j,n,H)} - \hat{\bm{\theta}}_{(j,n,H)}\right\Vert_2.
	\end{equation*}
	By construction of the sequence $\left\{ \tilde{\bm{\theta}}^{(k)}_{(j,n,H)} \right\}_{k \in \mathbb{N}}$ and since $\hat{\bm{\theta}}_{(j,n,H)}$ is a fixed point of $T_{(j,n,H)}$, we can derive that
	\begin{equation*}
	    \left\Vert T_{(j,n,H)}\left(\tilde{\bm{\theta}}^{(k)}_{(j,n,H)}\right) - T_{(j,n,H)}\left(\hat{\bm{\theta}}_{(j,n,H)}\right) \right\Vert_2 \leq \frac{\epsilon^k}{1-\epsilon} \left\Vert\tilde{\bm{\theta}}^{(0)}_{(j,n,H)} - \tilde{\bm{\theta}}^{(1)}_{(j,n,H)}\right\Vert_2  =\mathcal{O}_{\rm p}({p}^{\nicefrac{1}{2}}\epsilon^k),
	\end{equation*}
	which concludes the second part of the proof. \hfill \qedsymbol
	
	\section{Proof of Proposition \ref{THM:bias}}
	\label{proof:thm:bias}
	
	From (\ref{eq:indirectInf:hTimesN}) we have that $\hat{\bm{\theta}}_{(j,n,H)}$ is for $(j, n, H) \in \mathbb{N} \times \mathbb{N}^\ast \times \mathbb{N}^\ast$ such that:
	\begin{equation*}
	    \hat{\bm{\pi}}(\bm{\theta}_0, n, \bm{\omega}_0) = \frac{1}{H} \sum_{h = 1}^H  \hat{\bm{\pi}}(\hat{\bm{\theta}}_{(j,n,H)}, n, \bm{\omega}_{h+jH}).
	\end{equation*}
	Using Assumptions \ref{assum:B:1} and \ref{assum:D:2}, we expand each side of the above equation as follows:
	\begin{equation}
	    \begin{aligned}
	     \hat{\bm{\pi}}(\bm{\theta}_0, n, \bm{\omega}_0) &= \bm{\theta}_0 + \mathbf{a}(\bm{\theta}_0) + \mathbf{c}(n) + \mathbf{L}(n) \bm{\theta}_0 + \mathbf{r} (\bm{\theta}_0, n) + \mathbf{v} \left(\bm{\theta}_0, n, \bm{\omega}_0\right)\\
	     \frac{1}{H} \sum_{h = 1}^H  \hat{\bm{\pi}}(\hat{\bm{\theta}}_{(j,n,H)}, n, \bm{\omega}_{h+jH}) &= \hat{\bm{\theta}}_{(j,n,H)} + \mathbf{a}(\hat{\bm{\theta}}_{(j,n,H)}) +\mathbf{c}(n) + \mathbf{L}(n) \hat{\bm{\theta}}_{(j,n,H)} + \mathbf{r} (\hat{\bm{\theta}}_{(j,n,H)}, n) \\
	     &\quad +  \frac{1}{H} \sum_{h = 1}^H  \mathbf{v}(\hat{\bm{\theta}}_{(j,n,H)}, n, \bm{\omega}_{h+jH}).
	    \end{aligned}
	    \label{eq:proof:bias:inter:1}
	\end{equation}
	Therefore, we obtain
	\begin{equation*}
	    \begin{aligned}
	    &\mathbb{E}\left[\frac{1}{H} \sum_{h = 1}^H  \hat{\bm{\pi}}(\hat{\bm{\theta}}_{(j,n,H)}, n, \bm{\omega}_{h+jH}) - \hat{\bm{\pi}}(\bm{\theta}_0, n, \bm{\omega}_0) \right]\\
	    &= \left(\mathbf{I} + \mathbf{L}(n)\right)\mathbb{E}\left[\hat{\bm{\theta}}_{(j,n,H)} - \bm{\theta}_0\right] + \mathbb{E}\left[\mathbf{a}(\hat{\bm{\theta}}_{(j,n,H)}) - \mathbf{a}(\bm{\theta}_0)\right] + \mathbb{E} \left[\mathbf{r} (\hat{\bm{\theta}}_{(j,n,H)}, n) - \mathbf{r} (\bm{\theta}_0, n)\right]\\
	    &= \left(\mathbf{I} + \mathbf{L}(n) + \mathbf{M}\right) \mathbb{E}\left[\hat{\bm{\theta}}_{(j,n,H)} - \bm{\theta}_0\right] +\mathbb{E} \left[\mathbf{r} (\hat{\bm{\theta}}_{(j,n,H)}, n) - \mathbf{r} (\bm{\theta}_0, n)\right]\\
	    &= \mathbf{0}.
	    \end{aligned} 
	\end{equation*}
	By hypothesis, $n$ is such that the inverse of $\mathbf{B} \equiv \mathbf{I} + \mathbf{L}(n) + \mathbf{M}$ exists and we obtain
	\begin{equation}
	    \mathbb{E}\left[\hat{\bm{\theta}}_{(j,n,H)} - \bm{\theta}_0\right] = -\mathbf{B}^{-1} \mathbb{E} \left[\mathbf{r} (\hat{\bm{\theta}}_{(j,n,H)}, n) - \mathbf{r} (\bm{\theta}_0, n)\right]. 
	    \label{proof:bias:eq:inter2}
	\end{equation}
    By Assumptions \ref{assum:A:2} and \ref{assum:D:2}, $\mathbf{r} (\hat{\bm{\theta}}_{(j,n,H)}, n)$ is a bounded random variable on a compact set. Moreover, since $\mathbf{r}(\bm{\theta},n) = \mathcal{O}(n^{-\gamma})$ elementwise by Assumption \ref{assum:D:2}, we have $\mathbb{E} \left[\mathbf{r} (\hat{\bm{\theta}}_{(j,n,H)}, n) - \mathbf{r} (\bto,n) \right] = \mathcal{O}(n^{-\gamma})$ elementwise. Consequently, we deduce from \eqref{proof:bias:eq:inter2} that 
    \begin{equation}
        \mathbb{E} \left[ \hbt - \bto \right] = \mathcal{O}(pn^{-\gamma})
        \label{proof:bias:eq:inter3}
    \end{equation}
    elementwise. 
    
    The idea now is to re-evaluate $\mathbb{E} \left[\mathbf{r} (\hat{\bm{\theta}}_{(j,n,H)}, n) - \mathbf{r} (\bm{\theta}_0, n)\right]$ using the mean value theorem. This will allow us to make an induction that will show that for all $i\in\mathbb{N}$ 
	\begin{equation*}
	     \Big\lVert \mathbb{E}\left[ \hbt - \bto \right] \Big\rVert_2 = \mathcal{O}\left((p^2n^{-\gamma})^i\right),
	\end{equation*}
	and hence ends the proof.
    
    Applying the mean value theorem for vector-valued functions to $\mathbf{r} (\hat{\bm{\theta}}_{(j,n,H)}, n) - \mathbf{r} (\bm{\theta}_0, n)$ we have
	\begin{equation*}
	    \mathbf{r} (\hat{\bm{\theta}}_{(j,n,H)}, n)- \mathbf{r} (\bm{\theta}_0, n) = \mathbf{R}\left(\bt^{(\mathbf{r})},n\right) \left(\hat{\bm{\theta}}_{(j,n,H)} - \bm{\theta}_0 \right).
	\end{equation*}
	where the matrix $\mathbf{R}\left(\bt^{(\mathbf{r})},n\right)$ is defined in \eqref{def:MVT:multi_vari}. 
	By Assumptions \ref{assum:A:2} and \ref{assum:D:2}, $\mathbf{R}\left(\bt^{(\mathbf{r})},n\right)$ is also a bounded random variable. Moreover, we have $\mathbb{E}\left[\mathbf{R}\left(\bt^{(\mathbf{r})},n\right)\right] = \mathcal{O}(n^{-\gamma})$ elementwise since by Assumption \ref{assum:D:2} $\mathbf{r}(\bm{\theta},n) = \mathcal{O}(n^{-\gamma})$ elementwise, and 
	\begin{equation*}
	    \mathbf{r} (\bm{\theta}, n) = \mathbf{r} (\bm{\theta}_0, n) + \mathbf{R}\left(\bt^{(\mathbf{r})},n\right) \left(\bm{\theta} - \bm{\theta}_0 \right).
	\end{equation*}
	for all $\bm{\theta} \in \bm{\Theta}$.
	For simplicity, we denote 
	\begin{equation*}
	    \mathbf{R} \equiv \mathbf{R}\left(\bt^{(\mathbf{r})},n\right),\ \ \
	    \mathbf{\Delta}^{(\mathbf{r})} \equiv  \mathbf{r} (\hat{\bm{\theta}}_{(j,n,H)}, n) - \mathbf{r} (\bm{\theta}_0, n)
	    \ \ \ \text{and} \ \ \ \mathbf{\Delta} \equiv \hat{\bm{\theta}}_{(j,n,H)} - \bm{\theta}_0,
	\end{equation*} 
	which implies that $\mathbf{\Delta}^{(\mathbf{r})}=\mathbf{R}\mathbf{\Delta}$. Moreover, a consequence of \eqref{proof:bias:eq:inter3} is that 
		\begin{equation*}
	    \big|\mathbb{E}\left[\mathbf{\Delta}_k\right]\big|=\mathcal{O}(pn^{-\gamma}),
	\end{equation*}
	for any $k=1,\dots,p$ and hence 
	\begin{equation*}
	    \Vert \mathbb{E}\left[\mathbf{\Delta}\right] \Vert_2 = \mathcal{O}\left(p^2n^{-\gamma}\right).
	\end{equation*}
	Now, for all $l=1,\dots,p$ we have 
	\begin{equation}
	    \mathbf{\Delta}^{(\mathbf{r})}_l=\sum^{p}_{m=1}\mathbf{R}_{l,m}\mathbf{\Delta}_m \leq p \max_{m}\mathbf{R}_{l,m}\mathbf{\Delta}_m.
	    \label{proof:bias:eq:inter4}
	\end{equation}

  Using Cauchy-Schwarz inequality, we have 
   \begin{equation*}
       \mathbb{E} \left[ | \mathbf{R}_{l,m}\mathbf{\Delta}_m |\right] \leq  \mathbb{E} \left[ \mathbf{R}_{l,m}^2 \right]^{\nicefrac{1}{2}} \mathbb{E} \left[ \mathbf{\Delta}_m^2 \right]^{\nicefrac{1}{2}}.
    \end{equation*}
	Since $\mathbf{R}_{l,m}$ and $\mathbf{\Delta}_m$ are bounded random variables on a compact set, $\mathbb{E}\left[\mathbf{R}_{l,m}\right] = \mathcal{O}\left(n^{-\gamma}\right)$ and $\mathbb{E}\left[\mathbf{\Delta}_m\right] = \mathcal{O}\left(pn^{-\gamma}\right)$ we have 
	\begin{equation*}
	    \mathbb{E} \left[ \mathbf{R}_{l,m}^2 \right]^{\nicefrac{1}{2}} \mathbb{E} \left[ \mathbf{\Delta}_m^2 \right]^{\nicefrac{1}{2}} = \mathcal{O}\left(n^{-\gamma}\right)\mathcal{O}\left(p n^{-\gamma}\right)=\mathcal{O}\left(p n^{-2\gamma}\right).
	\end{equation*}
	Therefore, we obtain 
	\begin{equation*}
    \mathbb{E} \left[ | \mathbf{R}_{l,m}\mathbf{\Delta}_m |\right] \leq  \mathbb{E} \left[ \mathbf{R}_{l,m}^2 \right]^{\nicefrac{1}{2}} \mathbb{E} \left[ \mathbf{\Delta}_m^2 \right]^{\nicefrac{1}{2}} = \mathcal{O}\left(p n^{-2\gamma}\right),
    \end{equation*}
	and hence, using \eqref{proof:bias:eq:inter4} we deduce that
	\begin{equation*}
        \Big| \mathbb{E} \left[  \mathbf{\Delta}^{(\mathbf{r})}_l \right] \Big| \leq p \mathbb{E} \left[ \max_{m} | \mathbf{R}_{l,m}\mathbf{\Delta}_m | \right] = \mathcal{O}\left((pn^{-\gamma})^2\right). 
    \end{equation*}
	Considering this, and since $\mathbb{E}\left[\mathbf{\Delta}\right] = - \mathbf{B}^{-1}\mathbb{E}\left[\mathbf{\Delta}^{(\mathbf{r})}\right]$, we have, for all $k=1,\dots,p$
	\begin{equation*}
	   \big| \mathbb{E} \left[  \mathbf{\Delta}_k\right] \big| =  \mathcal{O}\left(p^3n^{-2\gamma}\right),
	\end{equation*}
	and consequently,
	\begin{equation*}
	   \lVert \mathbb{E} \left[  \mathbf{\Delta}\right] \rVert_2 =  \mathcal{O}\left((p^2n^{-\gamma})^2\right).
	\end{equation*}
	Since $\mathbb{E}\left[\mathbf{\Delta}\right] = -\mathbf{B}^{-1}\mathbb{E}\left[\mathbf{\Delta}^{(\mathbf{r})}\right]$ and $\mathbb{E}\left[\mathbf{\Delta}^{(\mathbf{r})}\right] = \mathbb{E}\left[\mathbf{R}\mathbf{\Delta}\right]$, one can repeat the same computations and deduce by induction that for all $i\in\mathbb{N}^*$
	\begin{equation*}
	    \Vert \mathbb{E}\left[\mathbf{\Delta}\right] \Vert_2 = \mathcal{O}\left((p^2n^{-\gamma})^i\right).
	\end{equation*}  \hfill \qedsymbol

    \section{Proof of Corollary \ref{coro:consist}}
    \label{proof:coro:consist}
    
    Fix $\varepsilon > 0$ and $\delta > 0$. We need to show that there exists a sample size $n^\ast \in\mathbb{N}^\ast$ such that for all $n \geq n^\ast$
	\begin{equation*}
	    \Pr \left(|| \hat{\bm{\theta}}_{(j,n,H)} - \bm{\theta}_0 ||_2 \geq \varepsilon \right) \leq   \delta.
	\end{equation*}
	From Chebyshev's inequality we have that
	\begin{equation*}
	    \Pr \left(|| \hat{\bm{\theta}}_{(j,n,H)} - \bm{\theta}_0 ||_2  \geq \varepsilon \right) \leq \frac{\mathbb{E}\left[|| \hat{\bm{\theta}}_{(j,n,H)} - \bm{\theta}_0 ||_2^2\right]}{\varepsilon^2}.
	\end{equation*}
	Therefore, we only need to show that there exists a sample size $n^\ast \in\mathbb{N}^\ast$ such that for all $n \geq n^\ast$
	\begin{equation*}
	    \mathbb{E}\left[|| \hat{\bm{\theta}}_{(j,n,H)} - \bm{\theta}_0 ||_2^2\right] \leq \varepsilon^2\delta .
	\end{equation*}

	Using the same notation as in the proof of Proposition \ref{THM:bias}, we have from (\ref{eq:proof:bias:inter:1}) and for $n$ sufficiently large that
	\begin{equation}
	\begin{aligned}
	   \mathbf{\Delta} & \equiv \hat{\bm{\theta}}_{(j,n,H)} - \bm{\theta}_0 \\
	    & = \mathbf{B}^{-1} \left(\mathbf{r} (\bm{\theta}_0, n) - \mathbf{r} (\hat{\bm{\theta}}_{(j,n,H)}, n) + \mathbf{v} \left(\bm{\theta}_0, n, \bm{\omega}_0\right) - \frac{1}{H} \sum_{h = 1}^H  \mathbf{v}(\hat{\bm{\theta}}_{(j,n,H)}, n, \bm{\omega}_{h+jH})\right)\\
	    & =  \mathbf{B}^{-1} \left( \mathbf{R}\mathbf{\Delta} + \mathbf{v}^{(H)} \right), 
	    \end{aligned}
	    \label{eq:comment:consist}
	\end{equation}
	where $\mathbf{v}^{(H)}$ is zero mean random vector of {\color{black}order $\mathcal{O}_{\rm p}\left(n^{-\alpha}\right)$} elementewise by Assumption \ref{assum:C:1}. Thus, by Assumption \ref{assum:D:2} we have that for all $k = 1,\dots,p$ 
	\begin{equation*}
	   \begin{aligned}
	        \mathbf{\Delta}_k &= \sum^p_{l=1} \left( \mathbf{B}^{-1} \right)_{k,l}\left(\mathbf{R}\mathbf{\Delta} + \mathbf{v}^{(H)}\right)_l = 
	        \sum^p_{l=1} \sum^p_{m=1} \left( \mathbf{B}^{-1} \right)_{k,l}\left(\mathbf{R}_{l,m}\mathbf{\Delta}_m + \mathbf{v}^{(H)}_m\right) \\
	         &= \sum^p_{l=1} \sum^p_{m=1} \mathcal{O_{\rm p}}(n^{-\gamma}) + \mathcal{O_{\rm p}}(n^{-\alpha}) =
	         \sum^p_{l=1} \sum^p_{m=1} \mathcal{O_{\rm p}}\left(\max(n^{-\alpha} , n^{-\gamma})\right) 
	         \\
	         &= \mathcal{O_{\rm p}}\left(p^2\max(n^{-\alpha} , n^{-\gamma})\right).
	   \end{aligned}
	\end{equation*}
	%
	%
	Therefore we have,
	%
	\begin{equation*}
	  \left\lVert \mathbf{\Delta} \right\rVert_2^2 = \sum_{k=1}^p \mathbf{\Delta}_k^2 = p\mathcal{O}_{\rm  p}\left(p^4\max(n^{-2\alpha}  , n^{-2\gamma})\right) = \mathcal{O}_{\rm  p}\left(p^5 n^{-2\min(\alpha, \gamma)}\right),
	\end{equation*}
	%
	and thus since $\bm{\Theta}$ is compact by Assumption \ref{assum:A:2}, we have
	\begin{equation*}
	\begin{aligned}
	    \mathbb{E}\left[\lVert \mathbf{\Delta} \rVert_2^2\right] = \mathcal{O}\left(p^5 n^{-2\min(\alpha, \gamma)}\right).
	    \end{aligned}
	\end{equation*}
	Using Assumptions \ref{assum:C:2} and \ref{assum:D:3} the last equality implies that there exists a sample size $n^\ast \in\mathbb{N}^\ast$ such that for all $n \in \mathbb{N}^*$ satisfying $n \geq n^\ast$, we obtain 
	\begin{equation*}
	\begin{aligned}
	    \mathbb{E}\left[\lVert \mathbf{\Delta} \rVert_2^2\right] \leq \varepsilon^2\delta.
	    \end{aligned}
	\end{equation*}

	\begin{flushright}
	\qedsymbol
	\end{flushright}

	\section{Proof of Proposition \ref{THM:consistency}}
	\label{proof:thm:2}

	This proof is directly obtained by verifying the conditions of Theorem 2.1 of \cite{newey1994large} on the functions $\hat{Q}(\bm{\theta}, n)$ and $Q(\bm{\theta})$ defined as follow:

	\begin{equation*}
	    \hat{Q}(\bm{\theta}, n) \equiv \Big\Vert \hat{\bm{\pi}}(\bm{\theta}_0, n, \bm{\omega}_0) - \frac{1}{H} \sum_{h = 1}^H  \hat{\bm{\pi}}(\bm{\theta}, n, \bm{\omega}_{h+jH}) \Big\Vert_2,
	\end{equation*}
	and 
	\begin{equation*}
	     {Q}(\bm{\theta}) = \big\Vert {\bm{\pi}}(\bm{\theta}_0) -  {\bm{\pi}}(\bm{\theta}) \big\Vert_2.
	\end{equation*}
	Reformulating the requirements of this theorem to our setting, we have to show that  (i) $\bm{\Theta}$ is compact, (ii) ${Q}(\bm{\theta})$ is continuous, (iii) ${Q}(\bm{\theta})$ is uniquely minimized at $\btheta_0$, (iv) $\hat{Q}(\bm{\theta}, n)$ converges uniformly in probability to $Q(\bm{\theta})$.

	On one hand, Assumptions \ref{assum:A:3} trivially ensures that $\bm{\Theta}$ is compact. On the other hand, Assumption \ref{assum:D:4} guarantees that ${Q}(\bm{\theta})$ is continuous and uniquely minimized at $\btheta_0$ since $\pi(\btheta)=\btheta+\mathbf{a}(\btheta)$ is required to be continuous and injective. What remains to be shown is that $\hat{Q}(\bm{\theta}, n)$ converges uniformly in probability to $Q(\bm{\theta})$, which is equivalent to show that: $\forall \varepsilon > 0$ and $\forall \delta > 0$, there exists a sample size $n^\ast \in\mathbb{N}^\ast$ such that for all $n \geq n^\ast$
	\begin{equation*}
	    \Pr \left( \sup_{\bm{\theta} \in \bm{\Theta}} \; \Big| \hat{Q}(\bm{\theta}, n) - Q(\bm{\theta}) \Big| \geq \varepsilon \right) \leq   \delta.
	\end{equation*}

	Fix $\varepsilon > 0$ and $\delta > 0$. Using the above definitions, we have that
	\begin{equation}
	    \sup_{\bm{\theta} \in \bm{\Theta}} \; \Big| \hat{Q}(\bm{\theta}, n) - Q(\bm{\theta}) \Big| \leq \sup_{\bm{\theta} \in \bm{\Theta}} \; \left[ \Big| \hat{Q}(\bm{\theta}, n) - Q(\bm{\theta}, n) \Big| + \Big| {Q}(\bm{\theta}, n) - Q(\bm{\theta}) \Big| \right],
	    \label{eq:convergen_proof_eq}
	\end{equation}
	where
	\begin{equation}
	        {Q}(\bm{\theta}, n) \equiv \big\Vert {\bm{\pi}}(\bm{\theta}_0, n) -  {\bm{\pi}}(\bm{\theta}, n) \big\Vert_2.
	\end{equation}
	We now consider each term of (\ref{eq:convergen_proof_eq}) separately. For simplicity, we define
		\begin{equation*}
	    \bar{\bm{\pi}}(\bm{\theta}, n) \equiv \frac{1}{H} \sum_{h = 1}^H  \hat{\bm{\pi}}(\bm{\theta}, n, \bm{\omega}_{h+jH}),
	\end{equation*}
	and, considering the first term of \eqref{eq:convergen_proof_eq}, we have
	\begin{equation*}
	   \begin{aligned}
	       \Big| \hat{Q}(\bm{\theta}, n) - Q(\bm{\theta}, n) \Big| & = \Big| \big\Vert \hat{\bm{\pi}}(\bm{\theta}_0, n, \bm{\omega}_0) - {\bm{\pi}}(\bm{\theta}_0, n) + {\bm{\pi}}(\bm{\theta}_0, n) - {\bm{\pi}}(\bm{\theta}, n)  + {\bm{\pi}}(\bm{\theta}, n) - \bar{\bm{\pi}}(\bm{\theta}, n) \big\Vert_2 \\
	       &- \big\Vert {\bm{\pi}}(\bm{\theta}_0, n) -  {\bm{\pi}}(\bm{\theta}, n) \big\Vert_2 \Big|\\
	       &\leq  \big\Vert \hat{\bm{\pi}}(\bm{\theta}_0, n, \bm{\omega}_0) - {\bm{\pi}}(\bm{\theta}_0, n) \big\Vert_2 + \big\Vert {\bm{\pi}}(\bm{\theta}, n) - \bar{\bm{\pi}}(\bm{\theta}, n) \big\Vert_2 \\
	       &= \big\Vert \mathbf{v} \left(\bm{\theta}_0, n, \bm{\omega}_0\right) \big\Vert_2 + \Big\Vert \frac{1}{H}\sum_{h=1}^H\mathbf{v} \left(\bm{\theta}, n, \bm{\omega}_{h+jH}\right) \Big\Vert_2 \\
	       &= \mathcal{O}_{\rm p} \left(\sqrt{p} n^{-\alpha} \right) + \mathcal{O}_{\rm p} \left(\sqrt{p} n^{-\alpha} H^{-\nicefrac{1}{2}}\right) = \mathcal{O}_{\rm p} \left(\sqrt{p} n^{-\alpha} \right),
	   \end{aligned} 
	\end{equation*}
	by Assumption \ref{assum:C:1}. Similarly, we have
	\begin{equation*}
	   \begin{aligned}
	    \Big| {Q}(\bm{\theta}, n) - Q(\bm{\theta}) \Big| &= \Big| \big\Vert {\bm{\pi}}(\bm{\theta}_0, n) - {\bm{\pi}}(\bm{\theta}_0) + {\bm{\pi}}(\bm{\theta}_0) - {\bm{\pi}}(\bm{\theta})  + {\bm{\pi}}(\bm{\theta}) - {\bm{\pi}}(\bm{\theta}, n) \big\Vert_2 \\
	       &- \big\Vert {\bm{\pi}}(\bm{\theta}_0) -  {\bm{\pi}}(\bm{\theta}) \big\Vert_2 \Big| \\
	       &\leq  \big\Vert {\bm{\pi}}(\bm{\theta}_0, n) - {\bm{\pi}}(\bm{\theta}_0) \big\Vert_2 + \big\Vert {\bm{\pi}}(\bm{\theta}) - {\bm{\pi}}(\bm{\theta}, n) \big\Vert_2 \\
	       &=   \big\Vert \mathbf{c}(n) + \mathbf{L}(n)\bto + \mathbf{r}(\bto,n) \big\Vert_2 + \big\Vert \mathbf{c}(n) + \mathbf{L}(n)\bm{\theta} + \mathbf{r}(\bm{\theta},n) \big\Vert_2 \\
	       &= \mathcal{O} \left(\sqrt{p} \max \left( c_n, p n^{-\beta}, n^{-\gamma} \right)\right),
	   \end{aligned} 
	\end{equation*}
	by Assumption \ref{assum:D:4}. Therefore, we obtain
    \begin{equation*}
        \sup_{\bm{\theta} \in \bm{\Theta}} \; \Big| \hat{Q}(\bm{\theta}, n) - Q(\bm{\theta}) \Big| = \mathcal{O}_{\rm p} \left(\sqrt{p} n^{-\alpha} \right) + \mathcal{O} \left(\sqrt{p} \max \left(c_n, p n^{-\beta}, n^{-\gamma} \right)\right).
    \end{equation*}
	By Assumptions \ref{assum:C:1} and \ref{assum:D:4}, there exists a sample size $n^\ast \in\mathbb{N}^\ast$ such that for all $n \in \mathbb{N}^*$ satisfying $n \geq n^\ast$ we have
	\begin{equation*}
	    \Pr \left( \sup_{\bm{\theta} \in \bm{\Theta}} \; \Big| \hat{Q}(\bm{\theta}, n) - Q(\bm{\theta}) \Big| \geq \varepsilon \right) \leq   \delta.
	\end{equation*}
	Therefore, the four condition of Theorem 2.1 of \cite{newey1994large} are verified implying the result.
	\begin{flushright}
	\qedsymbol
	\end{flushright}

\section{Proof of Proposition \ref{Thm:Gauss:approx}}
\label{proof:Gauss:approx}
Without loss of generality we can assume that the sample sizes $n^* \in \mathbb{N}^*$ of Assumption \ref{assum:C:3} and \ref{assum:D:5} are the same. Let $n \in \mathbb{N}^*$ be such that $n \geq n^*$.
By the definition of JIE in~\eqref{eq:indirectInf:hTimesN}, we have 
\begin{equation}\label{def:JIE}
     \hp{\bto}{\bwo} = \hpH{\hbt}.
\end{equation}
Since the auxiliary estimator may be expressed, for any $(\btheta, n, \bm{\omega}) \in \bm{\Theta} \times \mathbb{N}^* \times \bm{\Omega}$, as
\begin{equation}\label{eq:hat:pi}
    \hp{\bt}{\bw} = \bt + \a{\bt} + \cn + \Ln\bt + \r{\bt} + \v{\bt}{\bw}.
\end{equation}
we deduce from \eqref{def:JIE}, by using the mean value theorem on $\mathbf{a}(\btheta)$ and $\mathbf{r}(\btheta,n)$, that 
\begin{equation}
    \mathbf{B}\left(\btheta^{(\mathbf{a},\mathbf{r})},n\right)\bD= \v{\bto}{\bwo} - \vH{\hbt},
     \label{proof:prop4:B_THETA}
\end{equation}
where $\mathbf{B}\left(\btheta^{(\mathbf{a},\mathbf{r})},n\right) \equiv \mathbf{I} + \mathbf{A}(\btheta^{(\mathbf{a})}) + \mathbf{L}(n) + \mathbf{R}\left(\btheta^{(\mathbf{r})},n\right)$ and  $\bD \equiv \hbt - \bto$. Applying the mean value theorem on $\v{\hbt}{\bw_{h+jH}}$ we obtain for all $h=1,\dots,H$
\begin{equation*}
    \v{\hbt}{\bw_{h+jH}} = \v{\bto}{\bw_{h+jH}} + \mathbf{V}\left(\btheta^{(\mathbf{v})},n,\bw_{h+jH}\right)\bD.
\end{equation*}
Hence, let $\mathbf{u} \in \real^p$ be such that $|| \mathbf{u} ||_2 = 1$, then by Assumption \ref{assum:C:3} we have, 
\begin{equation*}
    \sqrt{n} \mathbf{u}^T \bm{\Sigma}(\bto,n)^{-\nicefrac{1}{2}}\v{\bto}{\bwo}   = Z(\bto, n, \bwo, \mathbf{u}) + \delta_n^{(1)} W(\bto, n, \bwo, \mathbf{u}) 
\end{equation*}
and for all $h=1,\dots,H$
\begin{equation*}
    \sqrt{n} \mathbf{u}^T \bm{\Sigma}(\bto,n)^{-\nicefrac{1}{2}}\v{\bto}{\bw_{h+jH}}   = Z(\bto, n, \bw_{h+jH}, \mathbf{u}) + \delta_n^{(h)} W(\bto, n, \bw_{h+jH}, \mathbf{u}) 
\end{equation*}
where $(\delta_n^{(0)}),(\delta_n^{(1)}),\dots, (\delta_n^{(H)}) \in \mathfrak{D}$ and 
\begin{equation*}
\begin{aligned}
        Z(\bto, n, \bwo, \mathbf{u})\sim \mathcal{N}(0,1), \;\;\; \;\;\;  Z(\bto, n, \bw_{h+jH}, \mathbf{u}) \sim \mathcal{N}(0,1), \\
        W(\bto, n, \bwo, \mathbf{u}) = \mathcal{O}_{\rm p}(1), \;\;\; \;\;\;  W(\bto, n, \bw_{h+jH}, \mathbf{u}) = \mathcal{O}_{\rm p}(1).
\end{aligned}
\end{equation*}
Without loss of generality, we can suppose that $(\delta_n) \equiv (\delta_n^{(0)})=(\delta_n^{(1)})=\dots=(\delta_n^{(H)})$ as one may simply modify $W(\bto, n, \bwo, \mathbf{u})$ and $W(\bto, n, \bw_{h+jH}, \mathbf{u})$ for all $h=1,\dots,H$. For simplicity, we write $\mathfrak{X}_{\bto,n,\mathbf{u}} \equiv \sqrt{n} \mathbf{u}^T \bm{\Sigma}(\bto,n)^{-\nicefrac{1}{2}}$, $Z_0 \equiv Z(\bto, n, \bwo, \mathbf{u})$ and $Z_h \equiv Z(\bto, n, \bw_{h+jH}, \mathbf{u})$.
Multiplying the right hand side of \eqref{proof:prop4:B_THETA} by $\mathfrak{X}_{\bto,n,\mathbf{u}}$, we have 
\begin{equation}
    \begin{aligned}
        &\mathfrak{X}_{\bto,n,\mathbf{u}}\left(\v{\bt}{\bwo} - \vH{\hbt}\right)\\ &=   
                Z_0 + \mathcal{O}_{\rm p}(\delta_n) -\dfrac{1}{H}\sum_{h=1}^{H} \big(Z_h + \mathcal{O}_{\rm p}(\delta_n)\big)  
         - \dfrac{1}{H}\sum_{h=1}^H \mathfrak{X}_{\bto,n,\mathbf{u}} \mathbf{V}\left(\btheta^{(\mathbf{v})},n,\bw_{h+jH}\right)\bD
    \end{aligned}
    \label{proof:prop4:XV_THETA}
\end{equation}
%

The first two terms of \eqref{proof:prop4:XV_THETA} are considered separately in the following computation 
\begin{equation}
    \begin{aligned}
        Z_0 + \mathcal{O}_{\rm p}(\delta_n) -\dfrac{1}{H}\sum_{h=1}^{H} \big(Z_h + \mathcal{O}_{\rm p}(\delta_n)\big) &=  Z_0 - \dfrac{1}{H}\sum_{h=1}^H Z_h + \mathcal{O}_{\rm p}(\delta_n) - \dfrac{1}{H}\sum_{h=1}^H \mathcal{O}_{\rm p}(\delta_n) 
        \\
        & \overset{d}{=} \left(1+\frac{1}{H}\right)^{\nicefrac{1}{2}} Z + \mathcal{O}_{\rm p}(\delta_n) - \mathcal{O}_{\rm p}\left(\frac{\delta_n}{\sqrt{H}}\right) \\
        &= \left(1+\frac{1}{H}\right)^{\nicefrac{1}{2}} Z + \mathcal{O}_{\rm p}(\delta_n),
    \end{aligned}
    \label{proof:prop4:XVI_THETA}
\end{equation}
where $Z \sim \mathcal{N}\left(0, 1\right)$. We now consider the last term of \eqref{proof:prop4:XV_THETA}. Since $\bm{\Sigma}(\btheta,n)^{-\nicefrac{1}{2}} = \mathcal{O}(1)$ and $\mathbf{V}\left(\btheta^{(\mathbf{v})},n,\bw_{h+jH}\right) = \mathcal{O}_p\left(n^{-\nicefrac{1}{2}}\right)$ elementwise by Assumption \ref{assum:C:3} and $|| \mathbf{u} ||_2 = 1$ we have for all $h=1,\dots,H$
\begin{equation*}
            \mathfrak{X}_{\bto,n,\mathbf{u}} \mathbf{V}\left(\btheta^{(\mathbf{v})},n,\bw_{h+jH}\right)\bD = \mathcal{O}_{\rm p}(p)\mathds{1}_p^T\bD,
\end{equation*}
where $\mathds{1}_p \in \real^p$ is the vector filled with ones. Since Assumptions \ref{assum:A:2}, \ref{assum:B:1}, \ref{assum:C:3} and \ref{assum:D:5} imply that $\hbt$ is consistent estimator of $\bto$ (see Figure \ref{Fig_assumptions-mess}), there exists a sequence $(\delta'_n) \in \mathfrak{D}$ such that $\bD=\mathcal{O}_{\rm p}(\delta'_n)$ for all $n \in \mathbb{N}^*$. Without loss of generality, we may assume that $(\delta'_n)=(\delta_n)$. Therefore, we have 
\begin{equation*}
\begin{aligned}
\mathfrak{X}_{\bto,n,\mathbf{u}} \mathbf{V}\left(\btheta^{(\mathbf{v})},n,\bw_{h+jH}\right)\bD 
& = 
\mathcal{O}_{\rm p}(p)\mathds{1}_p^T\bD =  
\mathcal{O}_{\rm p}(p)\mathds{1}_p^T\mathcal{O}_{\rm p}(\delta_n)\mathds{1}_p \\ 
& =
\mathcal{O}_{\rm p}(p\delta_n)\mathds{1}_p^T\mathds{1}_p
=
\mathcal{O}_{\rm p}(p^2\delta_n),
\end{aligned}
\end{equation*}
which leads to the computation for the last term of \eqref{proof:prop4:XV_THETA},
\begin{equation}
    \dfrac{1}{H}\sum_{h=1}^H \mathfrak{X}_{\bto,n,\mathbf{u}} \mathbf{V}\left(\btheta^{(\mathbf{v})},n,\bw_{h+jH}\right)\bD =
    \dfrac{1}{H}\sum_{h=1}^H \mathcal{O}_{\rm p}(p^2\delta_n) = \mathcal{O}_{\rm p}\left(\delta_n\frac{p^2}{\sqrt{H}}\right).
    \label{proof:prop4:XV_THETA_FINAL}
\end{equation}

Combining \eqref{proof:prop4:B_THETA}, \eqref{proof:prop4:XV_THETA} and \eqref{proof:prop4:XV_THETA_FINAL}, we obtain
\begin{equation*}
\begin{aligned}
   &\left(1+\frac{1}{H}\right)^{-\nicefrac{1}{2}}  \mathfrak{X}_{\bto,n,\mathbf{u}}\mathbf{B}\left(\btheta^{(\mathbf{a},\mathbf{r})},n\right)\bD   
   \\=
   &\left(1+\frac{1}{H}\right)^{-\nicefrac{1}{2}}  \mathfrak{X}_{\bto,n,\mathbf{u}}\left(\v{\bt}{\bwo} - \vH{\hbt}\right) 
    \\
    &\overset{d}{=}
    \left(1+\frac{1}{H}\right)^{-\nicefrac{1}{2}} \left(\left(1+\frac{1}{H}\right)^{\nicefrac{1}{2}} Z + \mathcal{O}_{\rm p}(\delta_n) + \mathcal{O}_{\rm p}\left(\delta_n\frac{p^2}{\sqrt{H}}\right)\right) 
    \\
    &\overset{d}{=}
    Z + \mathcal{O}_{\rm p}\left(\delta_n\sqrt{\frac{H}{H+1}}\right) + \mathcal{O}_{\rm p}\left(\delta_n\frac{p^2}{\sqrt{H+1}}\right)
    \\
    &=
    Z + \delta_n\mathcal{O}_{\rm p}\left(1\right) + \delta_n\mathcal{O}_{\rm p}\left(\frac{p^2}{\sqrt{H}}\right)
    =
    Z + \delta_n\mathcal{O}_{\rm p}\left(\max \left(1,\frac{p^2}{\sqrt{H}}\right)\right).
\end{aligned}
\end{equation*}
The proof is completed using \eqref{eq:assum:d:5} of Assumption \ref{assum:D:5}.

\begin{flushright}
\qedsymbol
\end{flushright}

\newpage
\bibliographystyle{plainnat}
\bibliography{biblio.bib}

\end{document}